\documentclass[review]{elsarticle}
\usepackage{lineno,hyperref}
\modulolinenumbers[5]

 \RequirePackage{geometry}
\def\norm#1{\|#1\|} 

\usepackage{algorithm}
\usepackage{algorithmic}
\usepackage{amsmath}
\usepackage{amsthm}
\usepackage{amssymb}
\usepackage{caption}
\usepackage{comment}
\usepackage{graphicx}
\usepackage{float}
\usepackage[latin1]{inputenc}
\usepackage{mathrsfs}
\usepackage{multicol}
\usepackage{multirow}
\usepackage{natbib}
\setcitestyle{numbers}
\usepackage{subcaption}
\usepackage{wrapfig}
\usepackage[dvipsnames]{xcolor}
\usepackage{lineno,hyperref}
\modulolinenumbers[5]
\usepackage[textsize=footnotesize]{todonotes}
\usepackage{siunitx} 

\usepackage{tikz}
\usetikzlibrary{shapes}
\usepackage{listings}

\def\genbox#1#2#3#4#5#6{
    \leavevmode\raise#4bp\hbox to#5bp{\vrule height#5bp depth0bp width0bp
    \pdfliteral{q .5 w \csname #2COLOR\endcsname\space RG
                       \csname #3PDF\endcsname{#5}{#6} S Q
             \ifx1#1 q \csname #2COLOR\endcsname\space rg 
                       \csname #3PDF\endcsname{#5}{#6} f Q\fi}\hss}}

\def\sqbox      #1#2{\genbox{#1}{#2}  {sq}       {0}   {4.5}  {2.25}}

\graphicspath{{./figures/}}

\newcommand{\beq} {\begin{equation}}
\newcommand{\eeq} {\end{equation}}
\newcommand{\bdm} {\begin{displaymath}}
\newcommand{\edm} {\end{displaymath}}
\newcommand{\bit}{\begin{itemize}}
\newcommand{\eit}{\end{itemize}}
\newcommand{\bde}{\begin{description}}
\newcommand{\ede}{\end{description}}
\newcommand{\bce}{\begin{center}}
\newcommand{\ece}{\end{center}}
\newcommand{\ben} {\begin{enumerate}}
\newcommand{\een} {\end{enumerate}}
\newcommand{\bea} {\begin{eqnarray}}
\newcommand{\eea} {\end{eqnarray}}
\newcommand{\barr} {\begin{array}}
\newcommand{\earr} {\end{array}}
\newcommand{\bean} {\begin{eqnarray*}}
\newcommand{\eean} {\end{eqnarray*}}
\newcommand{\edoc} {\end{document}}

\newcommand{\dd}[2] {\ensuremath{\frac{\partial {#1}}{\partial {#2}}}}
\newcommand{\DD}[2] {\ensuremath{\frac{d {#1}}{d {#2}}}}
\newcommand{\Grad} {\ensuremath{\nabla}}

\newcommand{\Div} {\ensuremath{\nabla\cdot}}

\newcommand{\eval}[2][\right]{\relax
  \ifx#1\right\relax \left.\fi#2#1\rvert}

\providecommand{\norm}[1]{\left\lVert#1\right\rVert}

\newcommand{\half} {\ensuremath{\frac{1}{2}}}

\newcommand{\mc}[1]{\mathcal{#1}}

\newcommand{\mb}[1]{\mathbf{#1}}

\newcommand{\nort}[1]{{\left\vert\kern-0.25ex\left\vert\kern-0.25ex\left\vert #1 
    \right\vert\kern-0.25ex\right\vert\kern-0.25ex\right\vert}}
\newcommand{\snor}[1]{\left| #1 \right|}
\newcommand{\LRp}[1]{\left( #1 \right)}
\newcommand{\LRs}[1]{\left[ #1 \right]}

\newcommand{\LRc}[1]{\left\{ #1 \right\}}

\newcounter{tablecell}

\newcommand{\figlab}[1]{\label{fig:#1}}
\newcommand{\eqnlab}[1]{\label{eq:#1}}

\newcommand{\tablab}[1]{\label{tab:#1}}

\newcommand{\figref}[1]{\ref{fig:#1}}

\newcommand{\eqnref}[1]{\eqref{eq:#1}}
\newcommand{\alglab}[1]{\label{alg:#1}}
\newcommand{\algref}[1]{\ref{alg:#1}}
\newcommand{\seclab}[1]{\label{sect:#1}}
\newcommand{\secref}[1]{\ref{sect:#1}}
\newcommand{\tabref}[1]{\ref{tab:#1}}

\renewcommand{\u}{u}

\newcommand{\q}{q}

\newcommand{\qb}{{\mb{\q}}}

\newcommand{\ub}{\mb{\u}}

\newcommand{\K}{K}

\newcommand{\pK}{{\partial \K}}

\renewcommand{\sb}{\mb{s}}

\newcommand{\Hb}{\mb{H}}

\newcommand{\Hbs}{{\mb{H}^*}}

\newcommand{\Gb}{\mb{G}}

\newcommand{\tb}{{\bf t}}
\newcommand{\nproc}{{N_p}}

\newcommand{\model}{{ \mathfrak{m} }}

\newcommand{\Ne}[1]{{N_{E_{#1}}}}

\newcommand{\Nsz}{{N_E^{\LRc{S}} }}
\newcommand{\Nfz}{{N_E^{\LRc{F}} }}

\newcommand{\Nbz}{{N_E^{\LRc{B}} }}
\newcommand{\Ntotal}{{N_E^{\texttt{total}} }}

\newcommand{\Nbeps}{{N_{\texttt{EPS}}^{\LRc{B}}}}
\newcommand{\Nszeps}{{N_{\texttt{EPS}}^{\LRc{S}}}}
\newcommand{\Nsps}{{N_{\texttt{EPS}}^{\LRc{S+B}} }}
\newcommand{\Nfps}{{N_{\texttt{EPS}}^{\LRc{F}} }}
\newcommand{\Ntotaleps}{{N_{\texttt{EPS}}^{\texttt{total}} }}

\newcommand{\Nze}[1]{{ N_{ze_{#1}} }}
\newcommand{\Nxe}[1]{{ N_{xe_{#1}} }}
\newcommand{\Nye}[1]{{ N_{ye_{#1}} }}

\newcommand{\jmhalf}{{j-\half}}
\newcommand{\jphalf}{{j+\half}}
\newcommand{\kmhalf}{{k-\half}}
\newcommand{\kphalf}{{k+\half}}

\newcommand{\imhalf}{{i-\half}}
\newcommand{\iphalf}{{i+\half}}

\newcommand{\dz}{{\triangle z}}

\newcommand{\xb}{{\bf x}}

\newcommand{\Fb}{{\bf F}}

\newcommand{\Fbs}{{{\bf F}^*}}

\newcommand{\Gbs}{{{\bf G}^*}}

\newcommand{\nb}{{\bf n}}

\newcommand{\Fcal}{\mathcal{F}}

\newcommand{\dtt}{\triangle t}

\newcommand{\Qbz}{{\bf Q}^{\LRc{z}} }
\newcommand{\qbz}{{\bf q}^{\LRc{z}} }

\newcommand{\Qbfz}{{\bf Q}^{\LRc{F}} }

\newcommand{\Qbsb}{{\bf Q}^{\LRc{B}} }
\newcommand{\Qbsz}{{\bf Q}^{\LRc{S}} }
\newcommand{\qbfz}{{\bf q}^{\LRc{F}} }

\newcommand{\qbsb}{{\bf q}^{\LRc{B}} }
\newcommand{\qbsz}{{\bf q}^{\LRc{S}} }
\newcommand{\Rbz}{{\bf R}^{\LRc{z}} }
\newcommand{\Rbfz}{{\bf R}^{\LRc{F}} }

\newcommand{\Rbsb}{{\bf R}^{\LRc{B}} }
\newcommand{\Rbsz}{{\bf R}^{\LRc{S}} }
\newcommand{\az}{{a}^{\LRc{z}} }
\newcommand{\bz}{{b}^{\LRc{z}} }
\newcommand{\cz}{{c}^{\LRc{z}} }

\newcommand{\dt}{{\triangle t}}
\newcommand{\dx}{{\triangle x}}
\newcommand{\dy}{{\triangle y}}

\newcommand{\pres}{{{p}}}

\newcommand{\Ical}{\mathcal{I}}

\newtheorem{remark}{Remark}

\begin{document}

\begin{frontmatter}

\title{Multirate Partitioned Runge--Kutta Methods for Coupled Navier--Stokes Equations}
\author[Argonne]{Shinhoo Kang\corref{mycorrespondingauthor}}
\cortext[mycorrespondingauthor]{Corresponding author}
\ead{shinhoo.kang@anl.gov}
\author[Argonne]{Alp Dener\fnref{alp}}
\fntext[alp]{This work was done by the author while at Argonne National Laboratory. The author's current affiliation is Siemens Digital Industries Software (Plano, TX, CA, USA). }
\ead{alp.dener@siemens.com}
\author[UDel]{Aidan Hamilton}
\ead{aidan@udel.edu}
\author[Argonne]{Hong Zhang}
\ead{hongzhang@anl.gov}
\author[Argonne]{Emil M. Constantinescu}
\ead{emconsta@anl.gov}
\author[Argonne-ESD]{Robert L. Jacob}
\ead{jacob@anl.gov}

\address[Argonne]{Mathematics and Computer Science Division, Argonne National Laboratory, Lemont, IL, USA}
\address[Argonne-ESD]{Environmental Science Division, Argonne National Laboratory, Lemont, IL, USA}
\address[UDel]{University of Delaware, Department of Mathematical Sciences, Newark, DE, USA}

\begin{abstract}

 Earth system models are complex integrated models of atmosphere, ocean, sea ice, and land surface. Coupling the components can be a significant challenge due to the difference in physics, temporal, and spatial scales. 
  This study explores multirate partitioned Runge--Kutta methods for the fluid-fluid interaction problem and demonstrates its parallel performance by using the PETSc library. 
  We consider compressible Navier--Stokes equations with gravity coupled through a rigid-lid interface. 
  Our large-scale numerical experiments reveal that multirate partitioned Runge--Kutta coupling schemes (1) can conserve total mass; 
  (2) have second-order accuracy in time; and 
  (3) provide favorable strong- and weak-scaling performance on modern computing architectures. 
  We also show that the speedup factors of multirate partitioned Runge--Kutta methods match theoretical expectations over their base (single-rate) method.
  
\end{abstract}


\end{frontmatter}


\section{Introduction}
Earth is an integrated system that consists of atmosphere, ocean, river,
land, and sea ice.  
Each submodel has different conservation laws, computational grid, and time step. Proper coupling between components is critical to maintain accuracy, mass conservation, and computational efficiency  \cite{jacob2005m,craig2012new,golaz2019doe}. 
In a previous study we proposed an implicit-explicit (IMEX) coupling method for coupled compressible Navier--Stokes systems \cite{KANG2021113988}.
 IMEX coupling methods can be expensive, however, because of the need to solve a linear or nonlinear system. 
This means that the computational cost highly depends on good preconditioning techniques.
For hyperbolic problems, 
developing such preconditioning methods is not a trivial task. Explicit multirate methods may therefore be a promising alternative.
A multirate method transforms an original problem into several subproblems while allowing different step sizes on each subproblem\cite{constantinescu2007multirate,sandu2019class,roberts2020coupled,hachtel2021multirate,gunther2021multirate,abdulle2022explicit}; on the contrary, a single-rate method is subject to severe step size limit due to the Courant--Friedrichs--Lewy condition if the subproblems have vastly different time scales. The computational advantages of multirate methods over single-rate methods have been demonstrated in many single-domain applications, including 
atmospheric~\cite{SkamarockKlemp08,seny2013multirate} and 
air pollution models~\cite{schlegel2012implementation},
hyperbolic problems \cite{lohner1984use,kirby2003convergence,constantinescu2007multirate},
Euler equations \cite{wensch2009multirate}, compressible Navier--Stokes equations \cite{mikida2019multi}, and adaptive mesh refined grids \cite{seny2013multirate}.
The work in \cite{seny2014efficient} presented
a strategy to parallelize explicit multirate schemes in the framework of discontinuous Galerkin methods for a single domain for shallow water and Euler equations. 

 Developing multirate methods for large-scale coupled climate models is challenging, however, not only because of the requirement in conservation and convergence but also because of the demand for scalability in a parallel framework. These coupled models consist of multiple domains and different sets of partial differential equations, and hence the interface treatment poses great numerical and computational difficulties. Cross-domain communication, transition between the coupling components, and load balancing need to be carefully considered.
 
In this study we apply multirate partitioned Runge--Kutta (MPRK) coupling methods \cite{constantinescu2007multirate} for fluid-fluid interaction problems, and demonstrate its strong and weak scalability on massively parallel computers by using the PETSc library. A buffer region at the interface is used for a coherent transition, which is key to achieving conservation and convergence. Domain decomposition is applied to each coupling problem as if the buffer region does not exist. Our strategy induces minimal additional complexity to an implementation that uses a single-rate method, while enjoying the full computational benefit of a multirate method. To evaluate the computational performance systematically, we develop a theoretical performance model for both serial and parallel cases. Our model takes the buffer region treatment into account, thus revealing its impact on the overall performance and shedding light on parallelization settings in practice. We describe several numerical 
experiments on coupled compressible Navier--Stokes equations (CNS) that demonstrate the suitability of using multirate methods for the coupling problem and the scalability of our implementation.

This paper is organized as follows.
In Section \secref{model-problems}
we describe the coupled systems and their spatial discretization. 
In Section \secref{MultirateCoupling} 
we explain multirate coupling methods.
In Section \secref{MPRKinPETSc} we discuss the speedup factors of the MPRK method over its base RK method. 
In Section \secref{NumericalResults} we demonstrate the performance of multirate coupling schemes through numerical examples. 
In Section \secref{Conclusion} we summarize our work.

\section{Model problems}
\seclab{model-problems}

The physical problem represents a simplified 
atmosphere-ocean system. The atmosphere and ocean have different thermodynamic properties because one
is a gas and the other is a liquid. For coupling under the rigid-lid assumption, the jumps in the velocity and temperature fields at the interface are needed for specifying the
bulk form \cite{KANG2021113988}. To simplify the coupling problem, we consider two fluids as
ideal gases but jumps of temperature and velocity across
the interface. 
Both the ocean and the atmosphere are described 
by compressible Navier--Stokes equations in the cell-centered second-order finite-volume (FV) spatial discretization on uniform grids with a rigid-lid coupling condition and gravity. 
To focus on the temporal challenges associated with the coupling, 
we use a conformal mesh at the interface. 
\footnote{
In general, the mesh used for ocean models is different from the mesh used for atmospheric model. To handle the non-matching grid at the interface, we can utilize the nonconforming mortar approach \cite{kopriva2002computation,bui2012analysis,kopriva2019free} or remapping algorithms \cite{trask2020compatible, mahadevan2022metrics,gmd-16-1537-2023}.
}


In Figure \figref{coupledmodel_cns}, two ideal gas fluids on $\Omega_1$ and $\Omega_2$ are separated by 
the interface $\Gamma = \overline{\Omega}_1 \cap \overline{\Omega}_2$.
Each fluid is governed by 
the CNS equations in $\Omega_\model$ ($\model\in\{1,2\}$), 
\begin{subequations}
\eqnlab{cns-gov}
\begin{align}
  \eqnlab{cns-mass}
  \dd{\rho_\model}{t}     + \Div \LRp{\rho \ub_\model} &= 0,\\
  \eqnlab{cns-momentum}
  \dd{\rho \ub_\model}{t} + \Div \LRp{\rho \ub_\model \otimes \ub_\model + \Ical \pres_\model} &= \Div \sigma_\model - \rho_\model {\bf g},\\
  \eqnlab{cns-energy}
  \dd{\rho E_\model}{t}   + \Div \LRp{\rho \ub_\model H_\model} &= \Div (\sigma_\model \ub_\model) - \Div \Pi_\model - \rho_\model {\bf g} \cdot \ub_\model,
\end{align}
\end{subequations}
\begin{wrapfigure}{r}{0.25\textwidth}
  \centering
  \includegraphics[width=0.24\textwidth]{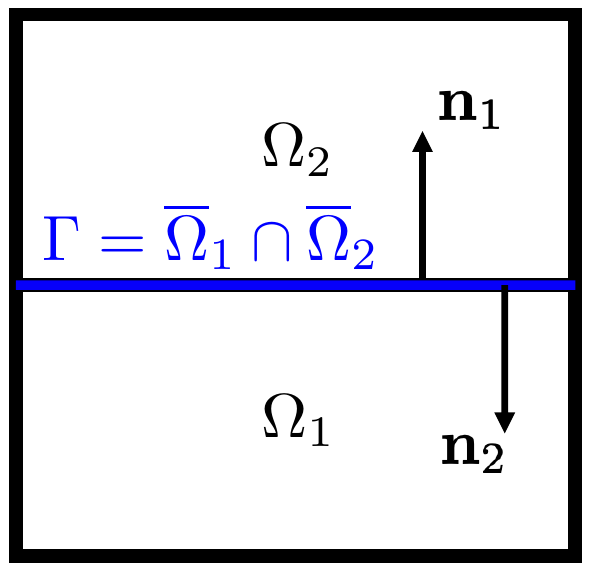}  
  \caption{Schematic of a coupled model.} 
  \figlab{coupledmodel_cns}
\end{wrapfigure}
where $\rho_\model$ is the density $[\si{\kilo\gram\per\cubic\meter}]$;
$\ub_\model$ is the velocity vector $[\si{\meter\per\second}]$;
$\pres_\model$ is the pressure $[\si{\newton\per\meter\squared}]$;
$\rho E_\model = \rho_\model e_\model + \half \rho_\model \norm{\ub_\model}^2 $ is the total energy $[\si{\joule\per\cubic\meter}]$;
$e_\model=\pres_\model \LRp{\rho_\model(\gamma_\model -1)}^{-1}$ 
is the internal energy $[\si{\joule\per\kilogram}]$;
$H_\model = E_\model + \pres_\model\LRp{\rho_\model}^{-1} = a_\model^2\LRp{\gamma_\model -1}^{-1}
+ \frac{1}{2}\norm{\ub_\model}^2$
is the total specific enthalpy $[\si{\joule\per\kilogram}]$;
$\bf g$ is the gravitational acceleration $[\si{\meter\per\second\squared}]$;
$\sigma_\model = \mu_\model \LRp{\Grad \ub_\model + \Grad (\ub_\model)^\top - \frac{2}{3} \Ical \Div \ub_\model } $ is the viscous stress tensor; $\Pi_\model = -\kappa_\model \Grad T_\model$ is the heat flux; 
$T_\model$ is the temperature; 
$\kappa_\model = \mu_\model (c_p)_\model Pr_\model^{-1}$ is the heat conductivity $[\si{\watt\per\meter\per\kelvin}]$;
$\mu_\model$ is the dynamic viscosity $[\si{\pascal\second}]$;
$Pr_\model$ is the Prantl number; 
$a_\model = \LRp{\gamma_\model\frac{\pres_\model}{\rho_\model}}^\half$
is the sound speed $[\si{\meter\per\second}]$ for ideal gas;
$\gamma_\model = \LRp{\frac{c_p}{c_v}}_\model$
is the ratio of the specific heats; 
and $(c_p)_\model$ and $(c_v)_\model$ are the specific heat capacities at constant pressure and at constant volume $[\si{\joule\per\kilogram \per\kelvin}]$, respectively. 
We write \eqnref{cns-gov} in a compact form, 
\begin{align}
\eqnlab{cns-gov-compact}
    \dd{\qb_\model}{t}  + \Div \Fcal_\model^I(\qb_\model) = \Div \Fcal_\model^V(\qb_\model) + \sb_\model (\qb_\model),
\end{align}
with the conservative variable $\qb_\model = (\rho_\model, \rho \ub_\model^\top, \rho E_\model)^\top$,
the source term $\sb_\model = (0,-\rho_\model {\bf g}, - \rho_\model {\bf g} \cdot \ub_\model)^\top$, 
the inviscid flux tensor 
$\Fcal_\model^I=(\rho \ub_\model, \rho \ub_\model \otimes \ub_\model + \Ical \pres_\model, \rho \ub_\model H_\model)^\top$, and 
the viscous flux tensor 
$\Fcal_\model^V=(0, \sigma_\model, \sigma_\model \ub_\model - \Pi_\model)^\top$  in $\Omega_\model$. 
We use the nondimensionalized form \eqnref{cns-gov-nondimension}, 
and we omit the superscript ($*$) from now on. 


\subsection{Finite-volume discretization}

We denote by $\Omega_{\model_h} := \cup_{\ell=1}^\Ne{\model} K_{\model_{\ell}}$ the mesh containing a finite
collection of nonoverlapping elements, $K_{\model_{\ell}}$, that partition $\Omega_\model$.
For example, in a three-dimensional Cartesian coordinate system, we have 
$\Ne{\model} = \Nxe{\model} \times \Nye{\model} \times \Nze{\model}$.
For clarity, we abbreviate the subscript $\model$ in this section.

By integrating \eqnref{cns-gov-compact} over elements, applying the divergence theorem,
and introducing a numerical flux, 
$\nb\cdot \Fcal^*$, we obtain a cell-centered FV scheme for the $\ell$ element,
\begin{align}
\eqnlab{cns-fv}
  \DD{\overline{\qb}_\ell}{t} - \overline{\sb}_\ell
  = - \frac{1}{\snor{\K_\ell}}\int_{\pK_\ell} \nb \cdot \Fcal^* d\pK
  = - \frac{1}{\snor{\K_\ell}}\sum_{f\in\pK_\ell}\int_{f} \nb \cdot \Fcal^* df ,
\end{align}
where $\nb$ is the outward unit normal vector on the boundary $\pK$ of the element $\K$,
$f$ is the elemental face of $\pK$, 
 $\overline{\qb}_\ell = \snor{\K_\ell}^{-1} \int_{\K_\ell} \qb d\K$ is the average state variable 
in $\K_\ell$, and $\snor{\K_\ell}$ is the Lebesgue measure of element $\K_\ell$. 
The numerical flux $\Fcal^*=\Fcal^{I^*} - \Fcal^{V^*}$ is composed of two parts: inviscid and viscous. 
For the inviscid part, we employ the Lax--Friedrich numerical flux,
\begin{align*}
\nb \cdot \Fcal^{I^*} 
 &= \half\LRp{\Fcal^{I}(\qb^l)  + \Fcal^{I}(\qb^r )} \cdot \nb
 + \frac{\lambda}{2}\LRp{\qb^r \nb^r + \qb^l \nb^l} \cdot \nb ,
\end{align*}
where 
$\qb^l$ and $\qb^r$ are the reconstructed values from the left and the right sides of the elemental face $f$, respectively; 
$|A|:=\mc{R} |\Lambda| \mc{R}^{-1}$;
$A:=\dd{\Fcal^I\cdot\nb}{\qb} = \mc{R} \Lambda \mc{R}^{-1}$ is the flux Jacobian; 
$\mc{R}$ and $\Lambda$ are eigenvectors and eigenvalues of the flux Jacobian; 
$\lambda=\max\LRp{|\Lambda(\qb^l)|,|\Lambda(\qb^r)|,|\Lambda(\qb^\texttt{ROE})|}$; 
and $\qb^\texttt{ROE}$ is the Roe average \cite{toro2013riemann}.
For the viscous part, inspired by \cite{nishikawa2011two},
we first compute the common velocity ($\hat{\ub}$), 
common velocity gradient ($\widehat{\Grad \ub}$), and 
common temperature gradient ($\widehat{\Grad T}$) at the elemental face $f$ 
and then evaluate the viscous flux 
$$
\nb \cdot \Fcal^{V^*} := \nb \cdot \Fcal^{V}\LRp{ \hat{\ub},\widehat{\Grad \ub},\widehat{\Grad T} }.
$$

The interface condition and the discretization of a two-dimensional model are explained in our previous study \cite{KANG2021113988}. Thus, we here describe the interface condition and the discretization of a three-dimensional model in \secref{FV3}. Nevertheless, in our three-dimensional model, we exchange only horizontal velocity components across the interface and set the normal velocity component to zero. This simple interface condition is sufficient for demonstrating the parallel performance of the MPRK coupling methods. 
%


\section{Multirate coupling framework}
\seclab{MultirateCoupling}

In this section we introduce the explicit multirate coupling method based on the second-order multirate partitioned Runge--Kutta (MPRK2) method \cite{constantinescu2007multirate}. 
The ocean model is typically solved implicitly because the speed of sound is infinite in incompressible systems, which allows for larger step size than for the atmospheric system. 
In our compressible flow settings, we explicitly solve both the models by using MPRK2 methods. 
To mimic the ratio of time step sizes, we assume that the atmospheric model is stiffer than the ocean model.
The step size used on the ocean side is $m$ times larger than the step size for the atmospheric partition in order to maintain the stability of the coupled explicit method    
($\dtt_{ocn} = m\dtt_{atm}$ with $m=1,2,\cdots$). 
In the multirate context, this implies that the atmosphere is considered as a fast region whereas the ocean is a slow region. 

We decompose the entire domain into three nonoverlapping regions:
fast region (F), slow buffer (B) region, and slow region (S). 
We assume that the stable step size for the slow region is $m$ times larger than that for the fast region; that is, $\dtt^{\LRc{S}} = m \dtt^{\LRc{F}}$.
By discretizing \eqnref{cns-gov-nondimension} using the FV method in \eqnref{cns-fv},
we have the following semidiscretized coupled systems on each region:
\begin{subequations}
\begin{align}
  \DD{\qbfz}{t} & = \Rbfz(\qbfz,\qbsb) \text{ on } \Omega^{\LRc{F}}, \\
  \DD{\qbsb}{t} & = \Rbsb(\qbsb,\qbfz,\qbsz) \text{ on } \Omega^{\LRc{B}}, \\
  \DD{\qbsz}{t} & = \Rbsz(\qbsz,\qbsb) \text{ on } \Omega^{\LRc{S}}. 
\end{align}
\eqnlab{semi-discretized-coupled-model}
\end{subequations}
We update the solutions in $z \in \LRc{S,B,F}$ using the MPRK2 method,
\begin{subequations}
   \eqnlab{MRK2}
  \begin{align}
   \eqnlab{MRK2-qi}
     \Qbz_{n,i} &= \qbz_n + \dtt \sum_{j=1}^{i-1}\az_{ij} \Rbz_{n,j},\quad i=1,2,\hdots,ms,\\
   \eqnlab{MRK2-qn}
     \qbz_{n+1} &= \qbz_n
       + \dtt\sum_{i=1}^{ms} \bz_{i} \Rbz_{n,i}, 
  \end{align}
\end{subequations}
where 
$s$ is the number of stages of the base method of MPRK2 methods, 
$\qbz_n \approx \qbz(t_n)$,
$\Qbz_{n,i} \approx \Qbz(t_n+\cz_i\dt)$, and
$\Rbz_{n,i} \approx \Rbz(t_n+\cz_i\dt)$. 
The scalar coefficients $\az_{ij}$, $\bz_i$, and $\cz_i$ determine all the properties of a given MPRK2 scheme. 
For example, Table \tabref{mrk2-butcher} shows the coefficients for the MPRK2 ($m=2$) coupling method.

We characterize the atmosphere as the fast region and the ocean as the slow region and place the buffer region in the ocean.  
The top layer of the ocean is identified as the slow buffer where 
communication occurs between the fast and the slow buffer regions
at every fast stage.
We choose the buffer size long enough so that 
$\Qbsb_{n,i}=\Qbsb_{n,\text{\texttt{mod}}(i-1,s)+1}$ holds for $i=1,2,\cdots,ms$ at the interface between the buffer and the slow regions. We note that the Butcher tableau for the buffer region requires only the first $s$ stage intermediate states (i.e., $\Qbsz_{n,i}$ for $i=1,2,\cdots,s$) for communication between the slow region and the buffer region.
This allows us to 
collapse $sm$ stages to $s$ stages for the slow region.

Figure \figref{coupling-diagram} illustrates the coupling diagrams, 
with $s=2$ for $m=2$ in Figure \figref{coupling-diagram-m2} and $m=4$ in Figure \figref{coupling-diagram-m4}. Wide orange boxes represent steps on the different regions, and small red and yellow boxes indicate internal stages. The red double arrow means a two-way coupling, where the information is exchanged between regions. The green single arrow denotes a one-way coupling, where an updated ocean solution is copied to the slow buffer region at stage $2k$ (for $k=1,2,\cdots,m)$ for the right-hand side (RHS) function evaluation of $\Rbsb_{n,2k}(\Qbfz_{n,2k},\Qbsb_{n,2k},\Qbsz_{n,2})$
for $k=1,2,\cdots,m$. 
At the $i=2k-1$ stage, we compute the RHS function for the buffer region using $\Qbsz_{n,1}$, that is, $\Rbsb_{n,i}(\Qbfz_{n,i},\Qbsb_{n,i},\Qbsz_{n,1})$. 
As illustrated in Figure \figref{coupling-diagram-m2}, four global stages are needed for one-step integration, $\dtt$. Two-way coupling occurs between regions at the first and the second stages. 
At the third and fourth stages, only the solutions in the fast and the buffer regions are updated. The solutions are advanced with the two-way coupling in the last stage. 
Similarly, in Figure \figref{coupling-diagram-m4}, eight global (fast) stages are required for one-step integration. The atmospheric model needs eight RHS evaluations; however, the ocean model needs only two RHS evaluations. The computational benefits of MPRK2 over its base RK method are evident through this evaluation ratio. 
The MPRK2 coupling algorithm is summarized in Algorithm \algref{MRK2-Coupling}.
  
\begin{remark}
The MPRK2 coupling method \eqnref{MRK2} is mass conservative if the base method is conservative when used to integrate \eqnref{cns-gov-compact} and the buffer size is long enough that $\Qbsb_{n,i}=\Qbsb_{n,\text{\texttt{mod}}(i-1,s)+1}$ holds for $i=1,2,\cdots,ms$ at the interface between the slow buffer and the slow regions \cite{constantinescu2007multirate}. 
\end{remark}
 
 \begin{remark} 
 We note that the base method can be any order of RK methods in the MPRK2 framework, 
but the overall temporal accuracy of MPRK2 is second order at most \cite{constantinescu2007multirate}. 
In this study we use RK2 as our base method for MPRK2, as described in Table \tabref{mrk2-butcher} for $m=2$. 
 \end{remark}


\begin{figure}[h!t!b!]
  \centering
    \begin{subfigure}{0.9\textwidth}
    \centering
      \includegraphics[trim=2.5cm 4cm 3.2cm 4.7cm,clip=true,
        width=0.8\columnwidth]{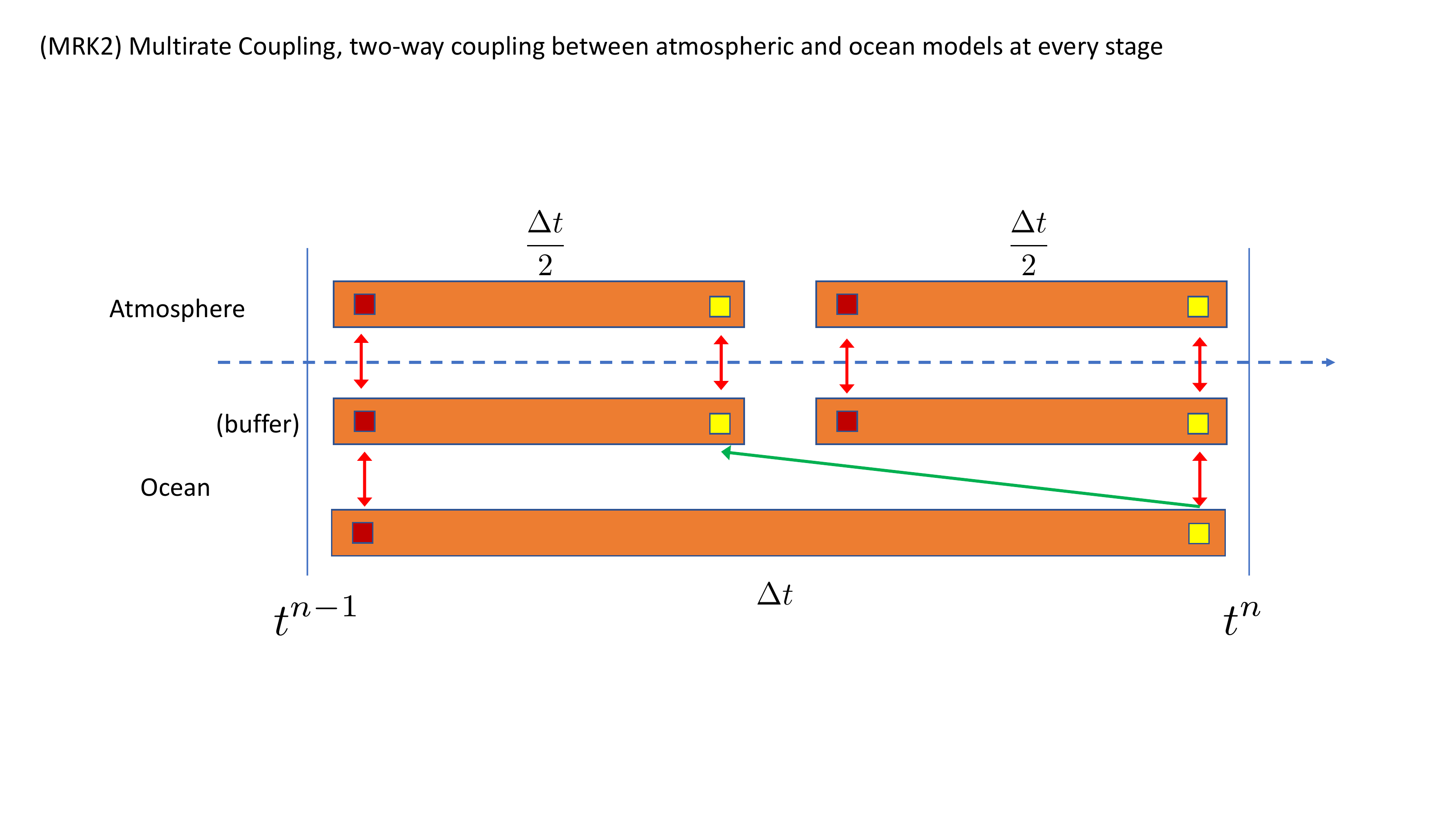}  
        \caption{$m=2$}
        \figlab{coupling-diagram-m2}
    \end{subfigure}%
    
    \begin{subfigure}{0.9\textwidth}
    \centering
      \includegraphics[trim=2.5cm 4cm 3.2cm 4.7cm,clip=true,
        width=0.8\columnwidth]{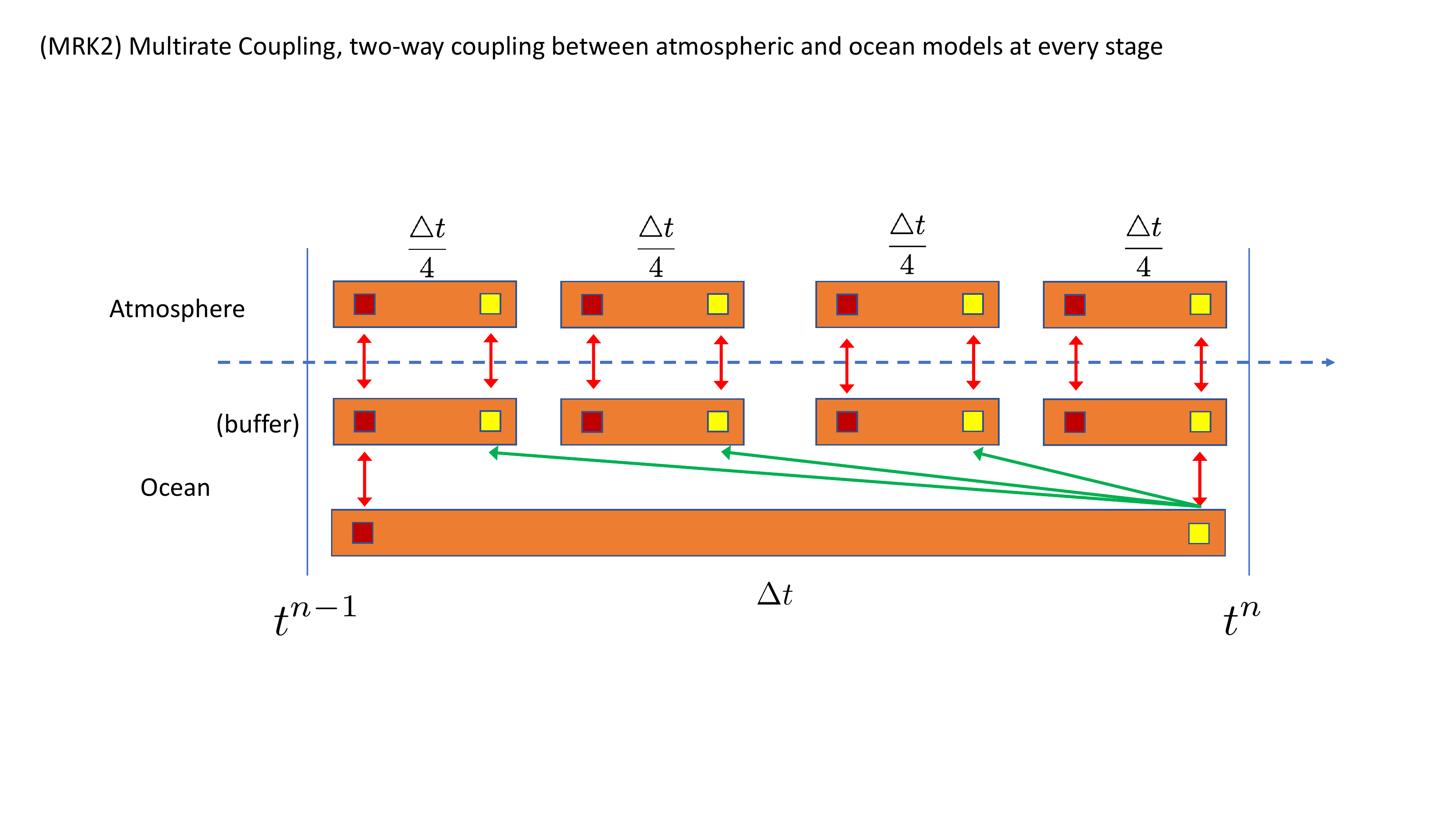}  
        \caption{$m=4$}
        \figlab{coupling-diagram-m4}
    \end{subfigure}%
      \caption{Coupling diagrams with s = 2:
      for (a) $m=2$ and (b) $m=4$. The wide
      (\sqbox1{corange}) boxes represent step sizes, and small (\sqbox1{cred}, \sqbox1{cyellow}) boxes indicate time stages. The red double arrow means a two-way coupling, where the information is exchanged between regions. The green single arrow denotes a one-way coupling, where the updated ocean solution is copied to the buffer for the RHS evaluation. 
      }
   \figlab{coupling-diagram}
\end{figure}

\begin{table}[ht]
    \caption{Butcher tableau for MPRK2 with $m=2$ rate and base method RK2. }
    \tablab{mrk2-butcher} 
    \begin{center}
    \begin{subtable}[c]{0.25\textwidth}    
        \begin{equation*}
            \begin{array}{c|cccc}
                  0 &     &     &     & \\
                \frac{1}{2} & \frac{1}{2} &     &     & \\
                \frac{1}{2} & \frac{1}{4} & \frac{1}{4} &     & \\
                  1 & \frac{1}{4} & \frac{1}{4} & \frac{1}{2} & \\
                \hline
                    & \frac{1}{4} & \frac{1}{4} & \frac{1}{4} & \frac{1}{4}
            \end{array}    
        \end{equation*}
        \caption{Fast (F): $a^{\{F\}}$, $b^{\{F\}}$}
    \end{subtable}
    \begin{subtable}[c]{0.25\textwidth}    
        \begin{equation*}
            \begin{array}{c|cccc}
                0 &   &   &   & \\
                1 & 1 &   &   & \\
                0 & 0 & 0 &   & \\
                1 & 0 & 0 & 1 & \\
                \hline
                & \frac{1}{4} & \frac{1}{4} & \frac{1}{4} & \frac{1}{4}
            \end{array}    
        \end{equation*}
        \caption{Slow buffer (B): $a^{\{B\}}$, $b^{\{B\}}$}
    \end{subtable}
    \begin{subtable}[c]{0.25\textwidth}    
        \begin{equation*}
            \begin{array}{c|cccc}
                0 &      &  &  &\\
                1 &     1&  &  &\\
                0 &     0&  0&  &\\
                1 &     1&  0&  0& \\
                \hline
                  & \frac{1}{2} & \frac{1}{2} & 0& 0  
            \end{array}    
        \end{equation*}
        \caption{Slow (S): $a^{\{S\}}$, $b^{\{S\}}$}
    \end{subtable}
    \begin{subtable}[c]{0.20\textwidth}    
        \begin{equation*}
            \begin{array}{c|cc}
                0 &      &  \\
                1 &     1&  \\
                \hline
                  & \frac{1}{2} & \frac{1}{2}
            \end{array}    
        \end{equation*}
        \caption{Base RK2}
    \end{subtable}
     
    \end{center}
\end{table}

\begin{algorithm}[h!t!b!]
  \begin{algorithmic}[1]
    \ENSURE Given solution state $\qb_n$, compute its next solution state $\qb_{n+1}$.
    The entire domain is decomposed into three regions:
fast region (F), slow buffer (B) region, and slow region (S). 
We assume that the slow step size is $m$ times larger than the fast step size, (i.e.,
$\dtt_{slow} = m \dtt_{fast}$).
      \FOR{$i=1$ to $sm$} 
        \STATE Update the stage value $\Qbsb_{n,i}$ for the buffer region using \eqnref{MRK2-qi}
        \IF {$i \in \LRc{1,2,\cdots,s}$ }
            \STATE Update the stage value $\Qbsz_{n,i}$ for the slow region  using \eqnref{MRK2-qi}
        \ELSE
          \STATE Copy $\Qbsz_{n,j}$ to $\Qbsz_{n,i}$ where $j=\texttt{mod}(i-1,s)+1$
        \ENDIF
        \STATE Update the stage value $\Qbfz_{n,i}$ for the fast region  using \eqnref{MRK2-qi}
        \STATE Compute RHS $\Rbsb_{n,i}(\Qbsb_{n,i},\Qbfz_{n,i},\Qbsz_{n,i})$ for SB
        \IF {$i \in \LRc{1,2,\cdots,s}$ }
            \STATE Compute RHS $\Rbsz_{n,i}(\Qbsz_{n,i},\Qbsb_{n,i})$ for the slow region
        \ENDIF
        \STATE Compute RHS $\Rbfz_{n,i}(\Qbfz_{n,i},\Qbsb_{n,i})$ for the fast region
      \ENDFOR
      \STATE Step completion using \eqnref{MRK2-qn} to compute $\qb_{n+1}$
  \end{algorithmic}
  \caption{Multirate paritioned Runge--Kutta coupling methods}
  \alglab{MRK2-Coupling}
\end{algorithm}

\section{MPRK implementation in a high-performance computing library}
\seclab{MPRKinPETSc}

Although MPRK schemes were proposed over a decade ago, implementation of these schemes had been missing from general-purpose ordinary differential equation solver packages, such as PETSc \cite{AbhyankarEtAl2018,petsc-user-ref} and SUNDIALS \cite{hindmarsh2005sundials}, because of the difficulties in the interface design and buffer region treatment. Nonetheless, recent efforts in both libraries have been deployed, which includes SUNDIALS support for multirate  infinitesimal generalized-structure additive Runge--Kutta \cite{gardner2022enabling}.
To take advantage of high-performance computing, we have implemented the MPRK methods in PETSc and used the callback interface for integrating the semidiscretized models described in Section \secref{MultirateCoupling}.

\subsection*{MPRK implementation in PETSc}

Thanks to the composable design of PETSc solvers, we compose several sub-timestepping (subTS) solvers into one coupled timestepping (TS) solver, with each subTS solver handling a separate model defined on the nonoverlapping regions. In the particular application addressed in this paper, the ocean and the atmosphere models are discretized on two differently structured meshes. 
The domain decomposition strategy in PETSc allows, but does not require, each mesh to be divided among all the processes.
Since three nonoverlapping regions are defined, three callback functions are needed to implement the RHS functions $\Rbfz$, $\Rbsz$ and $\Rbsb$ in \eqnref{semi-discretized-coupled-model} and provided to the main TS solver with the signature function
\newcommand{\func}[1]{\textcolor{Orange}{\textbf{#1}}}
\begin{lstlisting}[language=C++,basicstyle=\ttfamily\small,keywordstyle=\func, otherkeywords={TS,TSRHSSplitSetRHSFunction,Vec,TSRHSFunction}]
int TSRHSSplitSetRHSFunction(TS ts, const char[] name_tag, Vec r, 
                             TSRHSFunction user_callback,void *user_context);
\end{lstlisting}
Each region is assigned with a name tag for easy identification.
An index set needs to be provided for PETSc to access the subvectors (that store the states for each region) in the global solution vector (that stores the states for the combination of all the regions).
\begin{lstlisting}[language=C++,basicstyle=\ttfamily\small,keywordstyle=\func, otherkeywords={RHSSplitSet, IS, TS, Vec}]
int TSRHSSplitSetIS(TS ts, const char[] name_tag, IS index_set);
\end{lstlisting}

\begin{wrapfigure}{r}{0.5\textwidth}
  \centering
  \includegraphics[width=0.5\textwidth]{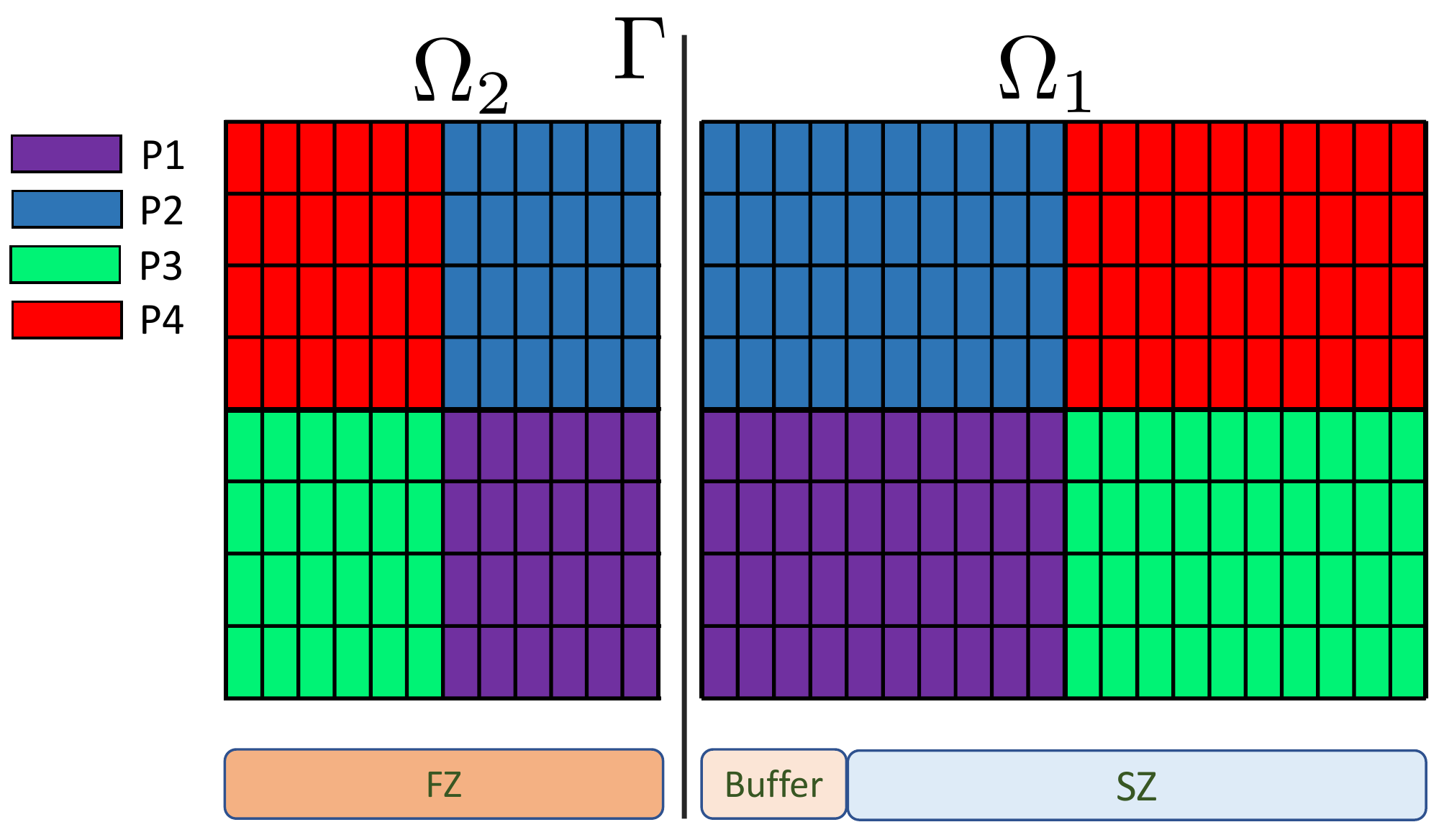}  
  \caption{Schematic of a MPRK2 partitioning with four processes.} 
  \figlab{mprk2-partition}
\end{wrapfigure}
When evaluated at run time, $\Rbfz$ uses all the processes, whereas $\Rbsz$ and $\Rbsb$ together use all the processes because the mesh of the ocean is conceptually split into the slow region and the slow buffer region and the domain decomposition strategy is associated with the mesh only.
For example, Figure \figref{mprk2-partition} shows the domain decomposition strategy used for MPRK2 on four processes.
The subdomains $\Omega_1$ and $\Omega_2$, which have $160$ and $96$ elements, respectively, are equally assigned to four processes, so that each process have $40$ $\Omega_1$ elements and $24$ $\Omega_2$ elements. The partitioning of each mesh is automatically handled by PETSc \texttt{DM}, and MPRK2 does not require any specific partitioning strategy to work properly. \footnote{
\texttt{DM} are abstract objects in PETSc that mediate between meshes, discretizations, the algebraic solvers, time integrators, and optimization algorithms.} The order of the processes can vary for each subdomain in Figure \figref{mprk2-partition}.
A two-way coupling between the two meshes is performed in the $\Rbsb$ callback, with the ocean to atmosphere communication overlapped with the actual $\Rbsb$ computation for efficiency. 

\subsection*{Performance analysis}

Next we analyze the performance of the MPRK2 method over its base RK method. We assume that the RHS function evaluation dominates the computational cost and is proportional to the number of elements in the domain. 
We also assume that the amount of work for the RHS function evaluation per element is the same for the slow, the buffer, and the fast regions. 

Since multirate methods have fewer RHS evaluations in the slow region than single-rate methods have,  a computational benefit can be obtained if the slow region dominates the computational cost. In the following, we first estimate the theoretical speedup for the MPRK2 coupling method over its single-rate (SR) counterpart represented by the base RK2 coupling method and then report the numerical results. 

\textbf{Serial case.}
We denote by $\Nsz$, $\Nbz$, $\Nfz$, and $\Ntotal$ the number of elements in the slow region, the slow buffer region, the fast region, and all the regions (including the slow, the fast, and the slow buffer regions), respectively. For each one-step global integration $\dtt$, SR takes $m$ subcycles with $\frac{\dtt}{m}$ step size in all regions, whereas MPRK2 has one cycle with step size $\dtt$ in the slow region and $m$ subcycles with $\frac{\dtt}{m}$ in the other regions. The number of RHS evaluations per element for MPRK2 and SR is 
\begin{align*}
  N_\texttt{MPRK2} &= 
   s~\Nsz + m s ~\Nbz + m  s ~\Nfz ,\\ 
  N_\texttt{SR}  &= 
   m s ~ \Ntotal. 
\end{align*} 
The theoretical speedup of MPRK2 over its base SR 
is estimated by 
\begin{align} 
 \eqnlab{speedup}
 \texttt{speedup} &= \frac{N_\texttt{SR} }{N_\texttt{MPRK2}}
 = 
 \LRp{1 + \LRp{\frac{1}{m} - 1} \frac{\Nsz}{\Ntotal}  }^{-1}.
\end{align}
This implies that the speedup factor of MPRK2 is proportional to the number of elements in the slow region and grows with increasing step ratio $m$. 

\textbf{Parallel case.}
The two domains that correspond to the fast region and to the slow region and the slow buffer region are decomposed into subdomains and distributed to all MPI processes. We let $N_p$ be the number of MPI processes. 
We denote by $\Nfps$, $\Nszeps$,   
and $\Nbeps$ the  
number of elements per process for the fast region, the 
number of elements per process for the slow region, and the 
number of elements per process for the slow buffer region, respectively. 
We also let $\Nsps$(=$\Nszeps$+$\Nbeps$) be the number of elements per process for both the slow and the buffer regions and let $\Ntotaleps$(=$\Nszeps+\Nbeps+\Nfps$) be the total number of elements per process. 
$\Nfps$ and $\Nsps$ are almost constant, but $\Nszeps$ and $\Nbeps$ can vary 
depending on $\nproc$, $\Nsz$, and $\Nbz$. 
We assume that each process can complete computations of $\varepsilon$ elements per second. 
The speedup does not depend on the value of $\varepsilon$ and $s$; they are introduced simply to facilitate the derivation and will be canceled out in the final formula. The CPU time spent on RHS evaluations for MPRK2 and SR at each global time step ($\dtt$) is
\begin{align*} 
T_\texttt{MPRK2} &= 
 s~ \Nszeps/\varepsilon + m s~ \Nbeps/\varepsilon + m s~ \Nfps/\varepsilon,\\
T_\texttt{SR}  &= 
  m s ~\Ntotaleps/\varepsilon.   
\end{align*}
Note that each process owns a part of the slow region and a part of the slow buffer region. We perform the RHS function evaluations for the slow buffer and the slow region simultaneously for efficiency whenever the RHS for the slow region is needed. This occurs $s$ times during one global step. 
The running time for these evaluations is determined  only by the MPI processes that own a part of the slow buffer, with each evaluation taking $\Nbeps/\varepsilon$ seconds.
Therefore, the theoretical speedup
can be estimated by
\begin{align} 
 \eqnlab{par_speedup}
 \texttt{speedup} = \frac{T_\texttt{SR} }{T_\texttt{MPRK2}} 
    = \LRp{1 + \LRp{\frac{1}{m}-1}\frac{\Nszeps}{\Ntotaleps}}^{-1}.
\end{align}
When considering one MPI process,  
\eqnref{par_speedup} falls back to \eqnref{speedup}. 
We note that as the number of MPI processes increases, 
$\Nszeps$ 
approaches zero gradually, diminishing the speedup. 
When $\Nszeps=0$,
the speedup becomes one, indicating no benefit of using MPRK for this case. 

Next we perform numerical experiments with 1 core to verify the speedup estimation in \eqnref{speedup} by using a three-dimensional wind-driven flow example \cite{KANG2021113988}. We uniformly discretize the whole domain $\Omega=(0,5)\times(0,5)\times(0,10)$ with the elements of $N_{xe}=50$, $N_{ye}=50$, and $N_{ze}=100$.
We vertically split the whole domain into the slow, the buffer, and the fast regions with the ratio of $\Nsz/\Ntotal \in \LRc{4,14,24,34,44,54,64,74,84}/100$ and $\Nbz=6$. 
We hypothetically consider that the fast region is stiffer than the slow region so that the SR method requires $\frac{\dtt}{m}$ step size for the entire domain, whereas the MPRK2 method uses $\dtt$ step size. 
\footnote{
Note that we treat both the ocean and the atmosphere as ideal gas fluids.
In this configuration 
it is hard to introduce scale-separable stiffness to one domain; thus, 
we intentionally reduce the step size of the SR method to $\dtt/m$ even if we can run the coupled model with the SR method of $\dtt$ step size. 
}

\begin{figure}[h!t!b!]
  \centering
    \begin{subfigure}{0.43\textwidth}
    \centering
      \includegraphics[trim=0cm 0.0cm 0cm 0cm,clip=true,
        width=0.95\columnwidth]{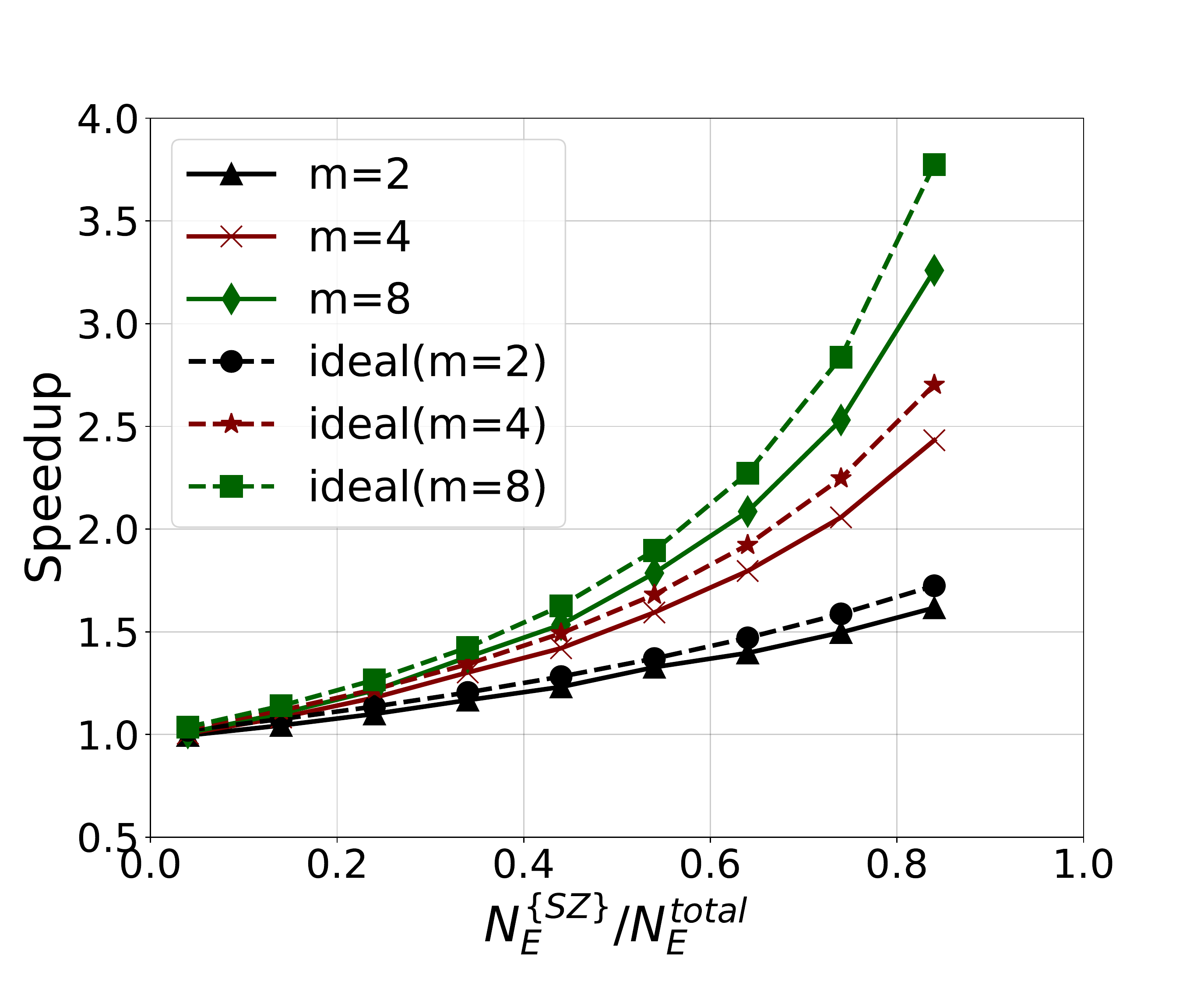}  
        \caption{1 core}
        \figlab{speedup-wdf-1core}
    \end{subfigure}  
    \begin{subfigure}{0.43\textwidth}
    \centering
      \includegraphics[trim=0.0cm 0.0cm 0cm 0cm,clip=true,
        width=0.95\columnwidth]{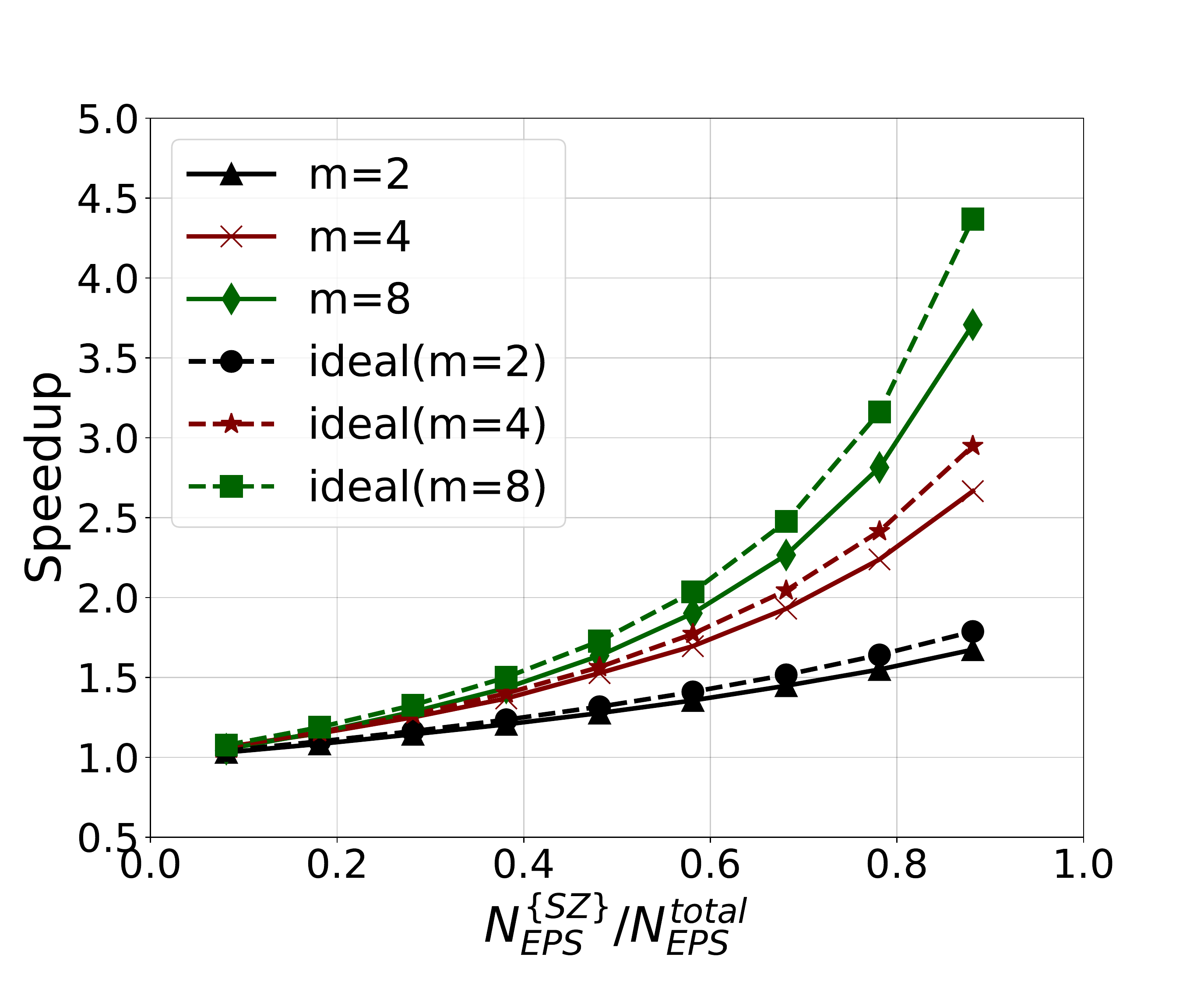}  
        \caption{32 cores}
        \figlab{speedup-wdf-32cores}
    \end{subfigure} 
      \caption{Speedup for multirate partitioned Runge--Kutta method over its base method with (a) 1 core and (b) 32 cores: 
       dashed lines indicate ideal speedup whereas solid lines represent measured speedup.}
  \figlab{multirate-speedup-wdf}
\end{figure}
 
Figure \figref{speedup-wdf-1core}
shows the theoretical speedup of MPRK methods in terms of $m$ and $\Nsz/\Ntotal$. 
As the number of elements in the slow region increases, the speedup grows exponentially. Given a $\Nsz/\Ntotal$ ratio, the speedup also increases with increasing $m$. 
We also report the wall-clock ratios of MPRK2 and its base RK2 in Table \tabref{wallclock-ratio-khi}, 
where the wall-clock ratio is defined by $\texttt{wcr}:=WC_{RK2}/WC_{MPRK2}$ and the ideal speedup by $\texttt{spd}$. 
The measured speedup is bounded above by the ideal speedup estimation in \eqnref{speedup}.
We note that the actual implementation involves  MPRK2 overhead (such as setup time and synchronization at each time level), but the theoretical estimation \eqnref{speedup} does not account for it, which explains the discrepancy between the two.

Similarly, we conduct numerical experiments with 32 cores to verify the speedup estimation \eqnref{par_speedup} in the parallel case. We uniformly discretize the whole domain $\Omega=(0,5)\times(0,5)\times(0,1024)$ with the elements of $N_{xe}=50$, $N_{ye}=50$, and $N_{ze}=10240$. 
We vertically split the whole domain into the slow, the buffer, and the fast regions with the ratio of $\Nsz/\Ntotal$ from $1018/10240$ to $9210/10240$ so that the ratio of $\Nszeps/\Ntotaleps$ varies 
from $65/800$ to $705/800$.
\footnote{
The total number of elements is $256\times10^5$. 
With 32 MPI processes, each core has $\Ntotaleps=8\times10^5$ elements. 
With the buffer size of $6$, the number of elements for the buffer region is $\Nbz=15000$, which is smaller than $\Ntotaleps$.
In our simulation, one core has all the buffer elements 
where $\Nbeps=\Nbz$ and $\Nszeps \in 10^3\times$ $\LRc{65, 145, 225, 305, 385, 465, 545, 625, 705}$.
}
The results are summarized in Figure \figref{speedup-wdf-32cores} and Table \tabref{wallclock-ratio-khi-32cores}. 
We observe similar behavior with 1 core.

We also note that the serial speedup estimation can be viewed as the upper bound of the parallel speedup in the sense of 
$\frac{\Nszeps}{\Ntotaleps}<\frac{\Nsps}{\Ntotaleps}\sim\frac{\Nsz}{\Ntotal}$.
In our implementation, 
$\Nszeps$ and $\Nbeps$ are automatically determined. Thus, we use the serial speedup estimation in the following numerical examples.

\begin{table}[t] 
  \caption{Speedup of MPRK2 over its base RK2 method for $m=\LRc{2,4,8}$ with 1 core, for the 3D wind-driven flow example from Section 4. 
  Wall clocks of MPRK2 and RK2 are summarized. We define the wall-clock ratio  by $\texttt{wcr}:=WC_{RK2}/WC_{MPRK2}$ and the theoretical speedup in \eqnref{speedup} by \texttt{spd}. 
  }
  \tablab{wallclock-ratio-khi}
  \begin{center} 
  \begin{tabular}{*{1}{c}|*{3}{c}|*{3}{c}|*{3}{c}} 
    \hline 
    
    \multirow{2}{*}{$\frac{\Nsz}{\Ntotal}$}
    
    & \multicolumn{3}{c|}{$ m=2  $} 
    
    & \multicolumn{3}{c|}{$ m=4  $} 
    
    & \multicolumn{3}{c}{$ m=8  $} \tabularnewline 
    
    &
    {\small MPRK2}&{\small RK2}&\texttt{wcr}(\texttt{spd}) &
    {\small MPRK2}&{\small RK2}&\texttt{wcr}(\texttt{spd}) &
    {\small MPRK2}&{\small RK2}&\texttt{wcr}(\texttt{spd}) 
    \tabularnewline 
      \hline\hline 

    4/100
    &    487 &    485 &      1.0    (1.0)
    &    964 &    968 &      1.0    (1.0)
    &   1926 &   1941 &      1.0    (1.0)
    \tabularnewline
    
    14/100
    &    464 &    484 &      1.0    (1.1)
    &    895 &    969 &      1.1    (1.1)
    &   1759 &   1939 &      1.1    (1.1)
    \tabularnewline
    
    24/100
    &    441 &    485 &      1.1    (1.1)
    &    826 &    974 &      1.2    (1.2)
    &   1603 &   1947 &      1.2    (1.3)
    \tabularnewline
    
    34/100
    &    429 &    501 &      1.2    (1.2)
    &    771 &   1004 &      1.3    (1.3)
    &   1457 &   2006 &      1.4    (1.4)
    \tabularnewline
    
    44/100
    &    394 &    485 &      1.2    (1.3)
    &    683 &    970 &      1.4    (1.5)
    &   1265 &   1941 &      1.5    (1.6)
    \tabularnewline
    
    54/100
    &    376 &    499 &      1.3    (1.4)
    &    621 &    990 &      1.6    (1.7)
    &   1122 &   2005 &      1.8    (1.9)
    \tabularnewline
    
    64/100
    &    349 &    487 &      1.4    (1.5)
    &    543 &    975 &      1.8    (1.9)
    &    935 &   1950 &      2.1    (2.3)
    \tabularnewline
    
    74/100
    &    323 &    484 &      1.5    (1.6)
    &    470 &    967 &      2.1    (2.2)
    &    763 &   1931 &      2.5    (2.8)
    \tabularnewline
    
    84/100
    &    299 &    484 &      1.6    (1.7)
    &    397 &    967 &      2.4    (2.7)
    &    594 &   1937 &      3.3    (3.8)
    \tabularnewline
    
    \hline\hline 
     
    \end{tabular} 
 
  \end{center} 
\end{table}

\begin{table}[t] 
  \caption{Speedup of MPRK2 over its base RK2 method for $m=\LRc{2,4,8}$ with 32 cores, for the 3D wind-driven flow example from Section 4. 
  Wall clocks of MPRK2 and RK2 are summarized. We define the wall-clock ratio  by $\texttt{wcr}:=WC_{RK2}/WC_{MPRK2}$ and the theoretical speedup in \eqnref{par_speedup} by $\texttt{pspd}$. 
  }
   
  \tablab{wallclock-ratio-khi-32cores} 
  \begin{center} 
  \begin{tabular}{*{1}{c}|*{3}{c}|*{3}{c}|*{3}{c}} 
    \hline 
    
    \multirow{2}{*}{$\frac{\Nszeps}{\Ntotaleps}$}
    
    & \multicolumn{3}{c|}{$ m=2  $} 
    
    & \multicolumn{3}{c|}{$ m=4  $} 
    
    & \multicolumn{3}{c}{$ m=8  $} \tabularnewline 
    
    &
    {\small MPRK2}&{\small RK2}&\texttt{wcr}(\texttt{pspd}) &
    {\small MPRK2}&{\small RK2}&\texttt{wcr}(\texttt{pspd}) &
    {\small MPRK2}&{\small RK2}&\texttt{wcr}(\texttt{pspd}) 
    \tabularnewline 
   
   \hline\hline 
     
    65/800
    
    &      315 &      325 &      1.0   (1.0)
    &      613 &      653 &      1.1   (1.1)
    &     1229 &     1292 &      1.1   (1.1)
    \tabularnewline

    145/800
    
    &      299 &      324 &      1.1   (1.1)
    &      567 &      652 &      1.2   (1.2)
    &     1117 &     1292 &      1.2   (1.2)
    \tabularnewline
    
    225/800
    
    &      284 &      325 &      1.1   (1.2)
    &      521 &      652 &      1.3   (1.3)
    &     1006 &     1293 &      1.3   (1.3)
    \tabularnewline
    
    305/800
    
    &      271 &      327 &      1.2   (1.2)
    &      478 &      654 &      1.4   (1.4)
    &      903 &     1295 &      1.4   (1.5)
    \tabularnewline
    
    385/800
    
    &      255 &      325 &      1.3   (1.3)
    &      429 &      655 &      1.5   (1.6)
    &      788 &     1289 &      1.6   (1.7)
    \tabularnewline
    
    465/800
     
    &      241 &      327 &      1.4   (1.4)
    &      386 &      654 &      1.7   (1.8)
    &      682 &     1296 &      1.9   (2.0)
    \tabularnewline
     
    545/800
     
    &      224 &      325 &      1.4   (1.5)
    &      338 &      653 &      1.9   (2.0)
    &      569 &     1290 &      2.3   (2.5)
    \tabularnewline
    
    625/800
     
    &      210 &      325 &      1.5   (1.6)
    &      292 &      653 &      2.2   (2.4)
    &      459 &     1293 &      2.8   (3.2)
    \tabularnewline
    
    705/800
     
    &      195 &      326 &      1.7   (1.8)
    &      246 &      655 &      2.7   (2.9)
    &      349 &     1294 &      3.7   (4.4)
    \tabularnewline
    
    \hline\hline 
     
    \end{tabular} 

  \end{center} 
\end{table}

\section{Numerical results}
\seclab{NumericalResults}

We perform several numerical experiments on coupled CNS systems to illustrate the properties of the multirate integration discussed above. 
To account for the gravity effect, 
we consider a neutrally stratified atmosphere and ocean with the background potential temperature $\bar{\theta}_0 = \SI{300}{\kelvin}$. 
We first consider Kelvin--Helmholtz examples to verify total mass conservation for MPRK2 coupling methods. 
 We then compare the performance of the MPRK2 method with its base RK2 method through thermal convection examples.
Next we investigate the parallel performance of MPRK2 methods by using three-dimensional wind-driven flow and thermal convection examples. Here we choose $\dtt_{RK2}$ such that two times $\dtt_{RK2}$ causes  numerical instability. We measure the $L_2$ error of $q$ by 
$$ \left\Vert q  - q_r   \right\Vert := 
\LRp{\sum_{\ell=1}^{N_E} \snor{\K_{\ell}} (q_\ell-q_{\ell r})^2 }^\half, $$
where $q_r$ can be an exact solution or a reference solution.

\subsection{Kelvin--Helmholtz instability}
\seclab{sec-khi}

 Kelvin--Helmholtz instability (KHI) is important in the initial process of turbulence or mixing of two fluids in the stratified atmosphere and ocean. 
KHI arises when two fluids have different densities and tangential velocities across the interface. Small disturbances such as waves at the interface grow exponentially, and the interface rolls up into KH rotors \cite{drazin2004hydrodynamic}. 
To see the nonlinear evolution of KHI, 
we add a jet to $\Omega_2$ and place a vortex on $\Omega_1$ to mimic ocean circulation. The initial and boundary conditions can be found in \cite{KANG2021113988}. 

We conduct numerical simulations for MPRK2 with the meshes of $160\times80$ elements on $\Omega_1$ and $160 \times \LRc{80,160,320,640}$ elements on $\Omega_2$ for $t\in[0,500]$. 
As we vertically refine the fast region (atmosphere), the Courant number in $\Omega_2$ grows. Thus, 
we increase the step ratio $m$ (from $1$ to $8$) to relax the geometrical stiffness while keeping the same MPRK2 time step $\dtt_{MPRK2}=0.025$. 
We take RK45 \cite{carpenter1994fourth} solutions as reference and measure the errors as shown in Table \tabref{khi-relerr}. For each of the four meshes of $160 \times 80$ elements on $\Omega_1$ and $160 \times \LRc{80,160,320,640}$ elements on $\Omega_2$, we use the four RK45 time steps ($\dtt_{RK45}=0.05 \LRc{1,\frac{1}{2},\frac{1}{4},\frac{1}{8}}$) to generate the reference solutions. 
We also see that the errors of density, momentum, and energy decrease as the mesh is being refined. 
Figure \figref{history-khi} shows the histories of mass and energy losses. The total mass and total energy changes are defined by $\snor{ \text{mass}(t) - \text{mass}(0)}$ and $\snor{ \text{energy}(t) - \text{energy}(0)}$, respectively. The total mass and the total energy are denoted as $\text{mass} = 
  \sum_{m=1}^2 \sum_{\ell=1}^{\Ne{m}} \overline{\rho }_{m_\ell} \snor{K_{m_\ell}}$ and $\text{energy} = 
  \sum_{m=1}^2 \sum_{\ell=1}^{\Ne{m}} \overline{\rho E}_{m_\ell} \snor{K_{m_\ell}}$.  
Regardless of the temporal rate $m$, the total mass loss is bounded within $\mc{O}(10^{-13})$. The total mass is conserved with MPRK2 methods.

\begin{figure}[h!t!b!]
  \centering
    \begin{subfigure}{0.4\textwidth}
    \centering
      \includegraphics[trim=0.3cm 0.3cm 0.4cm 0.2cm,clip=true,
        width=0.95\columnwidth]{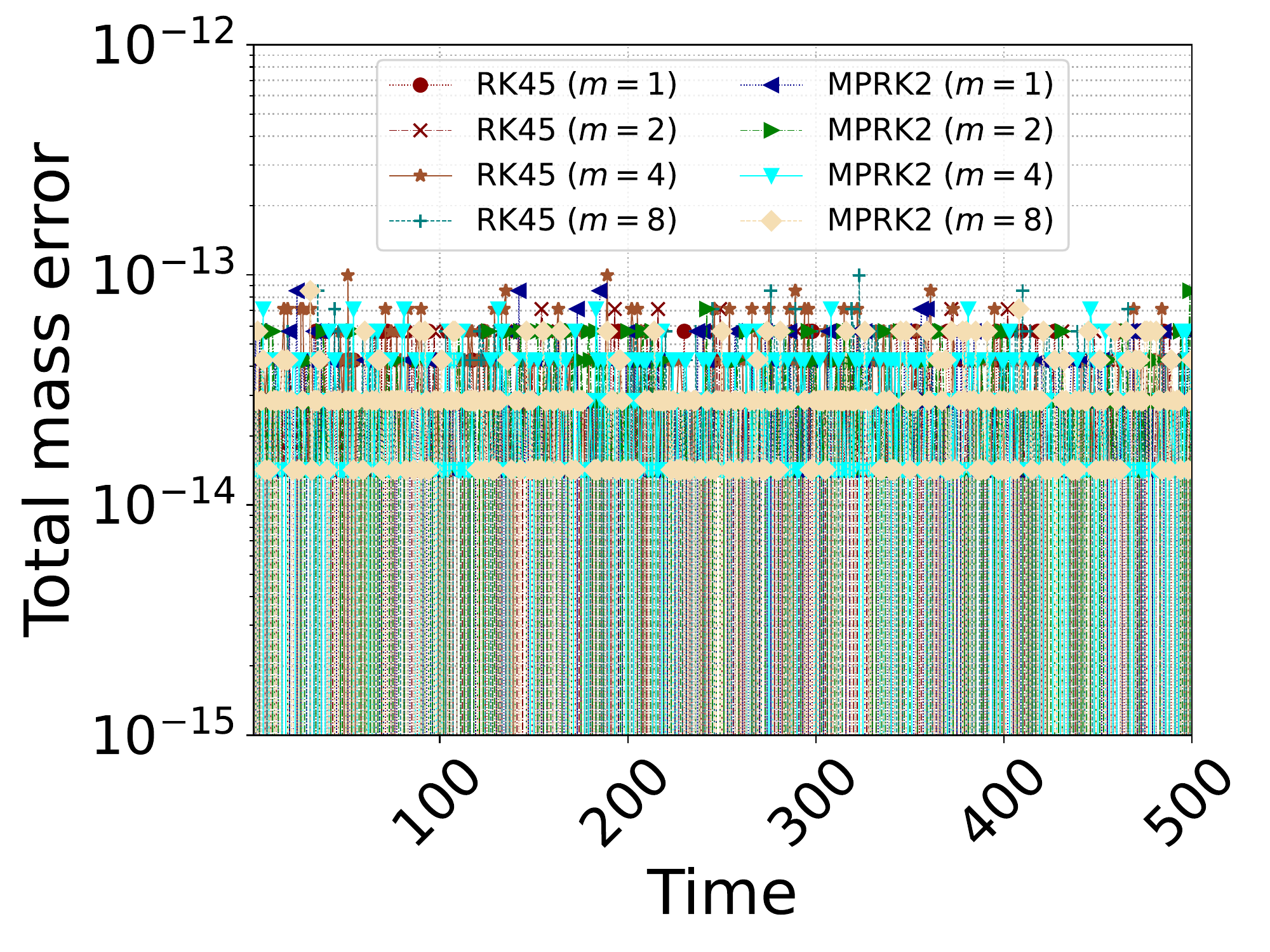}  
        \caption{Mass loss}
        \figlab{khi-mass}
    \end{subfigure}%
    \begin{subfigure}{0.4\textwidth}
    \centering
      \includegraphics[trim=0.3cm 0.3cm 0.4cm 0.2cm,clip=true,
        width=0.95\columnwidth]{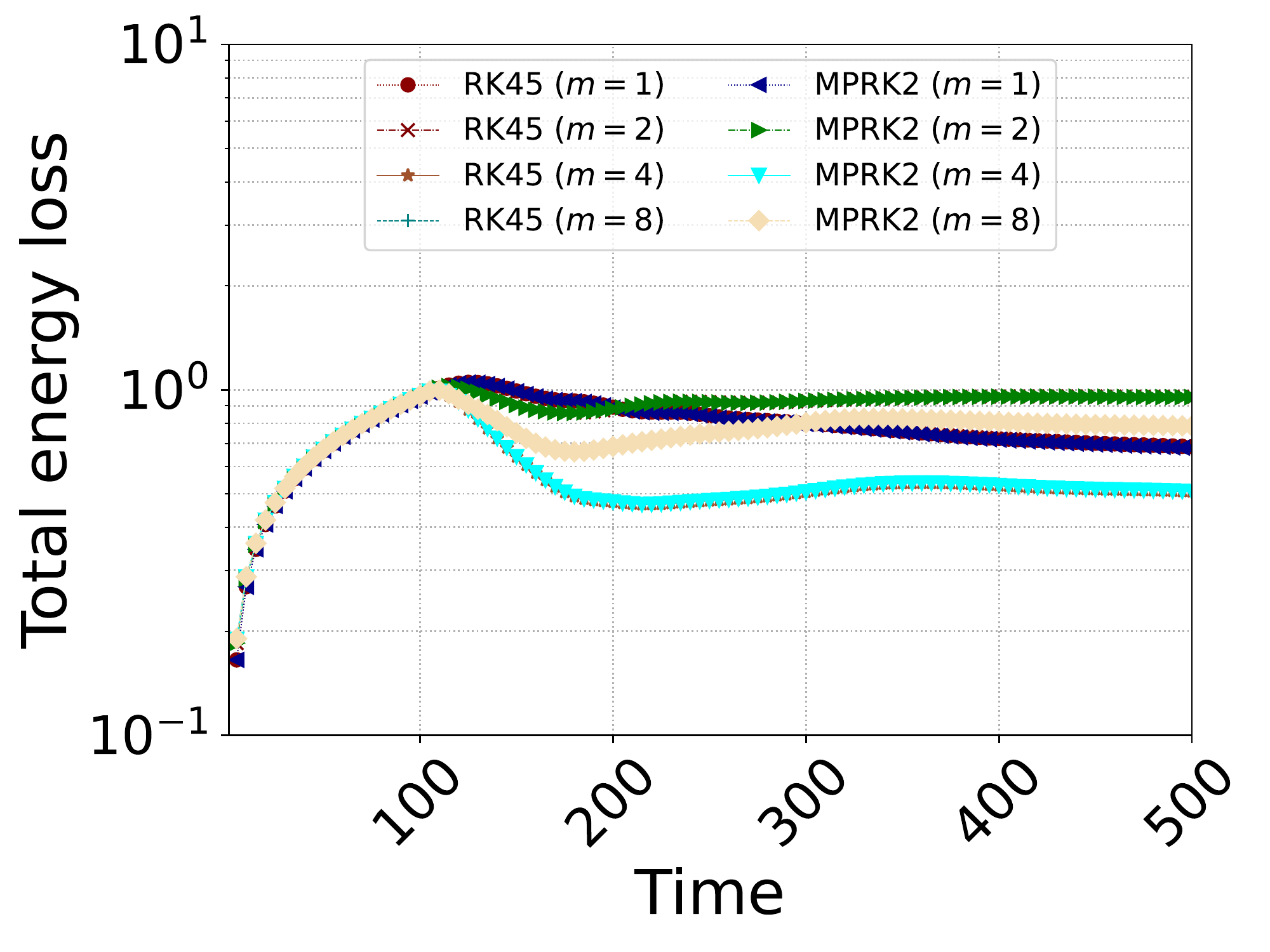}  
        \caption{Energy loss}
        \figlab{khi-energy}
    \end{subfigure}%
      \caption{Wind-driven flows with Kelvin--Helmholtz instability: histories of (a) total mass loss and (b) total energy loss for $t \in [0,500]$. Total mass loss is bounded within $\mc{O}(10^{-13})$. }
      
   \figlab{history-khi}
\end{figure}

\begin{table}[t] 
  \caption{Errors of MPRK2 for KHI example at $t=500$. RK45 solutions are used for reference solutions. $Cr_1$ and $Cr_2$ are the Courant numbers on $\Omega_1$ and $\Omega_2$, respectively; and the subscript $r$ indicates the reference solution.
  }
  \tablab{khi-relerr} 
  \begin{center} 
\begin{tabular}{*{1}{c}|*{1}{c}|*{1}{c}|*{1}{c}|*{1}{c}}  
\hline 
\multirow{2}{*}{ }
& \multirow{2}{*}{$\dtt (Cr_1,Cr_2)$}
& \multirow{2}{*}{$ \left\Vert \rho  - \rho_r \right\Vert  $}  
& \multirow{2}{*}{$ \left\Vert \rho {\bf u} - \rho {\bf u}_r \right\Vert  $}  
& \multirow{2}{*}{$ \left\Vert \rho E - \rho E_r \right\Vert  $}  
\tabularnewline 
& && &  \tabularnewline 
\hline\hline 
MPRK2 ($m=1$) & 0.25 (0.46,0.44)      &       2.52E-02&       2.14E-02&       3.08E-03 \tabularnewline
MPRK2 ($m=2$) & 0.25 (0.46,0.88)      &       5.67E-03&       3.23E-03&       7.14E-04 \tabularnewline
MPRK2 ($m=4$) & 0.25 (0.46,1.76)    &       1.05E-03&       5.63E-04&       2.20E-04 \tabularnewline
MPRK2 ($m=8$) & 0.25 (0.46,3.53)      &       1.54E-04&       1.00E-04&       2.43E-04 \tabularnewline
\hline\hline 
\end{tabular} 
  \end{center} 
\end{table}

\subsection{Thermal convection}
\seclab{sec-rtb}

Thermal convection is the transfer of heat due to the movement of fluid.
It is one of the major forces in the atmosphere and the ocean 
that directly influences the development of ocean currents, clouds, and storm systems. 
To mimic the convection system for a coupled
compressible Navier--Stokes equation, 
we place a warmed fluid perturbation in a neutrally stratified ocean. 

The whole domain is $\Omega =(-5,5) \times (-5,5)$,
comprising two subdomains:   
$\Omega_1 =(-5,5) \times (-5,0)$ and 
$\Omega_2 =(-5,5) \times (0,5)$.  
Adiabatic no-slip conditions are applied to all boundaries on $\Omega_1$ and $\Omega_2$.
The initial conditions are chosen as 
\begin{align*}
 T_\model &= \LRp{1 + \frac{\delta\theta_\model}{\bar{\theta}_0} } \Psi_\model,\\
 P_\model &= \frac{1}{\gamma} \Psi_\model^{\frac{\gamma}{\gamma-1}},\\
 \rho_\model  &= \frac{\bar{\theta}_0}{ \LRp{\bar{\theta}_0 + \delta\theta_\model } }  \Psi_\model^{\frac{1}{\gamma-1}},
\end{align*} 
where $\Psi_\model=1+\frac{gz_\model}{\tilde{c}_p (1 + \delta \theta_\model/\bar{\theta}_0 )}$ and $\bar{\theta}_0 = 300 [K]$ for $\Omega_\model$. 
Here we take 
$g=-0.008140864714$,
$\delta \theta_1 = 0.25 \LRp{1 + \cos (\pi r_1) } $ for $r_1 <= 2.5$, $\delta \theta_1=0$ for $r_1> 2.5$, $r_1:=\norm{\xb_1 - \xb_{1_c}}$, 
$\xb_{1_c}=(0,-2.5)$,$\gamma=1.4$,
 $Pr=0.72$, 
$\tilde{\mu}_1=20000^{-1}$, $\tilde{\mu}_2=5000^{-1}$, 
and  $\delta \theta_2 =0$. 
We note that $T_\model$,$P_\model$, and $\rho_\model$, then $g$ 
and $\Psi_\model$, 
are the normalized quantities.
 
We conduct the simulation with the MPRK2 ($m=2$) method of $\dtt_{MPRK2}=0.025$
over the mesh of $100\times100$ elements on $\Omega_1$ 
  and $100 \times 200$ elements on $\Omega_2$, as shown  in Figure \figref{evolution-rtb}.
The evolution of temperature fields is shown for $t\in\LRs{0,2000}$. 
Since the warmed perturbation in $\Omega_1$ is less dense than its surroundings, 
it rises upward because of buoyancy. 
When it reaches the interface between atmosphere and ocean, 
it moves laterally while exchanging heat and horizontal momentum fluxes, 
which induce horizontal air movement in the lower atmosphere. 
The horizontal-moving fluid hits the lateral wall, causing circulations in both $\Omega_1$ and $\Omega_2$.


\begin{figure}[h!t!b!]
  \centering
     
      \includegraphics[trim=2.5cm 1.5cm 2.2cm 2.2cm,
      clip=true,
        width=0.45\columnwidth]{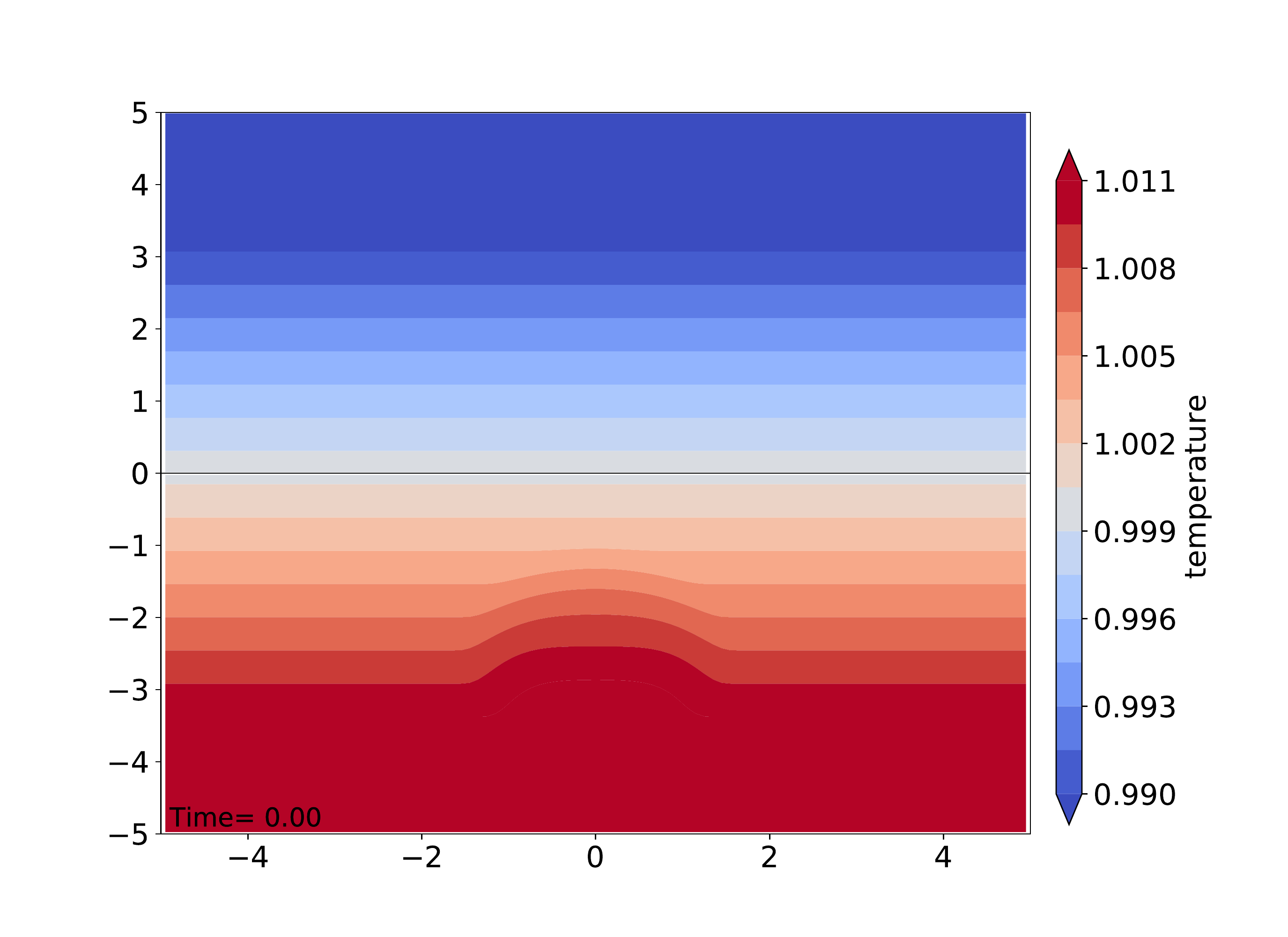} 
      \includegraphics[trim=2.5cm 1.5cm 2.2cm 2.2cm,
      clip=true,
        width=0.45\columnwidth]{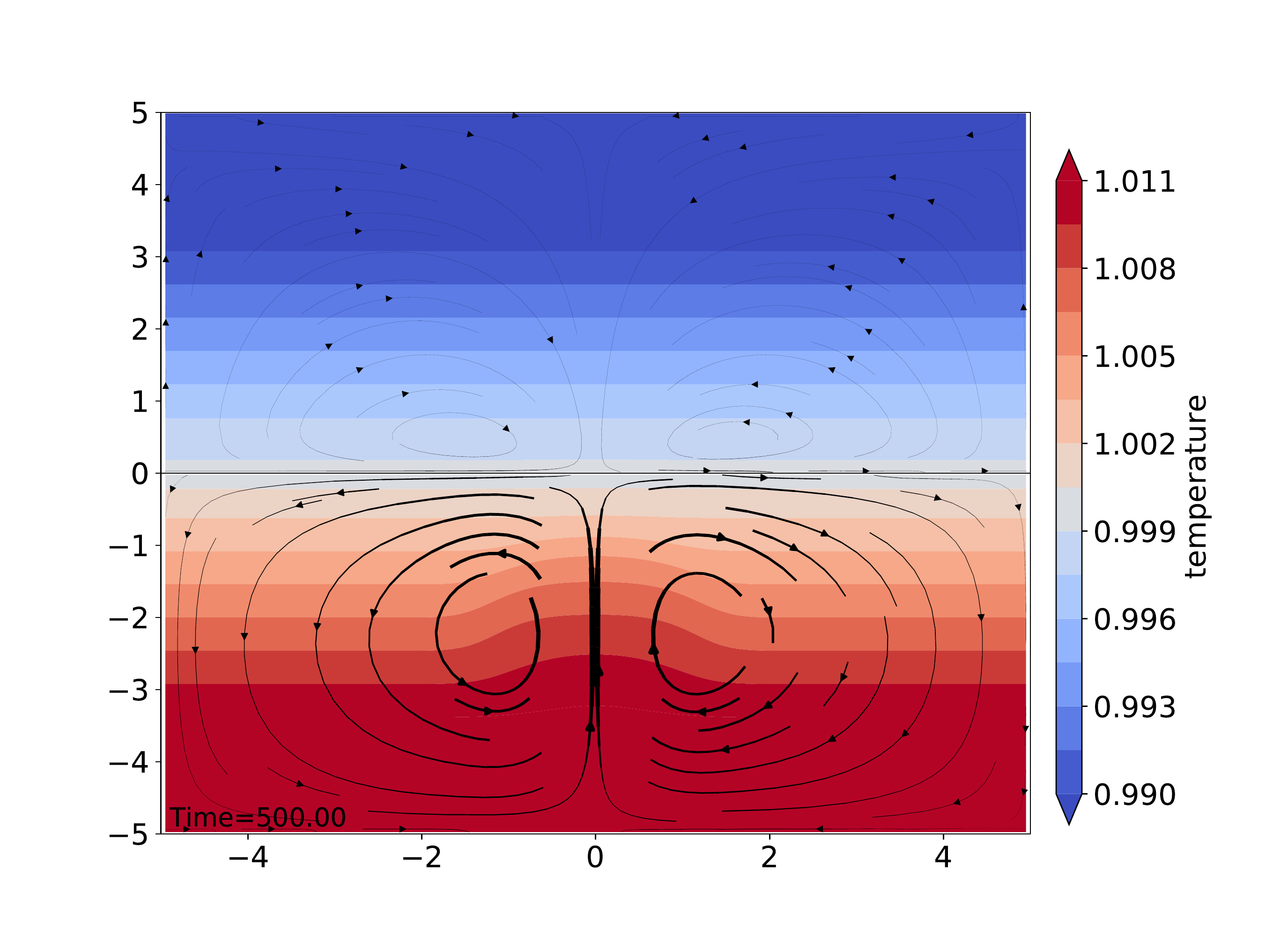} 
    \includegraphics[trim=2.5cm 1.5cm 2.2cm 2.2cm,
      clip=true,
        width=0.45\columnwidth]{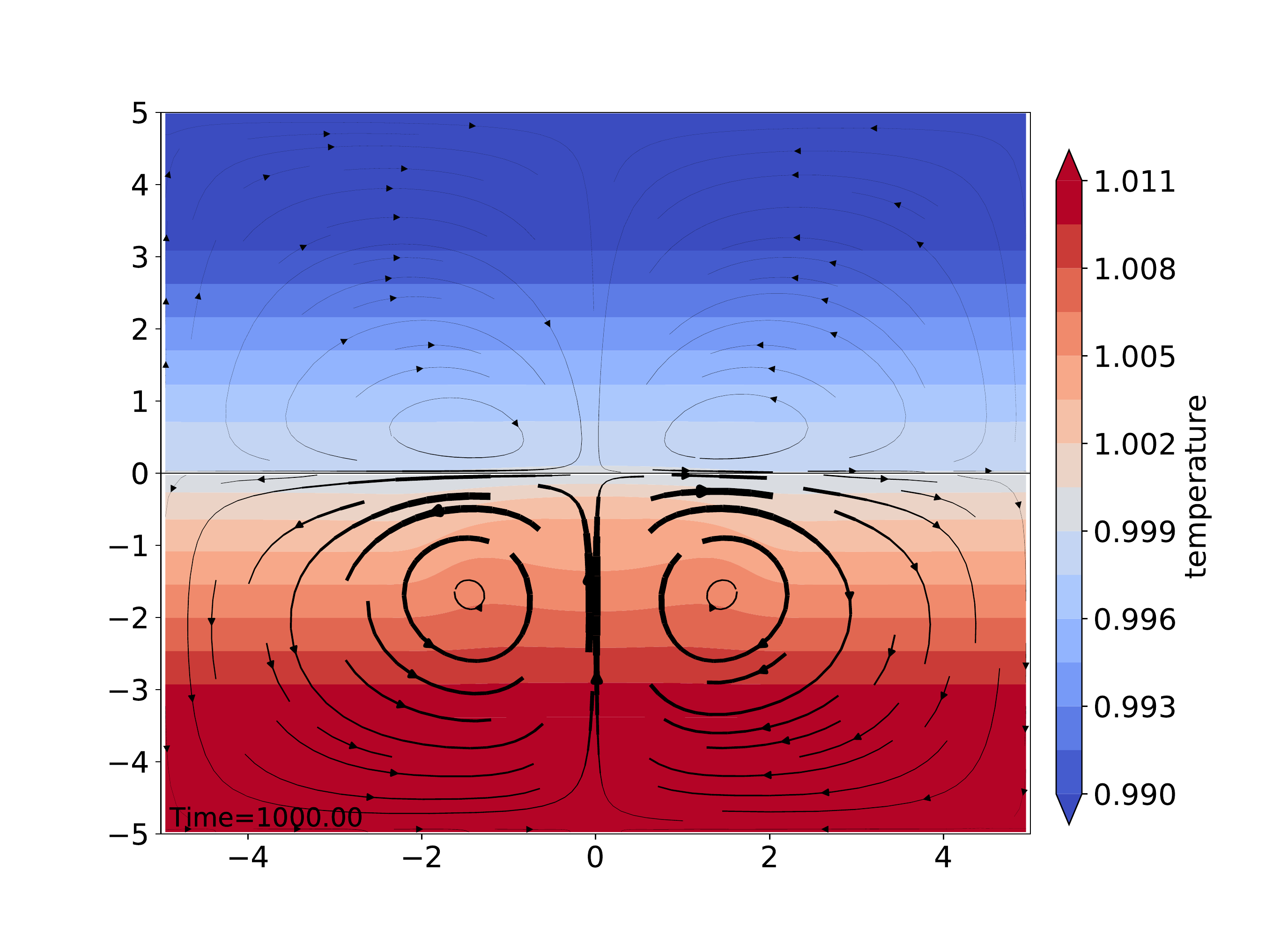}  
    \includegraphics[trim=2.5cm 1.5cm 2.2cm 2.2cm,
      clip=true,
        width=0.45\columnwidth]{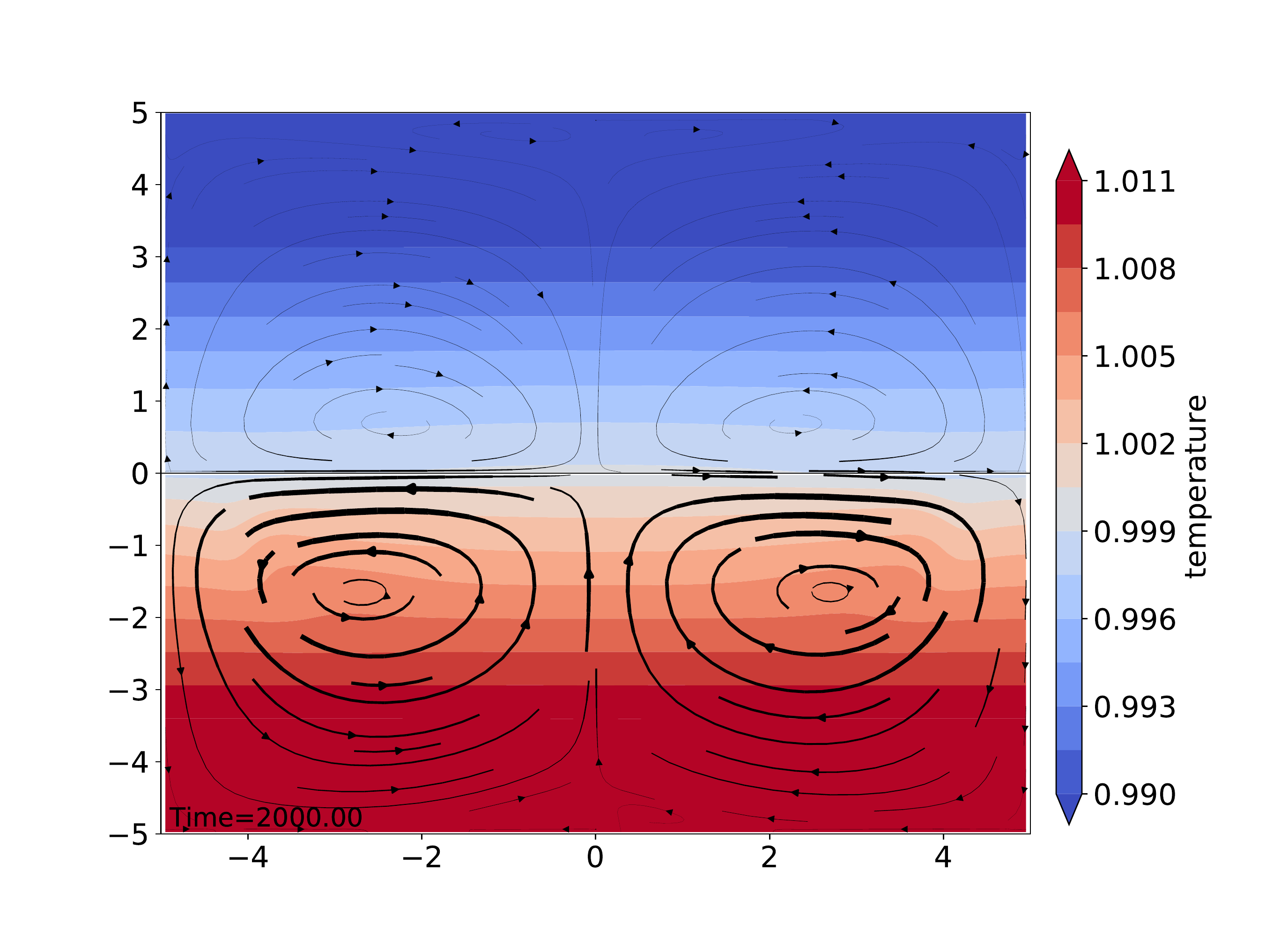}  
      \caption{Evolution of the temperature field for the thermal convection:
  simulation is conducted with the MPRK2 ($m=2$) method over 
 a mesh of $100\times100$ elements on $\Omega_1$ 
  and $100 \times 200$ elements on $\Omega_2$ for $t\in\{0,500,1000,2000\}$. 
  }
   \figlab{evolution-rtb}
\end{figure}

We compare the performance of MPRK2 ($m=2$) and its base method RK2 in terms of accuracy and wall clock. We integrate the systems with RK2 and RK45 with $\dtt_{RK2}=\dtt_{RK45}=0.0125$ over the elements of $100\times100$ on $\Omega_1$ and $100\times200$ on $\Omega_2$. We perform the simulation with the MPRK2 ($m=2$) method of $\dtt_{MPRK2}=0.025$ over the same mesh. Since this example has no exact solution, we take the RK45 solution as ``ground truth" and measure the errors of MPRK2 and RK2. 
Figure \figref{rtb-ss-difference} shows the temperature difference of MPRK2 ($m=2$) and RK45 at $t=2000$ and the temperature difference of RK2 and RK45. We observe that the difference of MPRK2 is larger than that of RK2 in $\Omega_1$ and is within $\mc{O}(10^{-8})$. This comes from the time truncation error associated with the slow part of MPRK2 ($\dt_{slow}=2\times\dt_{RK2}$). 
Table \tabref{rtb-relerr-a} reports the errors of density, momentum, and total energy and the wall clocks for MPRK2 and RK2. The errors of MPRK2 are within $\mc{O}(10^{-6})$. The wall clock of MPRK2 ($m=2$) is comparable to that of its single-rate counterpart. 


\begin{figure}[h!t!b!]
  \centering
    \begin{subfigure}{0.43\textwidth}
    \centering
      \includegraphics[trim=2.5cm 1.5cm 2.2cm 2.2cm,
      clip=true,
        width=0.95\columnwidth]{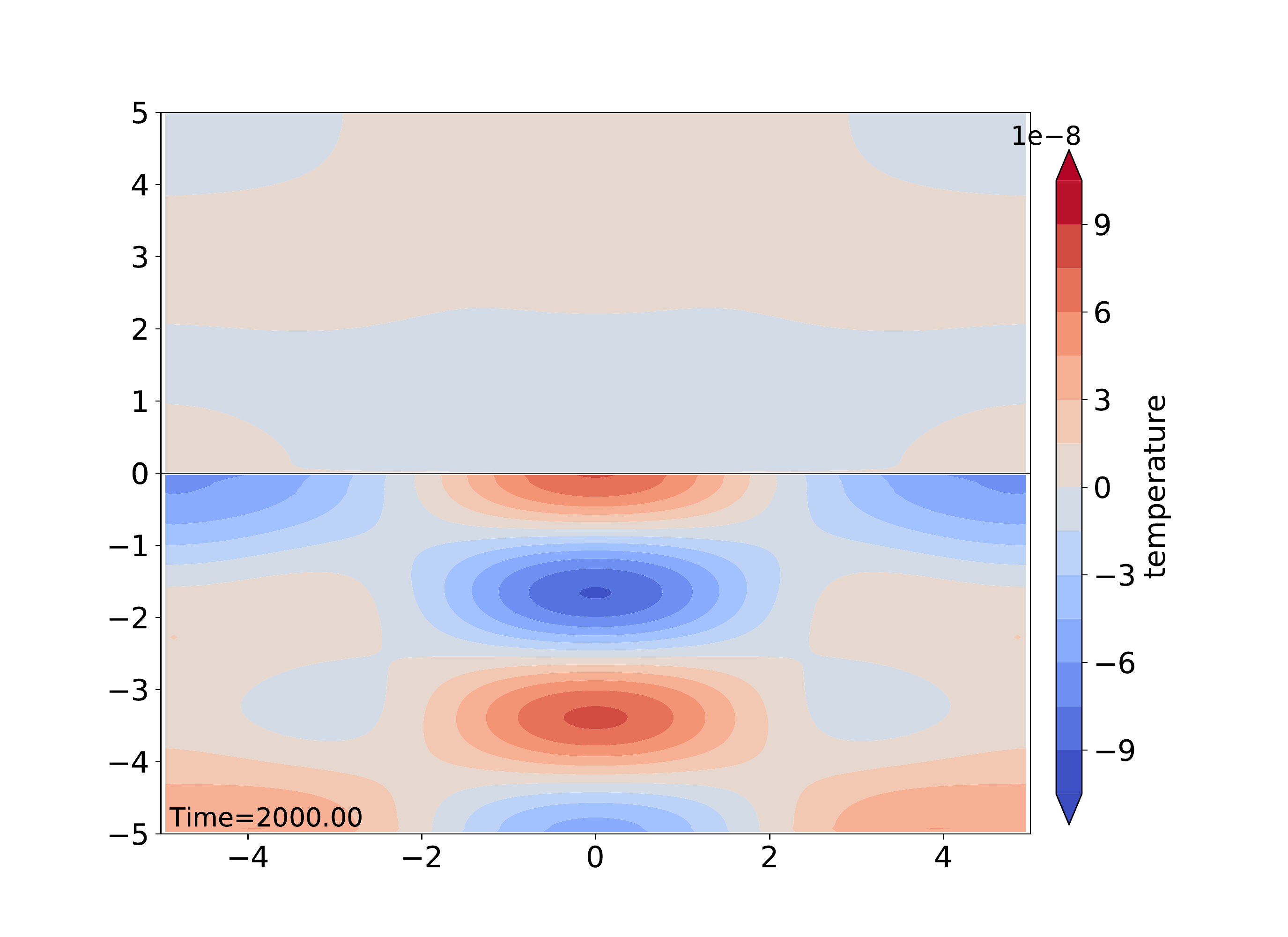}  
        \caption{MPRK2 - RK45}
        \figlab{rtb-ss-mprk2}
    \end{subfigure}%
    \begin{subfigure}{0.43\textwidth}
    \centering
      \includegraphics[trim=2.5cm 1.5cm 2.2cm 2.2cm,
      clip=true,
        width=0.95\columnwidth]{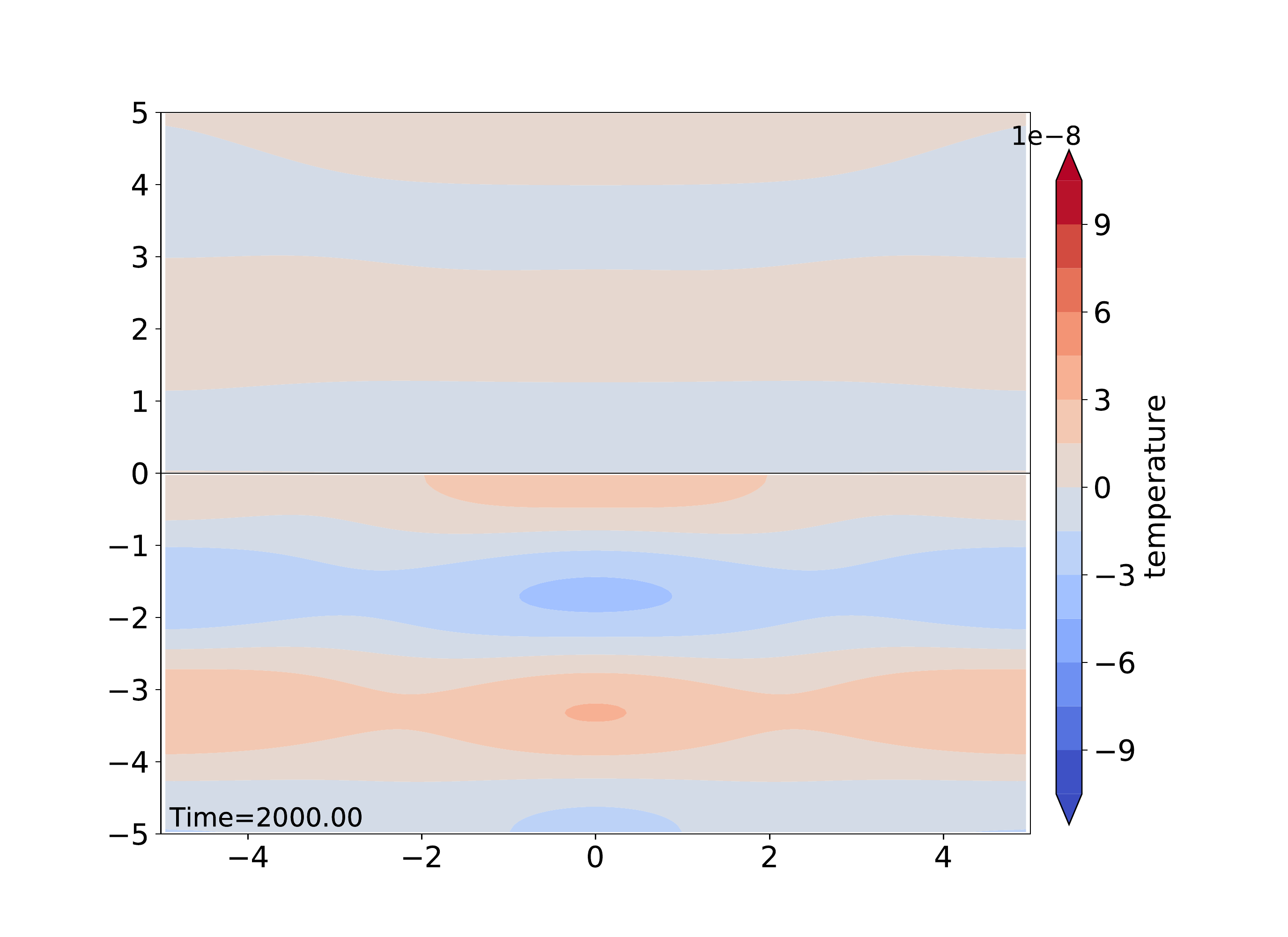}  
        \caption{RK2 - RK45}
        \figlab{rtb-ss-rk2}
    \end{subfigure}%
    \caption{Temperature difference of (a) MPRK2 ($m=2$) and (b) RK2 for the thermal convection example at $t=2000$ with respect to the RK45 solution of $\dtt_{RK45}=0.0125$. 
    We take $\dtt_{MPRK2} =0.025$ and $\dtt_{RK2}=0.0125$.
    }
   \figlab{rtb-ss-difference}
\end{figure}

\begin{table}[t] 
  \caption{Errors of MPRK2 ($m=2$) and RK2 for the thermal convection example at $t=2000$ with respect to the RK45 solution of $\dtt_{RK45}=0.0125$. $Cr_1$ and $Cr_2$ are the Courant numbers on $\Omega_1$ and $\Omega_2$, respectively, and the subscript $r$ indicates the reference solution.
  }
  \tablab{rtb-relerr-a} 
  \begin{center} 
     \begin{tabular}{*{1}{c}|*{1}{c}|*{1}{c}|*{1}{c}|*{1}{c}|*{1}{c}}  
    \hline 
    \multirow{2}{*}{ }
    & \multirow{2}{*}{$\dtt (Cr_1,Cr_2)$}
    & \multirow{2}{*}{$ \left\Vert \rho  - \rho_r \right\Vert  $}  
    & \multirow{2}{*}{$ \left\Vert \rho {\bf u} - \rho {\bf u}_r \right\Vert  $}  
    & \multirow{2}{*}{$ \left\Vert \rho E - \rho E_r \right\Vert  $}  
    & \multirow{2}{*}{wc[s]} \tabularnewline 
    & && & &  \tabularnewline 
    \hline\hline 
    RK2   & 0.0125 (0.25,0.5) &       2.83E-07&       4.44E-07&       7.12E-07& 3887\tabularnewline
    MPRK2 & 0.025 (0.5,1.0)&       5.54E-07&       1.81E-06&       1.40E-06& 3433\tabularnewline
    \hline\hline 
    \end{tabular} 
  \end{center} 
\end{table}


Next we set the domain size to be
$\Omega_1 =(-5,5) \times (-7,0)$ and 
$\Omega_2 =(-5,5) \times (0,3)$. 
We add the potential temperature perturbations in both $\Omega_1$ and $\Omega_2$. We take $\delta \theta_1 = 1.25 \LRp{1 + \cos (\pi r_1) } $ for $r_1 <= 2.5$, $\delta \theta_1=0$ for $r_1> 2.5$,  
$\delta \theta_2 = -7.5 \LRp{1 + \cos (\pi r_2) } $ for $r_2 <= 1$, and $\delta \theta_2=0$ for $r_2> 1$ with 
$\xb_{2_c}=(0,1.5)$.
We perform the simulation with the MPRK2 ($m=4$) method with $\dtt_{MPRK2}=0.025$
over the mesh of $100\times100$ elements on $\Omega_1$ 
  and $100 \times 200$ elements on $\Omega_2$, as shown in Figure \figref{evolution-rtbdc}.
The evolution of temperature fields is shown for $t\in\LRs{0,2000}$. 
The cold fluid parcel in the atmosphere drops down to the interface while the warm fluid in the ocean rises. The cold and warm perturbations move horizontally, balancing heat and momentum fluxes across the interface and hitting the walls. 
The cooled fluid on the ocean surface begins to sink and create circulations. 


\begin{figure}[h!t!b!]
  \centering
    
      \includegraphics[trim=2.5cm 1.5cm 2.2cm 2.2cm,
      clip=true,
        width=0.45\columnwidth]{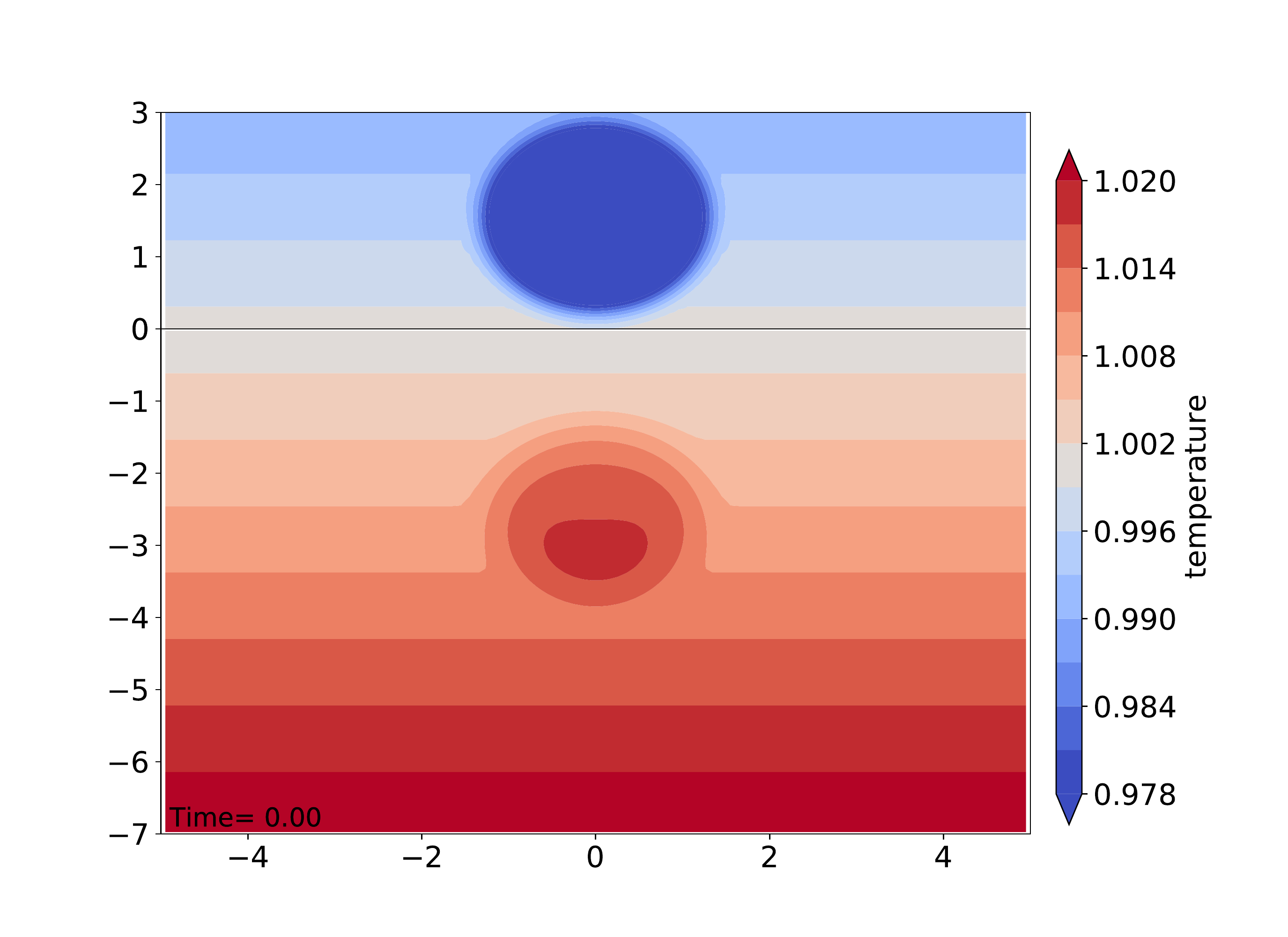}  
     \includegraphics[trim=2.5cm 1.5cm 2.2cm 2.2cm,
      clip=true,
        width=0.45\columnwidth]{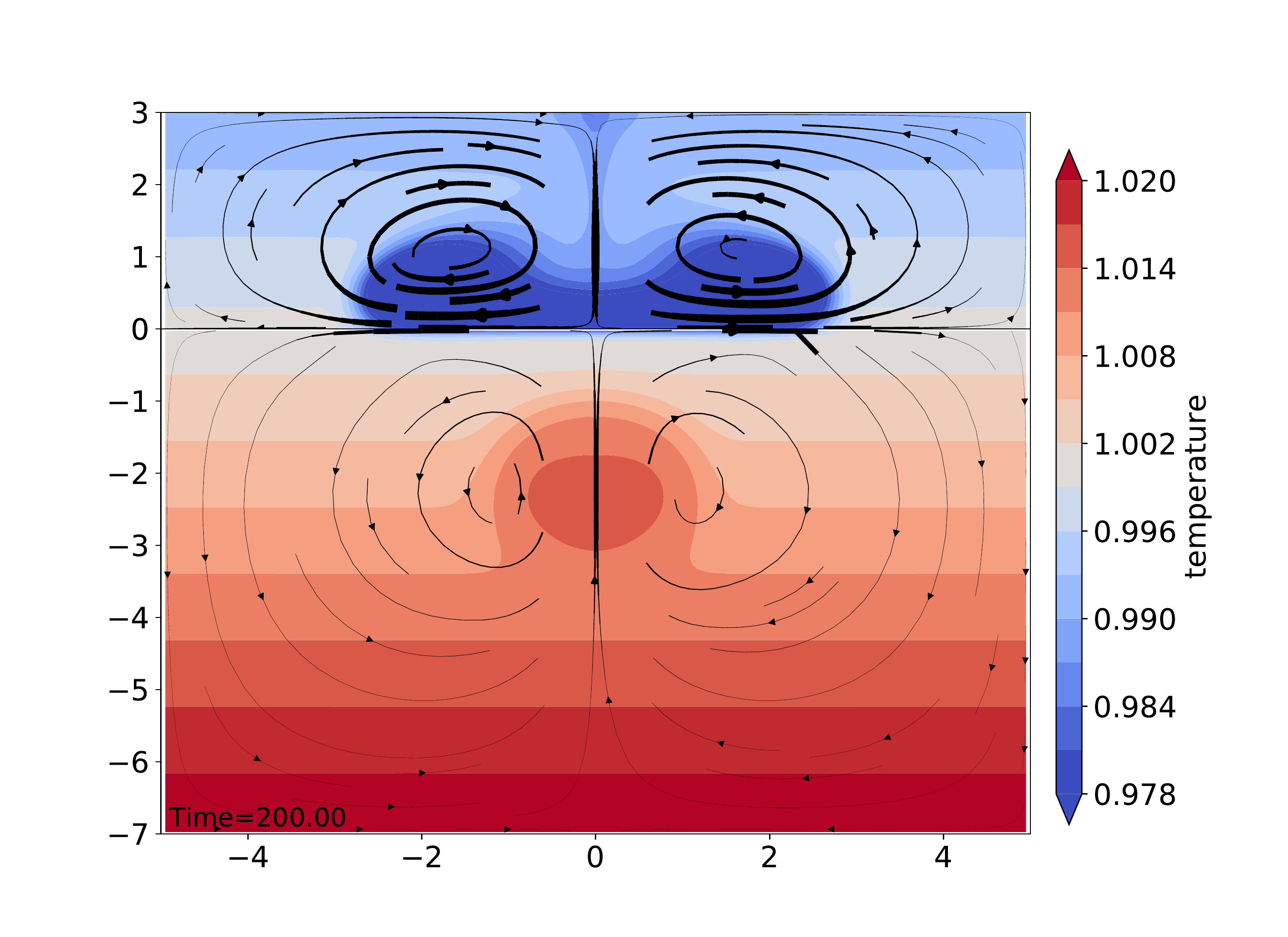} 
    \includegraphics[trim=2.5cm 1.5cm 2.2cm 2.2cm,
      clip=true,
        width=0.45\columnwidth]{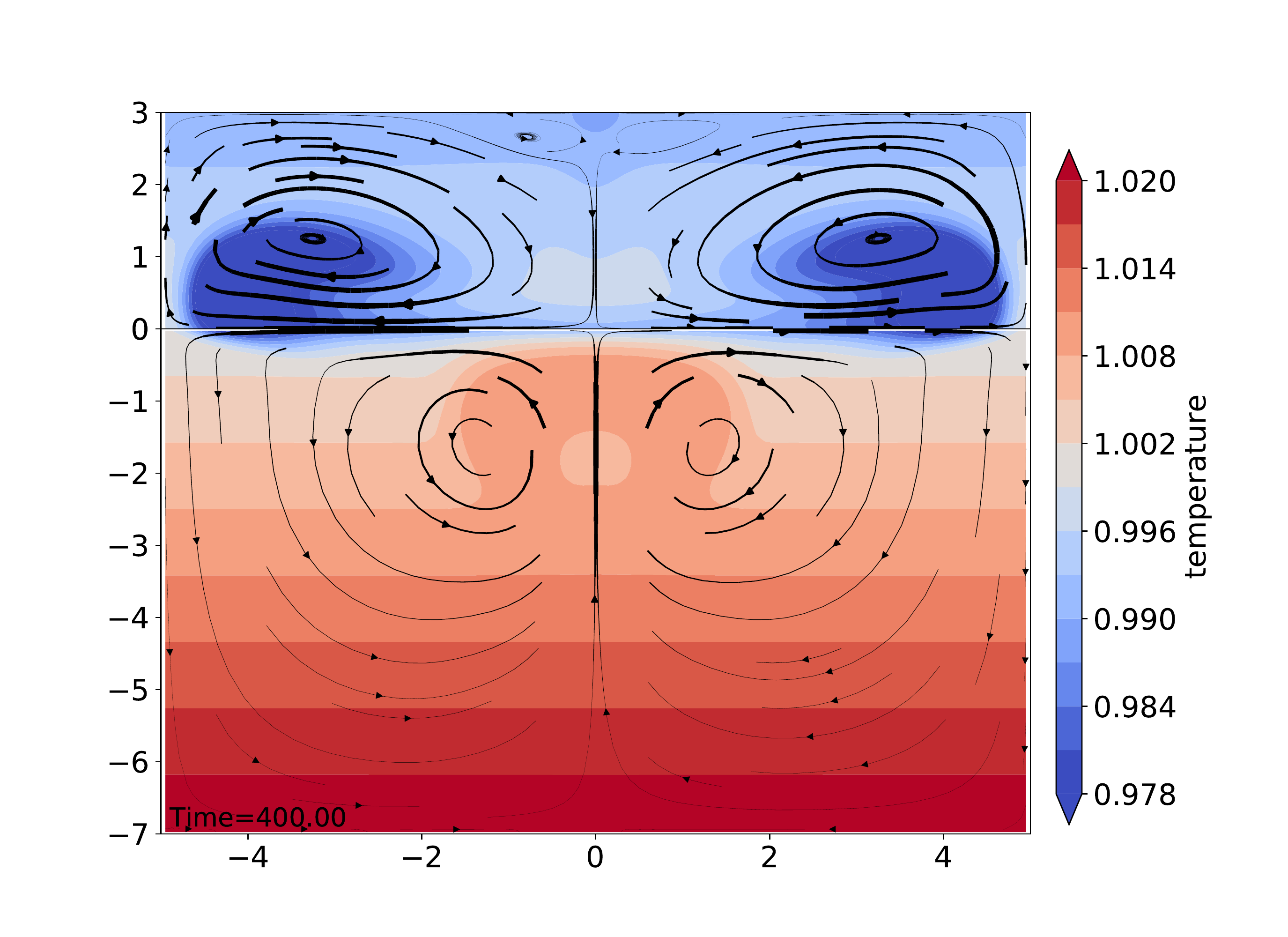}  
    \includegraphics[trim=2.5cm 1.5cm 2.2cm 2.2cm,
      clip=true,
        width=0.45\columnwidth]{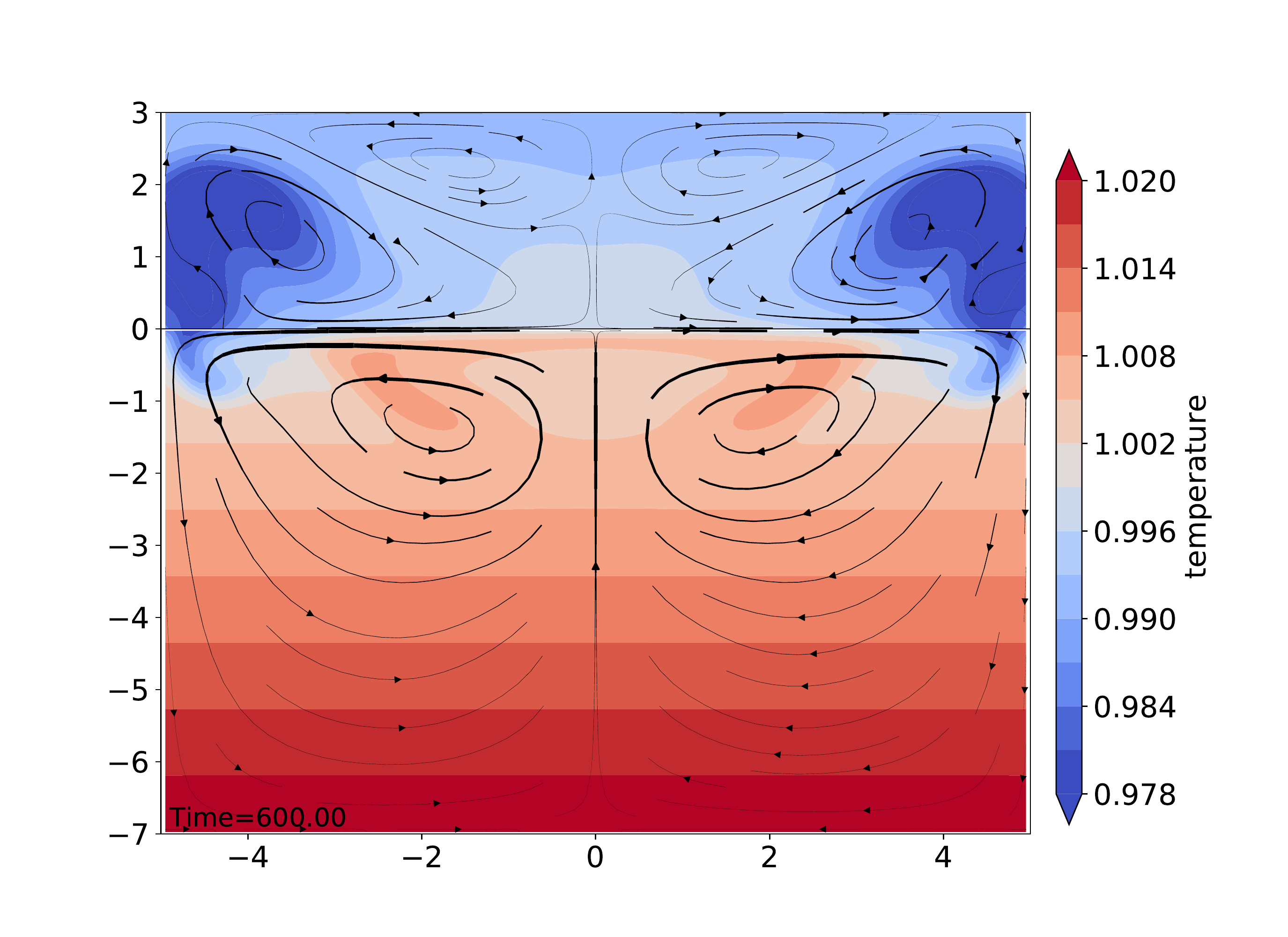}  
    \includegraphics[trim=2.5cm 1.5cm 2.2cm 2.2cm,
      clip=true,
        width=0.45\columnwidth]{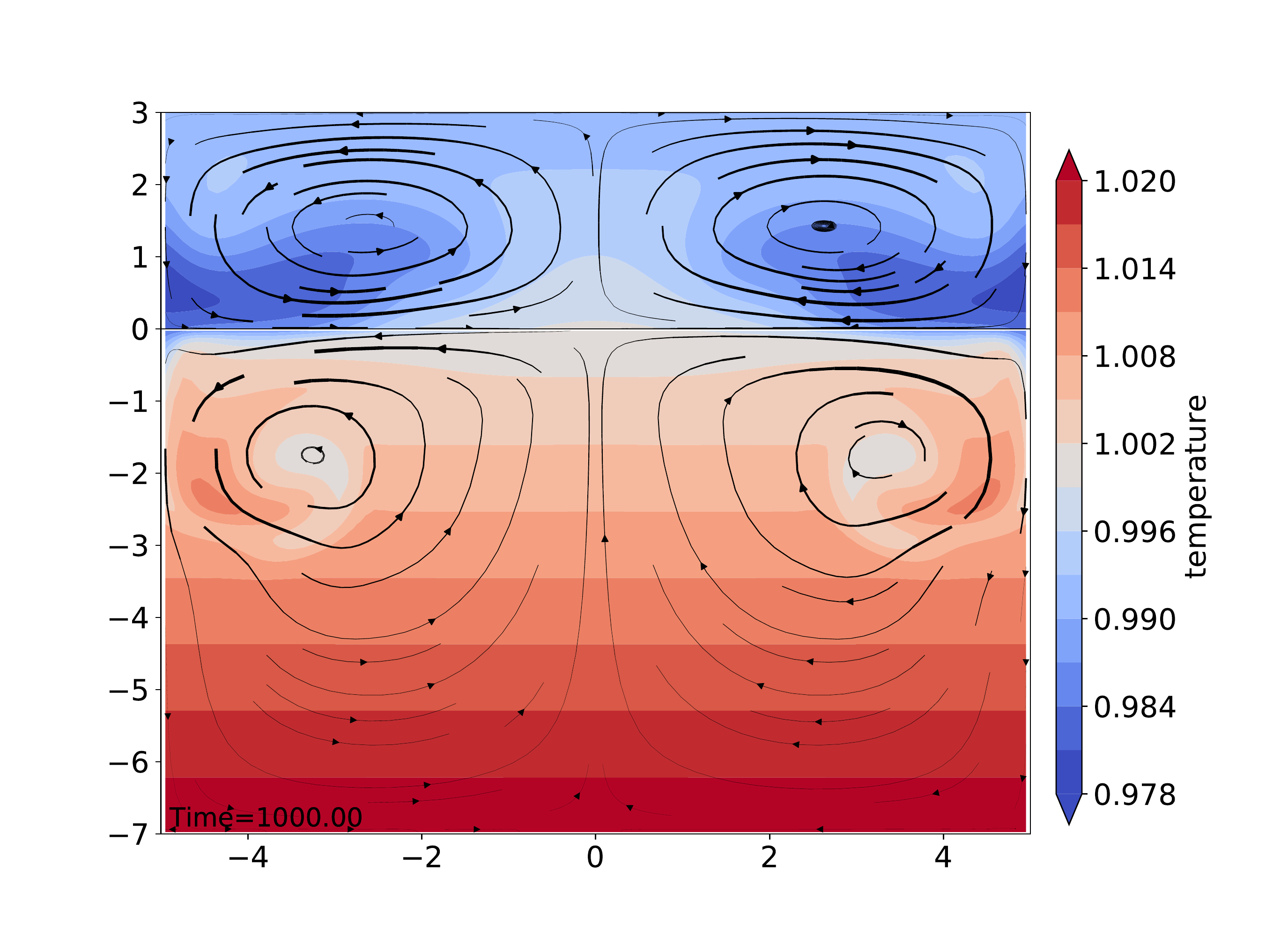}  
    \includegraphics[trim=2.5cm 1.5cm 2.2cm 2.2cm,
      clip=true,
        width=0.45\columnwidth]{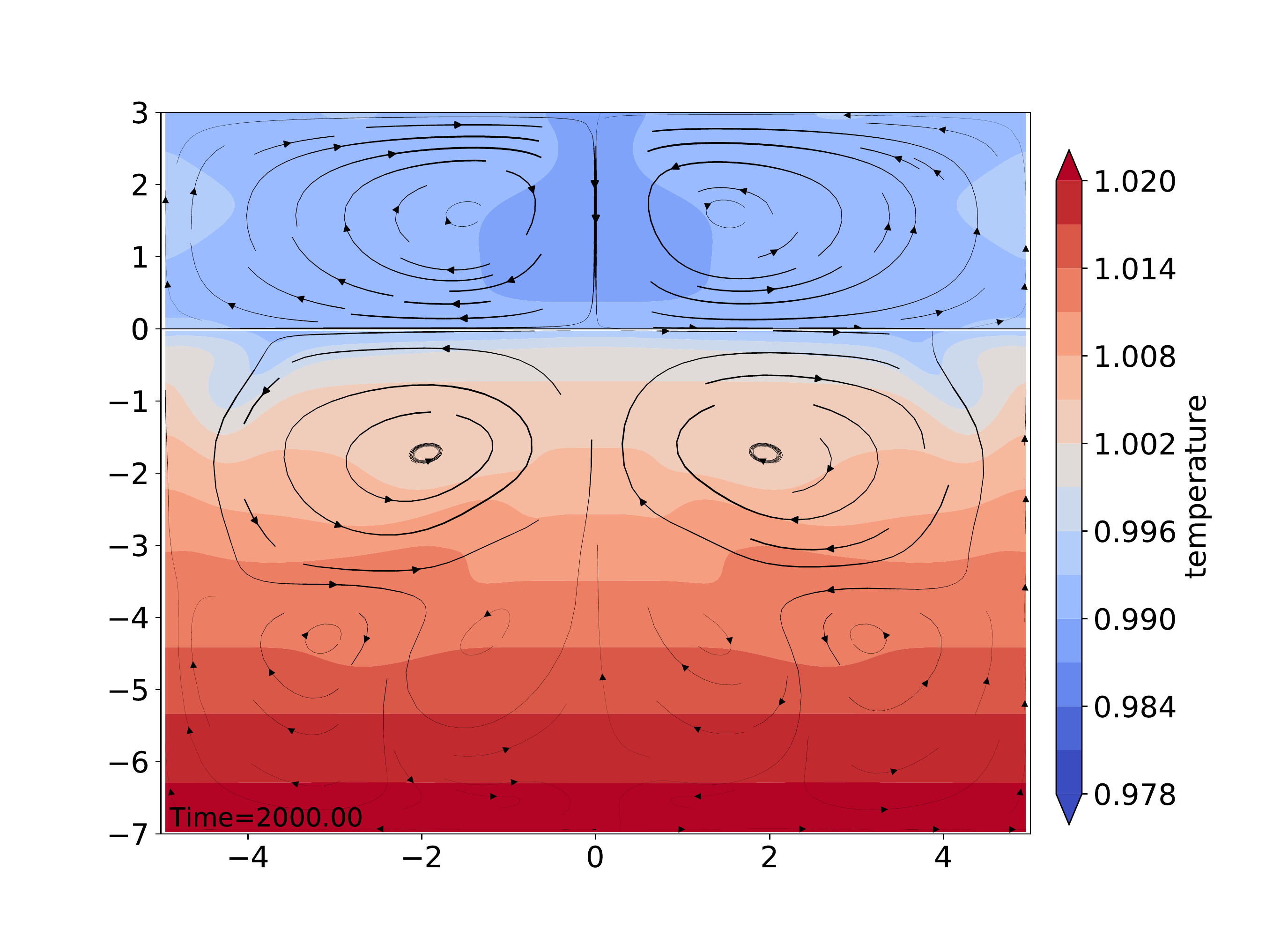}                  
      \caption{Evolution of the temperature field for the thermal convection:
  simulation is conducted with the MPRK2 ($m=4$) method over 
 a mesh of $100\times140$ elements on $\Omega_1$ 
  and $100 \times 240$ elements on $\Omega_2$ for $t\in\LRc{0,200,400,600,1000,2000}$. 
  }
   
   \figlab{evolution-rtbdc}
\end{figure}

\begin{table}[t] 
  \caption{Temporal convergence study of MPRK2 ($m=4$) for the thermal convection example. 
  We take the RK45 solution with $\dtt_{RK45}=0.00125$ as the ground truth solution, and we measure the $L_2$ errors of density, momentum, and total energy at $t=2.5$. The numerical solutions of MPRK2 converge to the reference solution with second-order accuracy. $Cr_1$ and $Cr_2$ are the Courant numbers on $\Omega_1$ and $\Omega_2$, respectively, and the subscript $r$ indicates the reference solution.}
  \tablab{tconv-rtb} 
  \begin{center} 
  \begin{tabular}{*{1}{c}|*{2}{c}|*{2}{c}|*{2}{c}} 
    \hline 
    \multirow{2}{*}{$\dtt (Cr_1,Cr_2)$}
    & \multicolumn{2}{c}{$ \left\Vert \rho   - \rho_r   \right\Vert  $} 
    & \multicolumn{2}{c}{$ \left\Vert \rho {\bf u} - \rho {\bf u}_r \right\Vert  $} 
    & \multicolumn{2}{c}{$ \left\Vert \rho E - \rho E_r \right\Vert  $} \tabularnewline 
    & error & order &error & order &error & order \tabularnewline 
    \hline\hline 
    0.025 (0.51,2.00)&         2.31E-06 & $-$&       2.33E-06 & $-$&       5.85E-06 & $-$\tabularnewline
    0.0125(0.25,1.00)&        5.61E-07 &    2.04&       5.64E-07 &    2.04&       1.42E-06 &    2.04\tabularnewline
    0.00625(0.13,0.51)&       1.40E-07 &    2.00&       1.41E-07 &    2.00&       3.54E-07 &    2.00\tabularnewline
    0.003125(0.07,0.25)&      3.49E-08 &    2.00&       3.51E-08 &    2.00&       8.84E-08 &    2.00\tabularnewline
    \hline\hline 
    \end{tabular} 
  \end{center} 
\end{table} 
For a temporal convergence study, we perform numerical simulations for MPRK2 ($m=4$) 
with $\dtt_{MPRK2} \in 0.025\LRc{\frac{1}{8},\frac{1}{4},\frac{1}{2},1}$. We measure the $L_2$ errors of density, momentum, and total energy at $t=2.5$ with the RK45 solution of $\dtt_{RK45}=0.00125$. The results are summarized in Table \tabref{tconv-rtb}. Second-order convergence rates are observed for the conservative variables as expected.

For the performance comparison, 
we conduct the simulations with RK45 and RK2 single rate methods  
with $\dtt_{RK45}=0.0125$ and $\dtt_{RK2}=0.00625$. 
Figure \figref{rtbdc-ss-difference} shows the temperature difference of MPRK2 ($m=4$) and RK45 at $t=2000$ and the temperature difference of RK2 and RK45. Similar to the previous example, we see that the difference of MPRK2 in the bottom fluid is larger than that of RK2 and is within $\mc{O}(10^{-7})$. 
In Table \tabref{rtbdc-relerr} we summarize the errors of density, momentum, and total energy and the wall clocks for MPRK2 and RK2. 
The errors of MPRK2 are within $\mc{O}(10^{-6})$. The wall clock of RK2 is about 1.5 times than that of the MPRK2 ($m=4$) counterpart. 
Compared with the previous example with $m=2$, MPRK2 ($m=4$) shows better improvement. This agrees with the speedup estimation in \eqnref{speedup}.


\begin{figure}[h!t!b!]
  \centering
    \begin{subfigure}{0.43\textwidth}
    \centering
      \includegraphics[trim=2.5cm 1.5cm 2.2cm 2.2cm,
      clip=true,
        width=0.95\columnwidth]{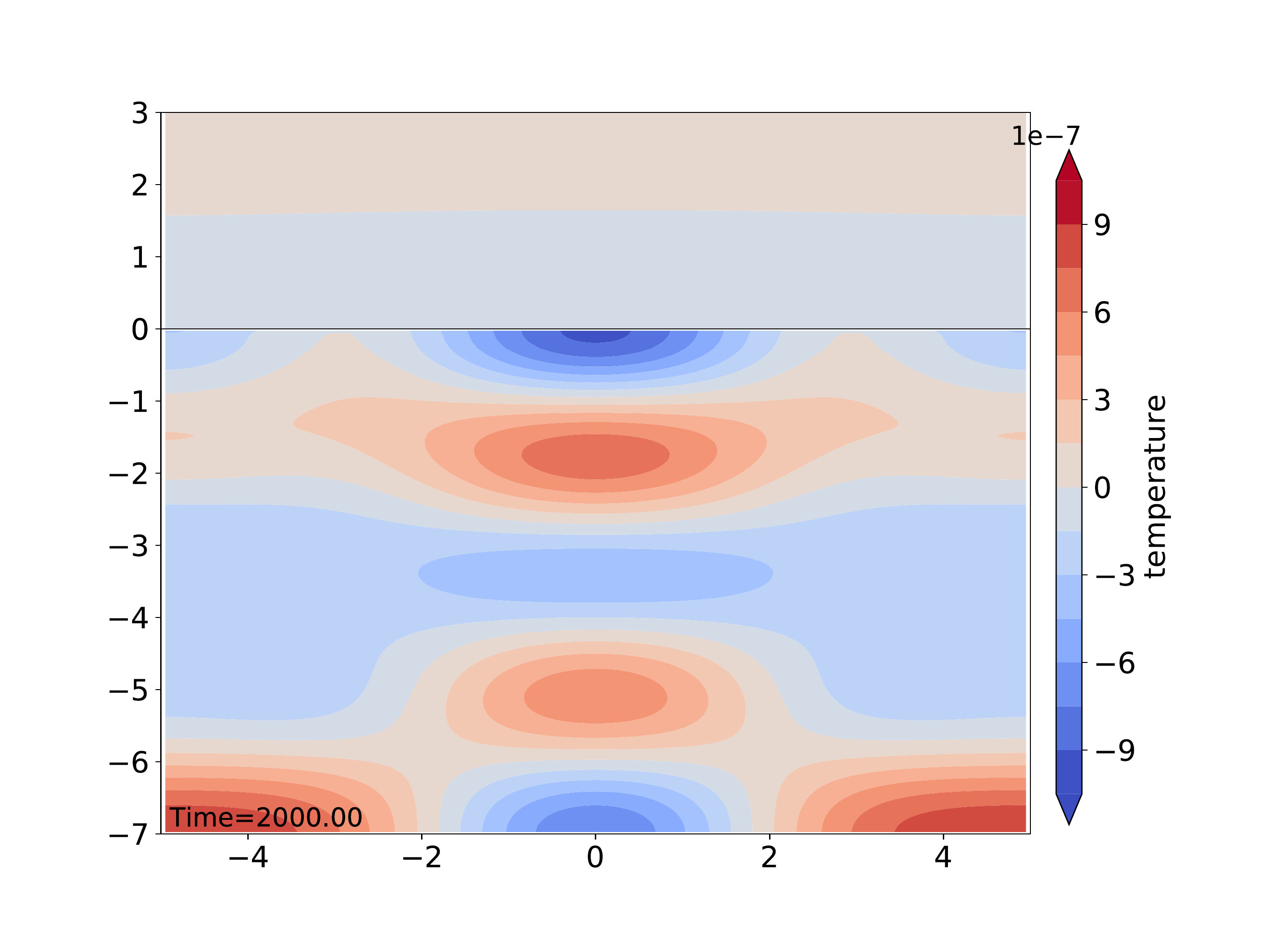}  
        \caption{MPRK2 - RK45}
        \figlab{rtbdc-ss-mprk2}
    \end{subfigure}%
    \begin{subfigure}{0.43\textwidth}
    \centering
      \includegraphics[trim=2.5cm 1.5cm 2.2cm 2.2cm,
      clip=true,
        width=0.95\columnwidth]{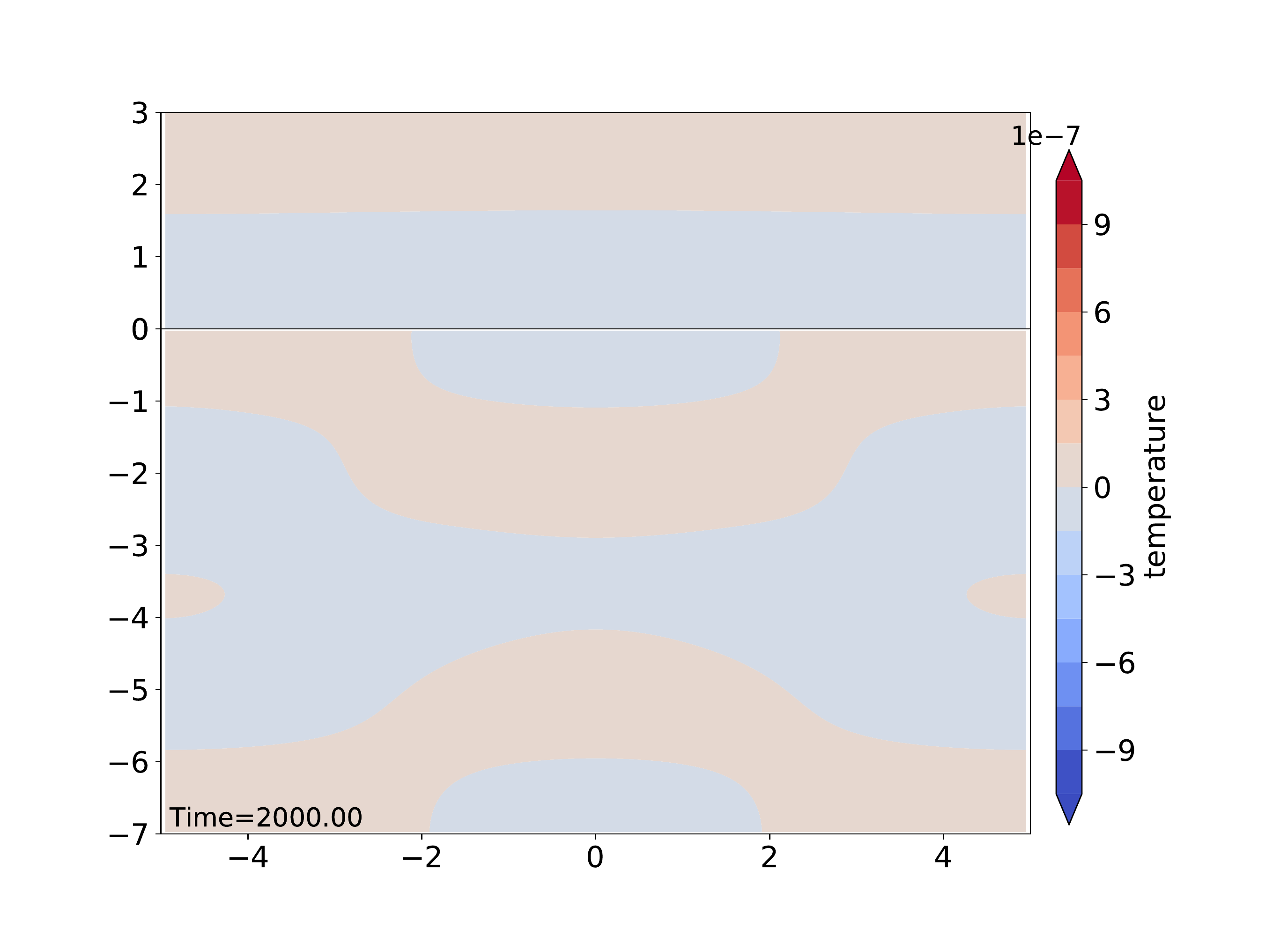}  
        \caption{RK2 - RK45}
        \figlab{rtbdc-ss-rk2}
    \end{subfigure}%
    \caption{Temperature difference of (a) MPRK2 ($m=4$) and (b) RK2 for the convection example at $t=2000$ with respect to the RK45 solution of $\dtt_{RK45}=0.0125$. 
    We take $\dtt_{MPRK2} =0.025$ and $\dtt_{RK2}=0.00625$. 
    } 
   \figlab{rtbdc-ss-difference}
\end{figure}

\begin{table}[t] 
  \caption{Errors of MPRK2 ($m=4$) and RK2 for the thermal convection example at $t=2000$ with respect to the RK45 solution of $\dtt_{RK45}=0.0125$. $Cr_1$ and $Cr_2$ are the Courant numbers on $\Omega_1$ and $\Omega_2$, respectively; wc is the wall clock; and the subscript $r$ indicates the reference solution.
  }
  \tablab{rtbdc-relerr} 
  \begin{center} 
    \begin{tabular}{*{1}{c}|*{1}{c}|*{1}{c}|*{1}{c}|*{1}{c}|*{1}{c}}  
    \hline 
    \multirow{2}{*}{ }
    & \multirow{2}{*}{$\dtt  (Cr_1,Cr_2)$}
    & \multirow{2}{*}{$ \left\Vert \rho  - \rho_r \right\Vert  $}  
    & \multirow{2}{*}{$ \left\Vert \rho {\bf u} - \rho {\bf u}_r \right\Vert  $}  
    & \multirow{2}{*}{$ \left\Vert \rho E - \rho E_r \right\Vert  $}  
    & \multirow{2}{*}{wc[s]} \tabularnewline 
    & && & &  \tabularnewline 
    \hline\hline 
    RK2 & 0.00625 (0.13,0.51) &       8.02E-07&       4.96E-07&       2.00E-06& 12946\tabularnewline
    MPRK2 & 0.025 (0.51,2.0)&       6.44E-06&       4.31E-06&       1.63E-05& 8028\tabularnewline
    \hline\hline 
    \end{tabular} 
  \end{center} 
\end{table}

\subsection{Parallel Performance of the MPRK method}

Now we study the parallel performance in terms of weak and strong scaling for three-dimensional coupled Navier--Stokes equations by using wind-driven flow and thermal convection examples. 
Figure \figref{evolution-tb3d} shows the evolution of the temperature field for the thermal convection. 
We conduct a simulation with the MPRK2 ($m=4$) method over a mesh of $100\times100\times200$ elements on $\Omega_1=(-5,5)\times(-5,5)\times(-16,0)$ and $100\times100\times100$ elements on $\Omega_2=(-5,5)\times(-5,5)\times(0,2)$ for $t\in\LRs{0,300}$. Here we take $\delta\theta_1=7.5(1+\cos(\pi r_1)$ for $r_1 \le 2.5$, $\delta\theta_1=0$ for $r_1>2.5$, $\delta\theta_2=-7.5(1+\cos(\pi r_2)$ for $r_2 \le 2.5$, and $\delta\theta_2=0$ for $r_2>2.5$. 


\begin{figure}[h!t!b!]
  \centering
      \includegraphics[trim=25cm 1.5cm 25cm 6cm,
      clip=true,
        width=0.45\columnwidth]{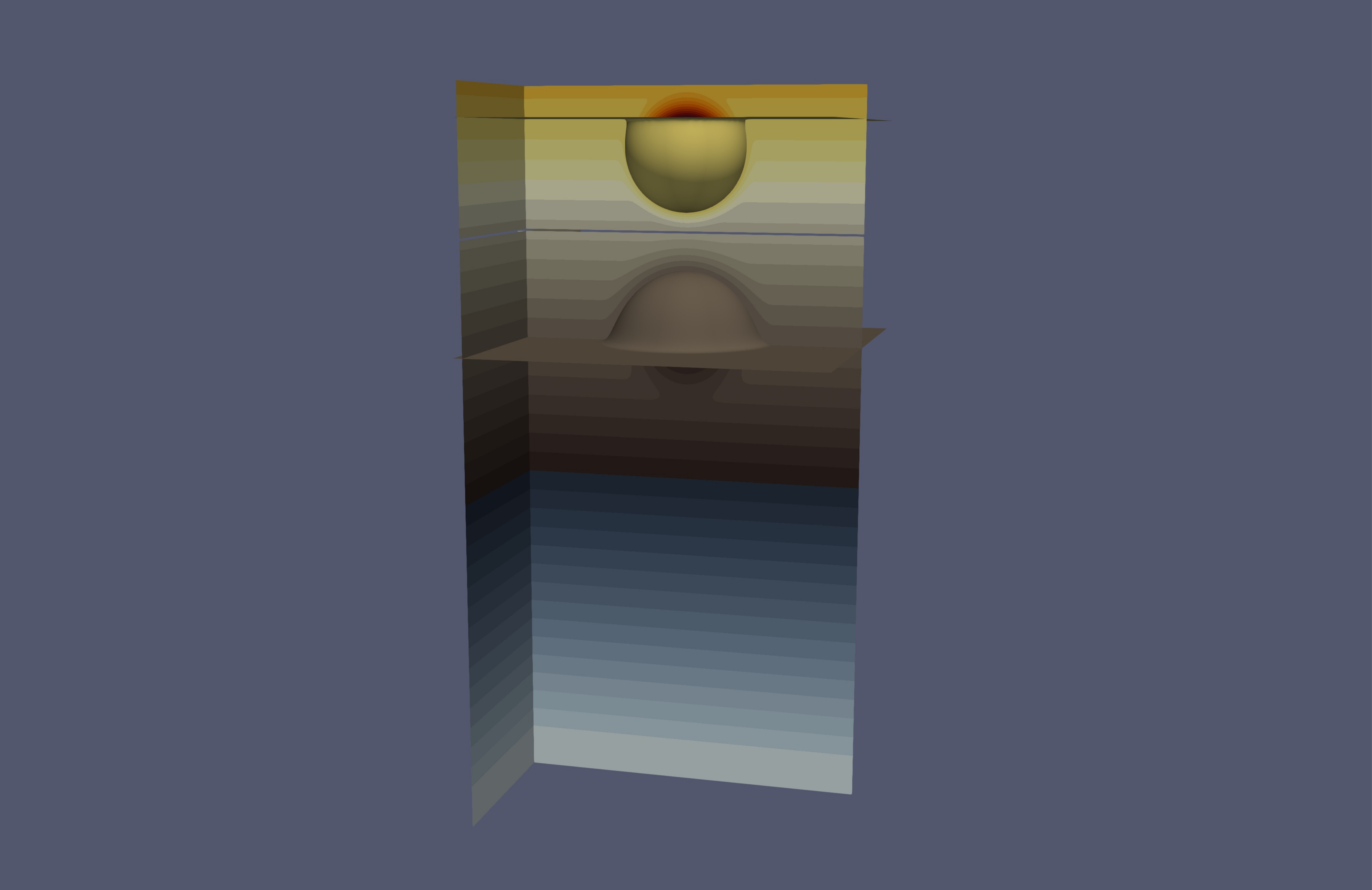}  
     \includegraphics[trim=25cm 1.5cm 25cm 6cm,
      clip=true,
        width=0.45\columnwidth]{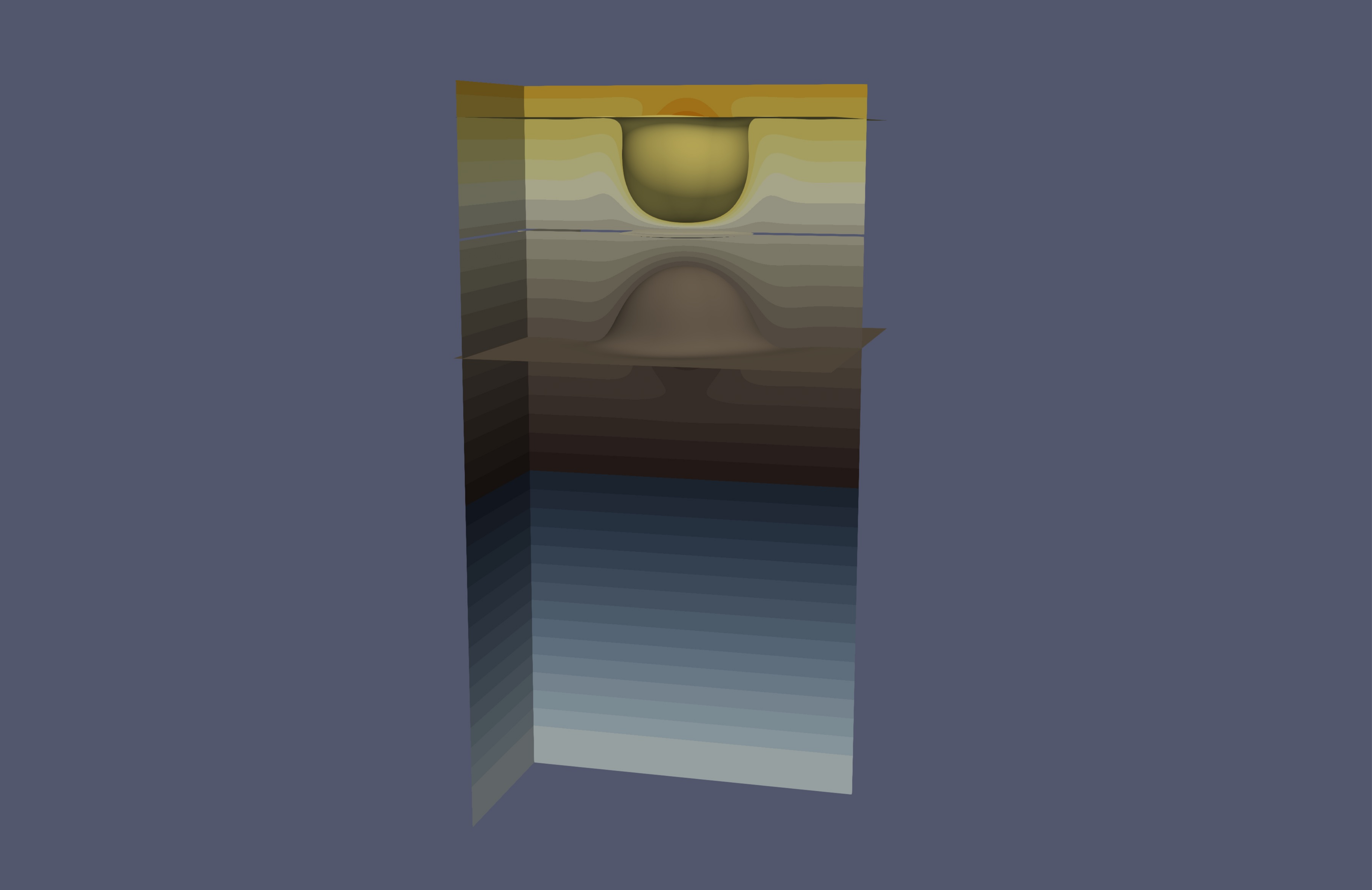} 
    \includegraphics[trim=25cm 1.5cm 25cm 6cm,
      clip=true,
        width=0.45\columnwidth]{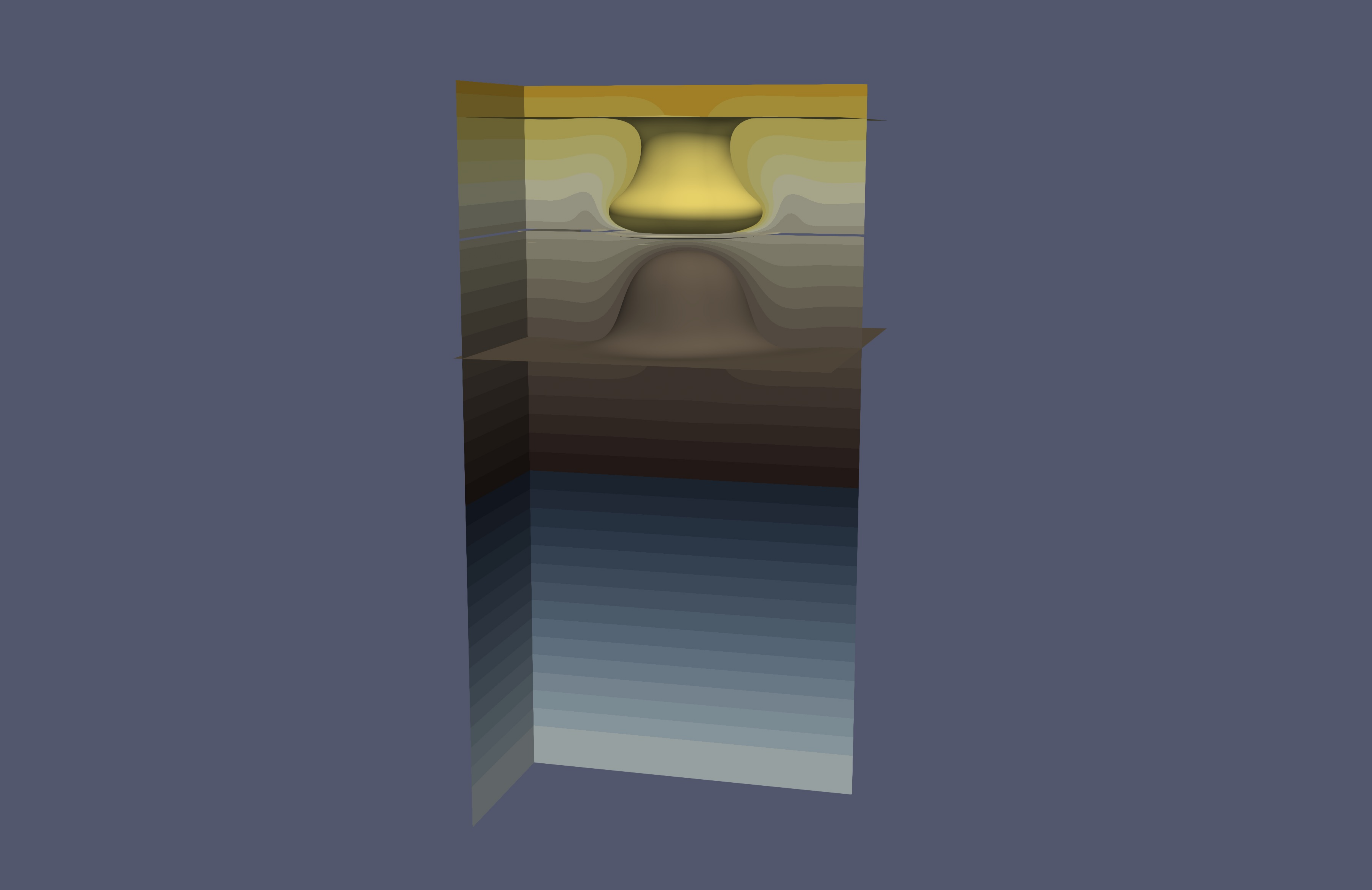}  
    \includegraphics[trim=25cm 1.5cm 25cm 6cm,
      clip=true,
        width=0.45\columnwidth]{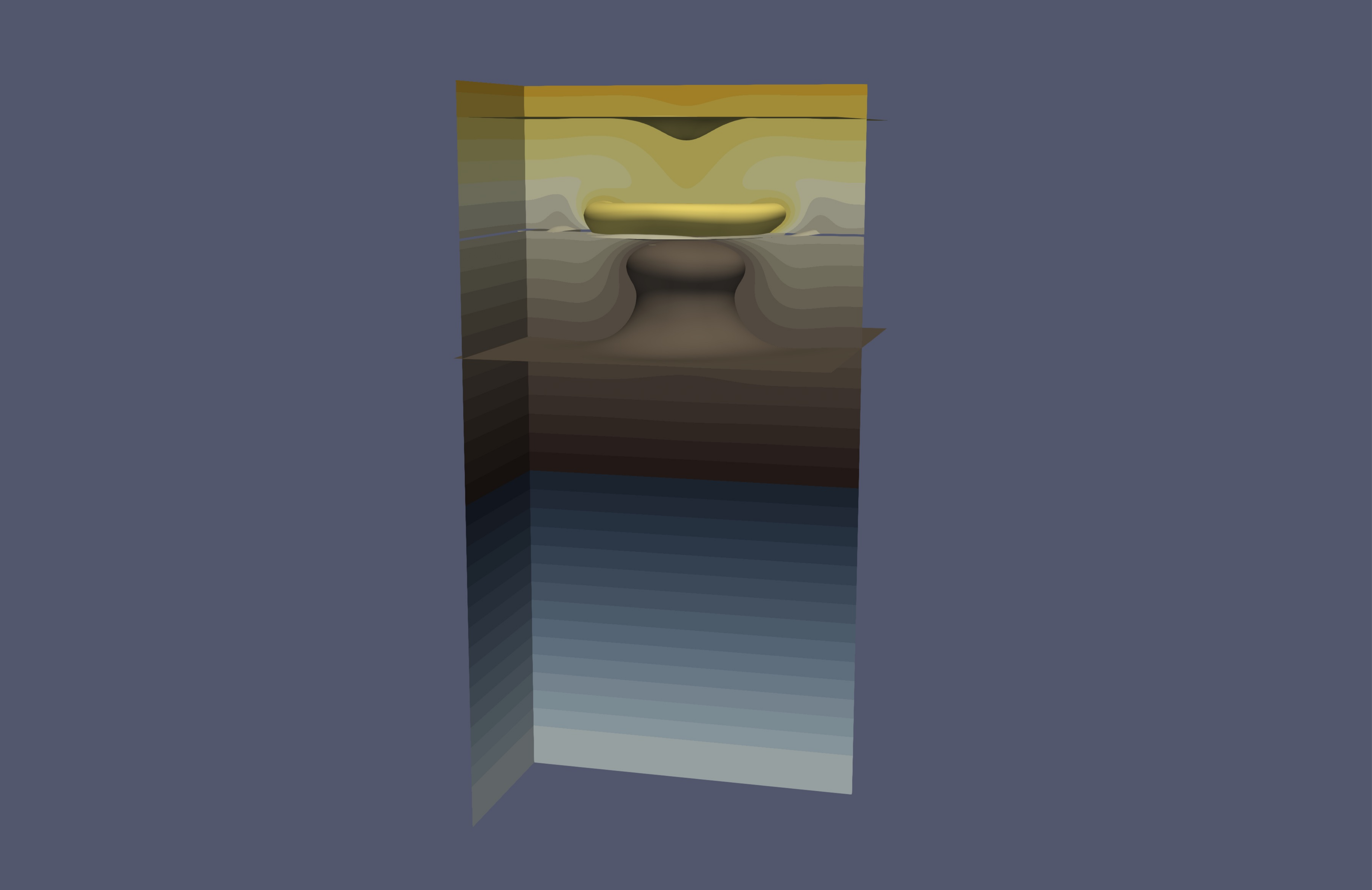}  
      \caption{Evolution of the temperature field for the thermal convection: 
  simulation is conducted with the MPRK2 ($m=4$) method over 
 a mesh of $100\times100\times100$ elements on $\Omega_1=(-5,5)\times(-5,5)\times(0,2)$ 
  and $100\times100\times200$ elements on $\Omega_2=(-5,5)\times(-5,5)\times(-16,0)$ for $t\in\LRc{0,100,200,300}$. 
  }
      
   \figlab{evolution-tb3d}
\end{figure}
 
The parallel simulations are conducted on the Cray XC40 Theta with Intel Knights Landing (KNL) CPUs at the Argonne Leadership Computing Facility. Each KNL compute node is composed of a single Xeon Phi 7230 1.3 GHz processor with 64 cores, 16 GB Multi-Channel DRAM, and 192 GB DDR4 memory. The interconnect topology is a dual place Dragonfly with ten groups. Each group consists of two cabinets or racks. The total bisection bandwidth is 7.2 TB/s.

\subsubsection*{Strong scaling} 

For a given simulation setup, the execution time is expected to decrease up to a certain point as the number of cores increases. 
In Figure \figref{strong-scaling} we present strong-scaling results with the two examples.
The measured speedup factors of MPRK2 over RK2 are also annotated.

For the wind-driven flow example,  
we take the mesh of $512\times512\times 1024$ elements on $\Omega_1$ and $512\times512\times 128$ on $\Omega_2$ with $\dx=\dy=0.1$,  $\dz_1=0.1$, and $\dz_2=0.025$. 
We run the simulations with $\dtt_{MPRK2}=0.05$ for MPRK2 ($m=4$) and 
$\dtt_{RK2}=0.0125$ for RK2 while increasing the number of cores from $1,024$ to $16,384$ for $t\in[0,0.25]$.  
The results are summarized in Figure \figref{wdf-strong-scaling}. 
Up to 8,192 cores, MPRK2 shows good strong scalability, after which the scalability slightly degrades. 
This is expected from the parallel speedup estimation in \eqnref{par_speedup}. 

For the thermal bubble example, 
we take the mesh of $512\times512\times 1024$ elements on $\Omega_1=(-6.4,6.4)\times(-6.4,6.4)\times(-16,0)$ and $512\times512\times 256$ on $\Omega_2=(-6.4,6.4)\times(-6.4,6.4)\times(0,2)$. 
This corresponds to $\dx=\dy=0.025$,  $\dz_1=0.015625$, and $\dz_2=0.00390625$.
We ran the simulations with $\dtt_{MPRK2}=0.00625 $ for MPRK2 ($m=4$) and with
$\dtt_{RK2}=0.0015625$ for RK2 while increasing the number of cores from $1,024$ to $16,384$ for $t\in[0,0.03125]$.  
MPRK2 shows favorable strong-scaling behavior in Figure \figref{tb-strong-scaling}.

In both cases, MPRK2 is faster than the RK2 counterpart. 
This is expected because the slow part dominates the computational cost. 
The estimated speedup factors of MPRK2 are 2.96 for wind-driven flow and 1.98 for thermal bubble examples according to \eqnref{speedup} with 6 buffer elements.
The measured speedup factors are bounded above the estimated factors. 
For example, the speedup factors are $2.1$ with $1,024$ cores in Figure \figref{wdf-strong-scaling}, 
and $1.6$ with $1,024$ cores in Figure \figref{tb-strong-scaling}. 
Also, we observe that the speedup factors tend to decrease toward 1 with increasing number of cores, which agrees with the speedup estimation of \eqnref{par_speedup} in the parallel case. 


\begin{figure}[h!t!b!]
  \centering
    \begin{subfigure}{0.43\textwidth}
    \centering
      \includegraphics[trim=1.0cm 0cm 3.5cm 3.0cm,clip=true, width=0.95\columnwidth]{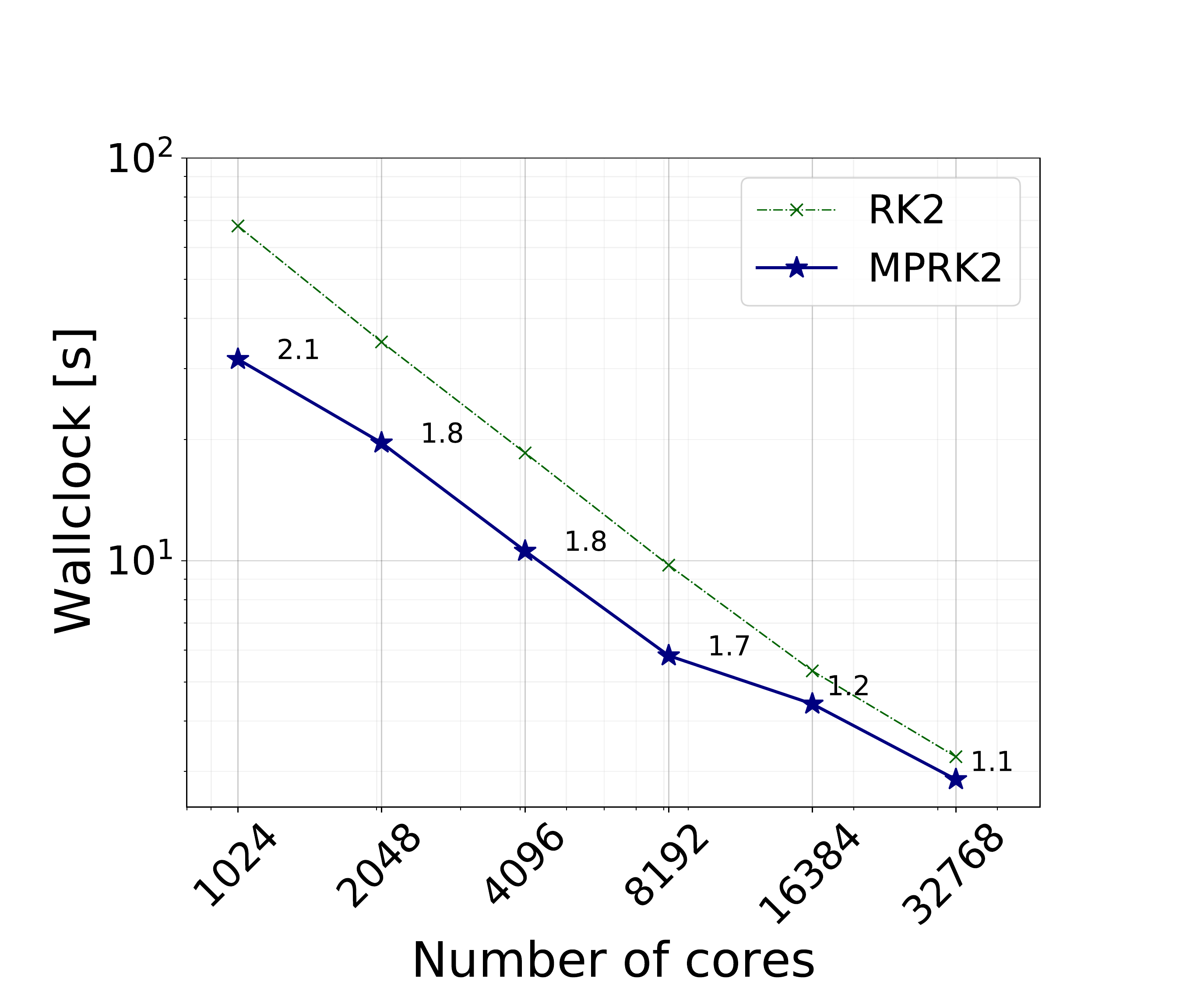}
        \caption{Wind-driven flow}
        \figlab{wdf-strong-scaling}
    \end{subfigure}%
    \begin{subfigure}{0.43\textwidth}
    \centering
      \includegraphics[trim=1.0cm 0cm 3.5cm 3.0cm,
      clip=true,
        width=0.95\columnwidth]{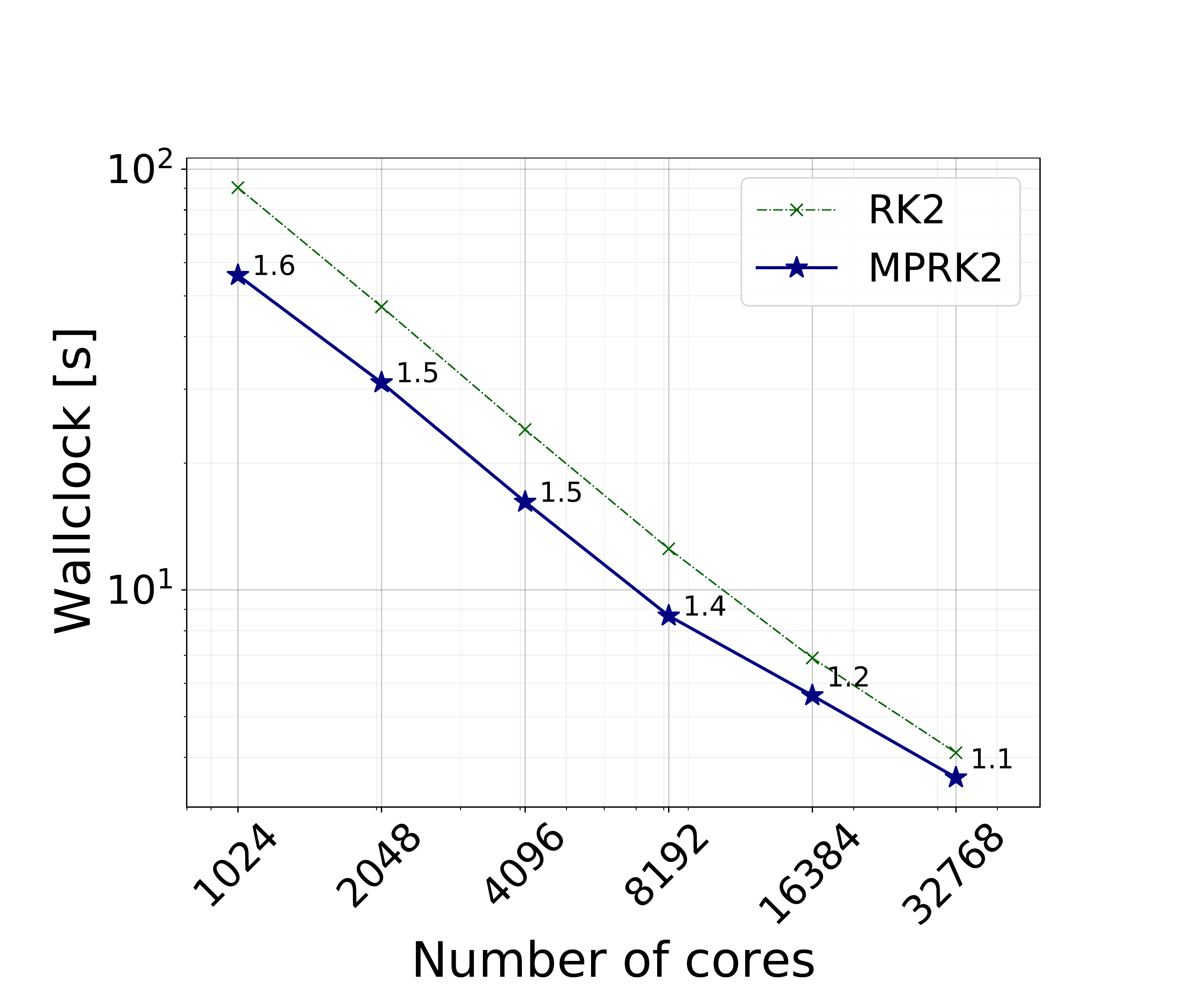}  
        \caption{Thermal bubble}
        \figlab{tb-strong-scaling}
    \end{subfigure}%
    \caption{Strong-scaling studies for three-dimensional (a) wind-driven flows and (b) thermal bubble examples.
    }
   \figlab{strong-scaling}
\end{figure}


\subsubsection*{Weak scaling}

For weak scaling we assign the same amount of work to each processor while increasing the number of processors. 
Ideally, the execution time will remain the same.
In practice, however,
the execution time tends to increase due to communication overhead. 
Figure \figref{weak-scaling}, Table \tabref{wdf-weak-scaling}, and Table \tabref{tb-weak-scaling} show two weak-scaling studies: wind-driven flows with 
9,216 elements per core in Figure \figref{wdf-weak-scaling} and thermal convection with 12,288 elements per core in Figure \figref{tb-weak-scaling}. For the former, we take $\dtt_{MPRK2}=0.05$ and $\dtt_{RK2}=0.0125$ and measure the wall-clock times at $T=0.25$. For the latter, we take $\dtt_{MPRK2}=0.00625$ and $\dtt_{RK2}=0.0015625$ and measure the wall-clock times at $T=0.03125$.
The number of processors is chosen in the set $\LRc{8,64,512,4096,32768}$, which is annotated in Figure \figref{weak-scaling}.
We observe that MPRK2 ($m=4$) shows favorable weak-scaling results and its time-to-solution is faster than that of its RK2 counterpart.

\begin{table}[t] 
  \caption{Weak-scaling test for wind-driven flow example: 
  wall clocks of MPRK2 and RK2 are summarized. 
  }
  \tablab{wdf-weak-scaling} 
  \begin{center} 
  \begin{tabular}{*{1}{c}|*{1}{c}|*{1}{c}|*{1}{c}|*{1}{c}} 

    \hline 
    
    $N_p$ & $N_E (\texttt{nxe,nye,nze1,nze2})$ & MPRK2 [s] & RK2 [s] & $\texttt{wcr}(\texttt{spd}$) \tabularnewline 
    
    \hline\hline 

           8
    & (      32,      32,      64,       8)
    &      1.5 &      2.1 &     1.37(2.53)
    \tabularnewline
    
          64
    & (      64,      64,     128,      16)
    &      1.9 &      2.7 &     1.38(2.74)
    \tabularnewline
    
         512
    & (     128,     128,     256,      32)
    &      2.1 &      2.8 &     1.36(2.87)
    \tabularnewline
    
        4096
    & (     256,     256,     512,      64)
    &      2.3 &      2.8 &     1.20(2.93)
    \tabularnewline
    
       32768
    & (     512,     512,    1024,     128)
    &      2.8 &      3.3 &     1.18(2.97)
    \tabularnewline
    
    \hline\hline 
    
    \end{tabular} 

  \end{center} 
\end{table}

\begin{table}[t] 
  \caption{Weak-scaling test for thermal bubble example: 
  wall clocks of MPRK2 and RK2 are summarized. 
  }
  \tablab{tb-weak-scaling} 
  \begin{center} 
  \begin{tabular}{*{1}{c}|*{1}{c}|*{1}{c}|*{1}{c}|*{1}{c}} 

    \hline 
    
    $N_p$ & $N_E (\texttt{nxe,nye,nze1,nze2})$ & MPRK2 [s] & RK2 [s] & $\texttt{wcr}(\texttt{spd}$)  \tabularnewline 
    
    \hline\hline 
    
           8
    & (      32,      32,      64,      32)
    &      2.1 &      2.7 &     1.26(1.83)
    \tabularnewline
    
          64
    & (      64,      64,     128,      64)
    &      2.6 &      3.4 &     1.33(1.91)
    \tabularnewline
    
         512
    & (     128,     128,     256,     128)
    &      2.7 &      3.6 &     1.32(1.95)
    \tabularnewline
    
        4096
    & (     256,     256,     512,     256)
    &      3.0 &      3.7 &     1.23(1.98)
    \tabularnewline
    
       32768
    & (     512,     512,    1024,     512)
    &      3.5 &      3.9 &     1.10(1.99)
    \tabularnewline
    
    \hline\hline 
    
    \end{tabular} 

  \end{center} 
\end{table}


\begin{figure}[h!t!b!]
  \centering
    \begin{subfigure}{0.43\textwidth}
    \centering
      \includegraphics[trim=0.0cm 0.0cm 1.5cm 1.8cm,clip=true, width=0.95\columnwidth]{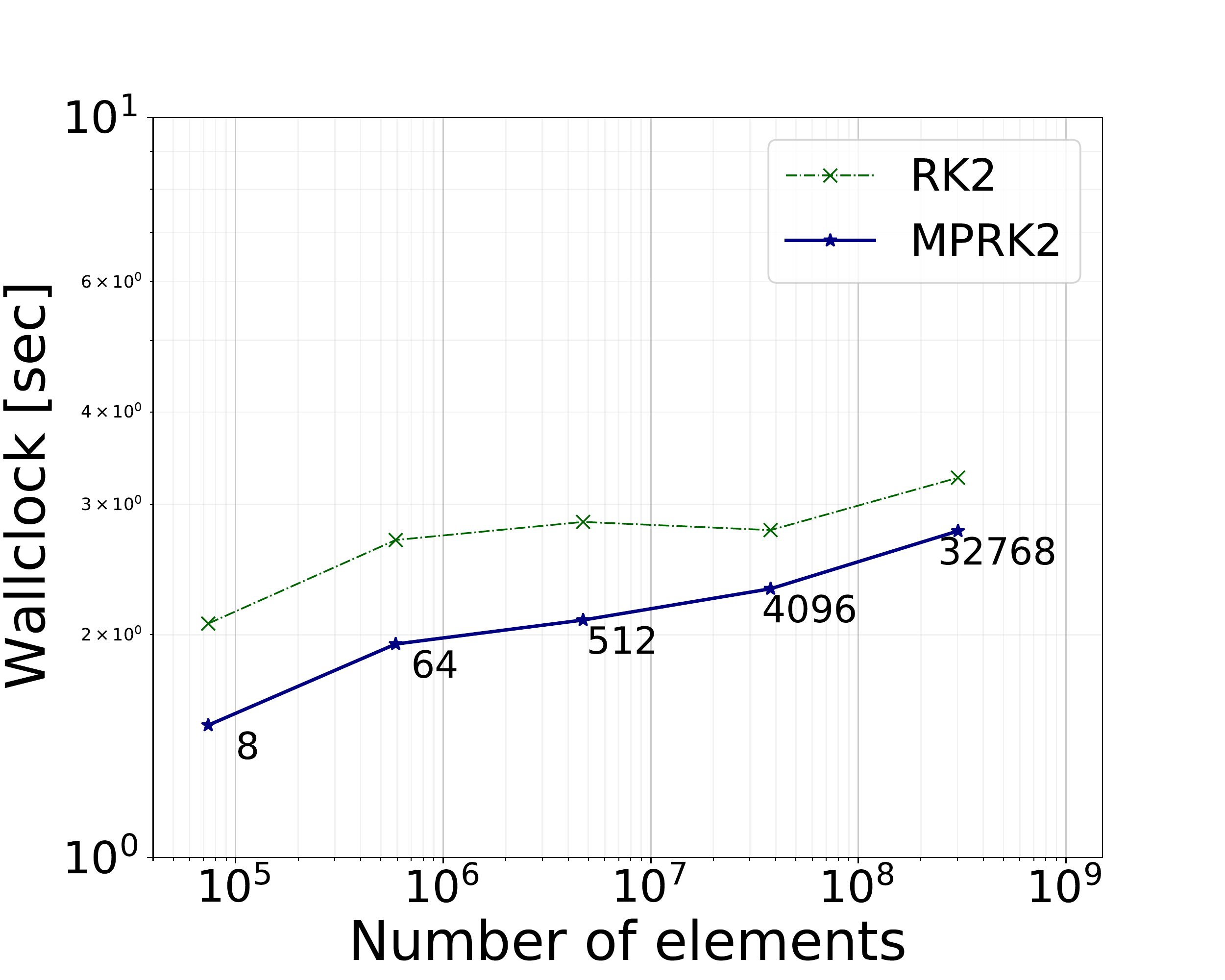}
        \caption{Wind-driven flow}
        \figlab{wdf-weak-scaling}
    \end{subfigure}%
    \begin{subfigure}{0.43\textwidth}
    \centering
      \includegraphics[trim=0.0cm 0.0cm 1.5cm 1.8cm,
      clip=true,
        width=0.95\columnwidth]{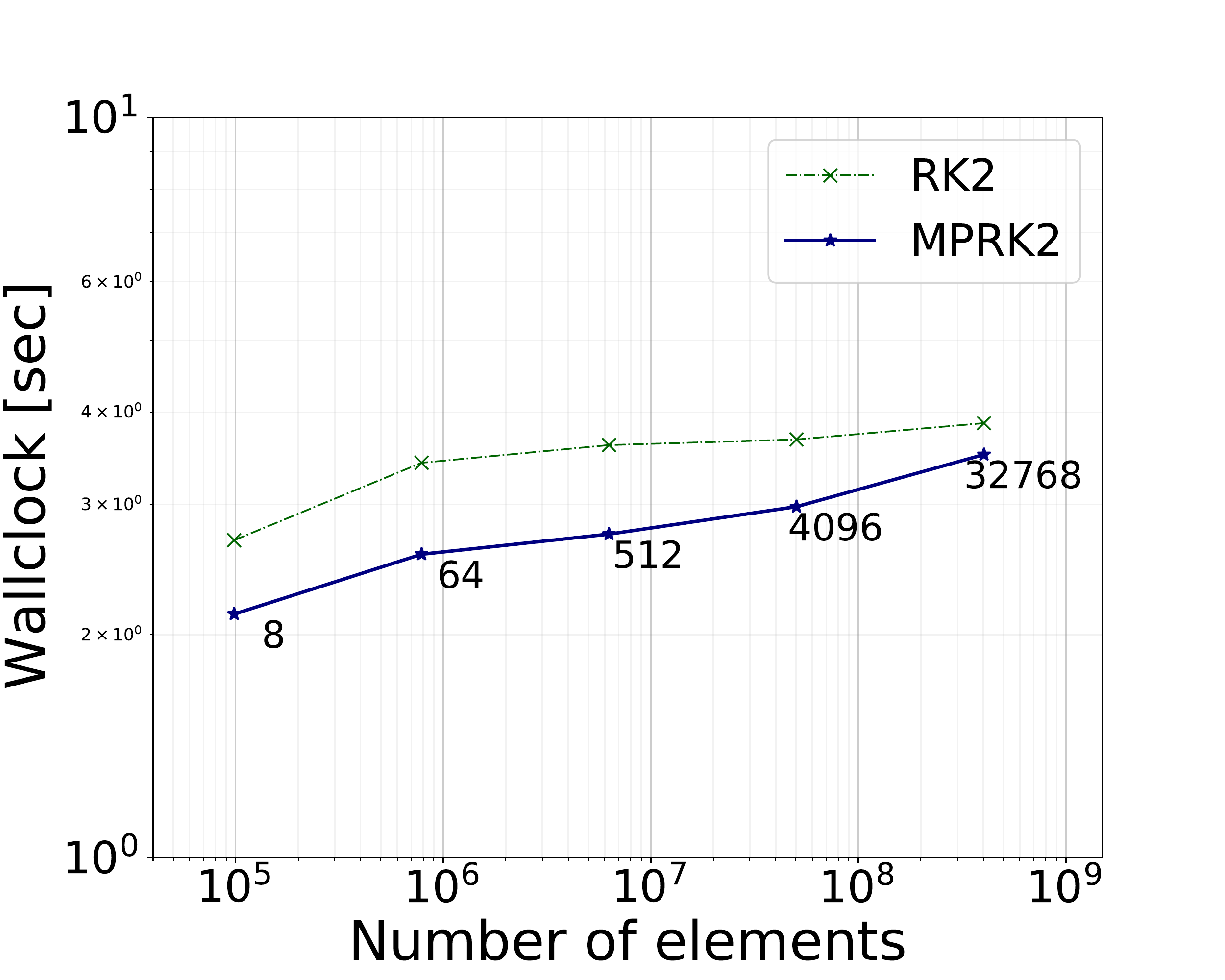}
        \caption{Thermal bubble}
        \figlab{tb-weak-scaling}
    \end{subfigure}%
    \caption{Weak-scaling studies for three-dimensional (a) wind-driven flow and (b) thermal convection examples.
    }
   \figlab{weak-scaling}
\end{figure}

\section{Conclusions}
\seclab{Conclusion}
 
 In this study, we have applied multirate partitioned Runge--Kutta methods to a fluid-fluid interaction problem and demonstrated its parallel performance by using the PETSc library.
 We considered the coupled compressible Navier--Stokes systems 
with gravity and the rigid-lid coupling condition arising from the atmosphere and the ocean interaction. 
Multirate partitioned Runge--Kutta (MPRK) coupling methods explicitly solve both the ocean and the atmospheric models with different step sizes, have a potential to be scalable on modern computing architectures. 
Moreover, by placing a buffer between the sub-models, MPRK methods naturally handle the coupling interface. 

The numerical examples confirm that MPRK2 methods conserve the total mass and have the second-order rate of convergence in time. 
In the Kelvin--Helmholtz instability example, we observe that the total mass loss for MPRK2 is less than $\mathcal{O}(10^{-13})$ regardless of the temporal rate $m$. That is, 
the total mass is conserved as expected. 
In the thermal convection examples, 
we observe the second-order convergence rate of the MPRK2 coupling method.

A theoretical speedup analysis shows that the MPRK2 method has an advantage over its base method in terms of computational cost when a larger number of elements is located in the slow regions, which are time-stepped with the temporal factor $m$ larger than in the fast regions. 
The theoretical estimates are verified numerically by using both a single core and 32 cores and computing the speedup of the MPRK2 method over its  base (single rate).
We also studied the parallel performance of the MPRK2 method using three-dimensional coupled compressible Navier--Stokes equations. 
Thanks to its explicit nature, the MPRK2 coupling method shows favorable strong- and weak-scaling results for the three-dimensional wind-driven flow and the thermal convection examples. 

In our current work, we only consider matching grids at the interface. In general, however, ocean and atmospheric model have non-matching grid at the interface. To handle this, the nonconforming mortar element methods can be employed. We will address the issue in our future work. In addition, coupling incompressible and compressible Navier--Stokes system, and high-order multirate methods will be considered as well.


\section*{Acknowledgments}
This material is based upon work supported by the U.S. Department of Energy, Office of Science, Office of Advanced Scientific Computing Research and Office of Biological and Environmental Research, Scientific Discovery through Advanced Computing (SciDAC) program under Contract DE-AC02-06CH11357 through the Coupling Approaches for Next-Generation Architectures (CANGA) Project and the FASTMath institute. We also gratefully acknowledge the use of Theta in the resources of the Argonne Leadership Computing Facility, which is a DOE Office of Science User Facility supported under Contract DE-AC02-06CH11357.

\appendix

\section{Nondimensionalization}

With the nondimensional variables, 
$ 
\rho^*_\model = \frac{\rho_\model}{\rho_r}
$,
$
\pres^*_\model = \frac{\pres_\model}{\rho_r u_r^2}
$,
$
\ub^*_\model = \frac{\ub_\model}{u_r}
$,
$
x^* = \frac{x}{L}
$,
$
t^* = \frac{t}{L/u_r}
$,
$
\mu^*_\model = \frac{\mu_\model}{\mu_r}
$,
$
{\bf g}^* = \frac{\bf g}{u_r^2/L}
$,
 and 
$
T^*_\model = \frac{T_\model}{T_r},
$
\noindent
we rewrite the governing equation \eqnref{cns-gov} as
\begin{subequations}
\eqnlab{cns-gov-nondimension}
\begin{align}
  \dd{\rho_\model^*}{t^*}     + \Grad^*\cdot \LRp{\rho \ub_\model^*} &= 0,\\
  \dd{\rho \ub_\model^*}{t^*} + \Grad^*\cdot \LRp{\rho \ub_\model^* \otimes \ub_\model^* + \Ical \pres_\model^*} &= \Grad^* \cdot \sigma_\model^* - \rho_\model^* {\bf g}^*,\\
  \dd{\rho E_\model^*}{t^*}   + \Grad^*\cdot \LRp{\rho \ub_\model^* H_\model^*} &= \Grad^*\cdot (\sigma_\model^* \ub_\model^* - \Pi_\model^* ) - \rho_\model^*{\bf g}^* \cdot \ub^*_\model,
\end{align}
\end{subequations}
where 
$\sigma_\model^* = \tilde{\mu}_\model \LRp{\Grad^* \ub^*_\model + \Grad^* (\ub^*_\model)^\top -\frac{2}{3}\Ical \Grad^* \cdot \ub^*_\model } $,  
$\Pi^*_\model = \frac{\tilde{c}_p\tilde{\mu}_\model}{Pr} \Grad^* T^*_\model $,  $\Grad^*=\frac{1}{L}\Grad$,
$\tilde{\mu}_\model := \frac{\mu^*_\model}{Re_r}$, 
$Re_r := \frac{\rho_r u_r L}{\mu_r}$,
and $\tilde{c}_p := \frac{T_r c_p}{u_r^2}$.
We choose $\rho_r = \rho_{\infty_2}$, $\mu_r=\mu_{\infty_2}$, $T_r=T_{\infty_2}$, and $u_r = a_{\infty_2} = \sqrt{\gamma R T_{\infty_2}}$.
This leads to $\tilde{c}_p = \frac{1}{\gamma-1}$ 
  and $ Re_r =  \frac{Re_{\infty_2}}{M_{\infty_2}}$.
The normalized equation of state for an ideal gas is $\pres^*_\model = \gamma^{-1} \rho^*_\model T^*_\model = \rho^*_\model e^*_\model(\gamma-1)$.

\section{FV discretization in a three-dimensional ($x,y,z$) uniform grid}
\seclab{FV3}

With a three-dimensional uniform mesh, \eqnref{cns-fv} for the $\ell=(i,j,k)$ element becomes 
\begin{multline}
  \DD{\overline{\qb}_{(i,j,k)}}{t} - \overline{\sb}_{(i,j,k)}
  =  
   - \frac{1}{\dx}\LRp{\Fbs_{(\iphalf,j,k)}  - \Fbs_{(\imhalf,j,k)}} 
   \\
   - \frac{1}{\dy}\LRp{\Gbs_{(i,\jphalf,k)}  - \Gbs_{(i,\jmhalf,k)}} 
   - \frac{1}{\dz}\LRp{\Hbs_{(i,j,\kphalf)}  - \Hbs_{(i,j,\kmhalf)}},
  \eqnlab{ns-fv-cartesian-3d}   
\end{multline}
where 
$\Fbs = \Fb^{I^*} - \Fb^{V^*}$,
$\Gbs = \Gb^{I^*} - \Gb^{V^*}$,
$\Hbs = \Hb^{I^*} - \Hb^{V^*}$,
$$
\Fb^{I^*} =
\begin{pmatrix}
\rho u \\ 
\rho u u + \pres \\ 
\rho u v \\ 
\rho u w \\ 
\rho u H
\end{pmatrix}^*,
\quad \Gb^{I^*} = 
\begin{pmatrix}
\rho v \\ 
\rho v u  \\ 
\rho v v + \pres \\ 
\rho v w \\ 
\rho v H
\end{pmatrix}^*,
\quad \Hb^{I^*} = 
\begin{pmatrix}
\rho w \\ 
\rho w u  \\ 
\rho w v  \\ 
\rho w w + \pres \\ 
\rho w H
\end{pmatrix}^*,
$$

$$
\Fb^{V^*} =
\begin{pmatrix}
0 \\ 
\sigma_{xx} \\ 
\sigma_{xy} \\ 
\sigma_{xz} \\ 
\ub\boldsymbol{\sigma}_x
- \Pi_x
\end{pmatrix}^*,
\quad \Gb^{V^*} =
\begin{pmatrix}
0 \\ 
\sigma_{yx} \\ 
\sigma_{yy} \\ 
\sigma_{yz} \\ 
\ub\boldsymbol{\sigma}_y - \Pi_y
\end{pmatrix}^*,
\quad \Hb^{V^*} =
\begin{pmatrix}
0 \\ 
\sigma_{zx} \\ 
\sigma_{zy} \\ 
\sigma_{zz} \\ 
\ub\boldsymbol{\sigma}_z - \Pi_z
\end{pmatrix}^*,
$$
$\ub\boldsymbol{\sigma}_x=u\sigma_{xx} + v\sigma_{yx} + w\sigma_{zx}$,
$\ub\boldsymbol{\sigma}_y=u\sigma_{xy} + v\sigma_{yy} + w\sigma_{zy}$, and 
$\ub\boldsymbol{\sigma}_z=u\sigma_{xz} + v\sigma_{yz} + w\sigma_{zz}$.

For the cell-centered second-order FV, 
we compute the numerical fluxes in \eqnref{ns-fv-cartesian-3d} using linearly reconstructed variables. 
Since each element has only cell-averaged values, 
the reconstruction requires information from adjacent elements. 
For the inviscid flux ($\Fb^{I^*}, \Gb^{I^*} \text{, and }\Hb^{I^*}$), 
we compute the cell-centered gradients $\overline{\Grad \qb}$ using the least-square (LS) scheme \cite{syrakos2017critical}.\footnote{
  In uniform mesh, LS schemes approximate $x$, $y$, and $z$ directional gradients by 
  \begin{align*}
    \overline{\dd{\qb}{x}}_{(i,j,k)} & \approx \frac{\overline{\qb}_{(i+1,j,k)} - \overline{\qb}_{(i-1,j,k)}  }{2\dx} + \mathcal{O}((\dx)^2),\\ 
    \overline{\dd{\qb}{y}}_{(i,j,k)} & \approx \frac{\overline{\qb}_{(i,j+1,k)} - \overline{\qb}_{(i,j+1,k)}  }{2\dy} + \mathcal{O}((\dy)^2), \\
    \overline{\dd{\qb}{z}}_{(i,j,k)} & \approx \frac{\overline{\qb}_{(i,j,k+1)} - \overline{\qb}_{(i,j,k-1)}  }{2\dz} + \mathcal{O}((\dz)^2).
  \end{align*}
} 
With the gradients, the left and the right conserved variables ($\qb^l$ and $\qb^r$) 
are obtained by 
\begin{align*}
  \qb_{(i+\half,j,k)}^l &= \overline{\qb}_{(i,j,k)}   + \overline{\dd{\qb}{x}}_{(i,j,k)} \frac{\dx}{2}, \quad
  \qb_{(i+\half,j,k)}^r = \overline{\qb}_{(i+1,j,k)} - \overline{\dd{\qb}{x}}_{(i+1,j,k)} \frac{\dx}{2}  
\end{align*}
at the $x$-face of $(i+\half,j,k)$, 
\begin{align*}
  \qb_{(i,j+\half,k)}^l &= \overline{\qb}_{(i,j,k)}   + \overline{\dd{\qb}{y}}_{(i,j,k)} \frac{\dy}{2}, \quad
  \qb_{(i,j+\half,k)}^r = \overline{\qb}_{(i,j+1,k)} - \overline{\dd{\qb}{y}}_{(i,j+1,k)} \frac{\dy}{2} 
\end{align*}
at the $y$-face of $(i,j+\half,k)$, 
and 
\begin{align*}
  \qb_{(i,j,k+\half)}^l &= \overline{\qb}_{(i,j,k)}   + \overline{\dd{\qb}{z}}_{(i,j,k)} \frac{\dz}{2}, \quad
  \qb_{(i,j,k+\half)}^r = \overline{\qb}_{(i,j,k+1)} - \overline{\dd{\qb}{z}}_{(i,j,k+1)} \frac{\dz}{2} 
\end{align*}
at the $z$-face of $(i,j,k+\half)$. 

For the viscous flux ($\Fb^{V^*}, \Gb^{V^*}\text{, and }\Hb^{V^*}$), 
we compute the cell-centered gradients of velocity and temperature ($\overline{\Grad \ub}$ and $\overline{\Grad T}$) by using the LS scheme, and  we compute
the common velocity of $\hat{\ub}$,
\begin{align*}
  \hat{\ub}_{(i+\half,j,k)} = \half\LRp{ \overline{\ub}_{(i,j,k)} + \overline{\ub}_{(i+1,j,k)} }, 
  \hat{\ub}_{(i,j+\half,k)} = \half\LRp{ \overline{\ub}_{(i,j,k)} + \overline{\ub}_{(i,j+1,k)} }, \text{and},
  \hat{\ub}_{(i,j,k+\half)} = \half\LRp{ \overline{\ub}_{(i,j,k)} + \overline{\ub}_{(i,j,k+1)} },
\end{align*}
by taking an arithmetic average. 
As for the common gradient of $\widehat{\Grad T}$, 
we obtain normal gradients by finite difference approximation
\begin{align*}
  \widehat{\dd{T}{x}}_{(i+\half,j,k)} &= \frac{1}{\dx}\LRp{ \overline{T}_{(i+1,j,k)} - \overline{T}_{(i,j,k)} },\\
  \widehat{\dd{T}{y}}_{(i,j+\half,k)} &= \frac{1}{\dy}\LRp{ \overline{T}_{(i,j+1,k)} - \overline{T}_{(i,j,k)} }, \\
  \widehat{\dd{T}{z}}_{(i,j,k+\half)} &= \frac{1}{\dz}\LRp{ \overline{T}_{(i,j,k+1)} - \overline{T}_{(i,j,k)} }, 
\end{align*}
and tangential gradients by averaging two adjacent cell-centered gradients,
\begin{align*}
  \widehat{\dd{T}{y}}_{(i+\half,j,k)} &= \half\LRp{ \overline{\dd{T}{y}}_{(i,j,k)} + \overline{\dd{T}{y}}_{(i+1,j,k)} },&
  \widehat{\dd{T}{z}}_{(i+\half,j,k)} &= \half\LRp{ \overline{\dd{T}{z}}_{(i,j,k)} + \overline{\dd{T}{z}}_{(i+1,j,k)} }, 
  \\
  \widehat{\dd{T}{x}}_{(i,j+\half,k)} &= \half\LRp{ \overline{\dd{T}{x}}_{(i,j,k)} + \overline{\dd{T}{x}}_{(i,j+1,k)} },&
  \widehat{\dd{T}{z}}_{(i,j+\half,k)} &= \half\LRp{ \overline{\dd{T}{z}}_{(i,j,k)} + \overline{\dd{T}{z}}_{(i,j+1,k)} }, 
  \\
  \widehat{\dd{T}{x}}_{(i,j,k+\half)} &= \half\LRp{ \overline{\dd{T}{x}}_{(i,j,k)} + \overline{\dd{T}{x}}_{(i,j,k+1)} },&
  \widehat{\dd{T}{y}}_{(i,j,k+\half)} &= \half\LRp{ \overline{\dd{T}{y}}_{(i,j,k)} + \overline{\dd{T}{y}}_{(i,j,k+1)} }.
\end{align*}
Similarly, the common gradients of $\widehat{\Grad \ub}$ are computed. 

The rigid-lid interface condition in \cite{KANG2021113988}, 
\begin{subequations}
\eqnlab{rigid-lid-cond}
\begin{align*}
  \nb_m \cdot \ub_m  &= 0,\\
  \nb_m \cdot \LRp{\sigma_m \nb_m} &= 0,\\
  \tb \cdot \LRp{\sigma_1 \nb_1 - \sigma_2 \nb_2}  &= 0,\\
  \nb_1 \cdot \Pi_1 + \nb_2 \cdot \Pi_2 &=0,
\end{align*}
\end{subequations}
is simplified as 
\begin{align*}
  \hat{w}_1 &= \hat{w}_2 = 0,\\
  \hat{\sigma}_{1zz} &= \hat{\sigma}_{2zz} = 0,\\
  \mu_1 \dd{u_1}{z} &= \mu_2 \dd{u_2}{z} =: \hat{\sigma}_{xz} ,\\ 
  \mu_1 \dd{v_1}{z} &= \mu_2 \dd{v_2}{z} =: \hat{\sigma}_{yz} ,\\ 
   - \kappa_1 \dd{T_1}{z} &= - \kappa_2 \dd{T_2}{z} =: \hat{\Pi}_z,
\end{align*}
with $\nb_1 = (0,0,1)^T$, $\nb_2=(0,0,-1)^T$, $\tb=(1,0,0)^T$, and $\tb=(0,1,0)^T$. 
According to previous work in \cite{liu1979bulk,smith1988coefficients,fairall1996bulk,bao2000numerical}, 
 we can associate temperature and velocity with the heat and the horizontal momentum fluxes: 
\begin{align}
  \eqnlab{bulk-flux}
  \hat{\sigma}_{xz} := b_u (u_2 - u_1) \text{, }
  \hat{\sigma}_{yz} := b_v (v_2 - v_1) \text{, and }
  \hat{\Pi}_{z} := - b_T (T_2 - T_1) 
\end{align}
by introducing the bulk coefficients $b_u$, $b_v$, and $b_T$. Based on the finite difference (FD) approximation of heat and horizontal momentum fluxes, the linear bulk coefficients (with constant $\mu_1$,$\mu_2$,$\kappa_1$ and $\kappa_2$) can be obtained by \begin{align*}
  b_u = b_v = \frac{2 \mu_1 \mu_2}{\dz_2\mu_1 + \dz_1\mu_2}\text{, and }
  b_T = \frac{2 \kappa_1 \kappa_2}{\dz_2\kappa_1 + \dz_1\kappa_2}. 
\end{align*}

Once we have computed the heat and the momentum fluxes at the interface by using \eqnref{bulk-flux}, 
we estimate the isothermal wall boundary states of $u_w$, $v_w$, and $T_w$ for $\Omega_1$ and $\Omega_2$, respectively,
\begin{align*}
  u_{w_1} = u_1 + \frac{\hat{\sigma}_{xz} \dz_1}{2\mu_1 }\text{, } 
  v_{w_1} = v_1 + \frac{\hat{\sigma}_{yz} \dz_1}{2\mu_1 }\text{, and } T_{w_1} = T_1 - \frac{\hat{\Pi}_z \dz_1}{2\kappa_1 },\\
  u_{w_2} = u_2 - \frac{\hat{\sigma}_{xz} \dz_2}{2\mu_2 }\text{, }
  v_{w_2} = v_2 - \frac{\hat{\sigma}_{yz} \dz_2}{2\mu_2 } \text{, and } T_{w_2} = T_2 + \frac{\hat{\Pi}_z \dz_2}{2\kappa_2 }.
\end{align*}
 
\section*{References}

\bibliography{main}

\begin{thebibliography}{10}
\expandafter\ifx\csname url\endcsname\relax
  \def\url#1{\texttt{#1}}\fi
\expandafter\ifx\csname urlprefix\endcsname\relax\def\urlprefix{URL }\fi
\expandafter\ifx\csname href\endcsname\relax
  \def\href#1#2{#2} \def\path#1{#1}\fi

\bibitem{jacob2005m}
R.~Jacob, J.~Larson, E.~Ong, M$\times${N} communication and parallel
  interpolation in {C}ommunity {C}limate {S}ystem {M}odel version 3 using the
  {M}odel {C}oupling {T}oolkit, The International Journal of High Performance
  Computing Applications 19~(3) (2005) 293--307.

\bibitem{craig2012new}
A.~P. Craig, M.~Vertenstein, R.~Jacob, A new flexible coupler for earth system
  modeling developed for {CCSM4} and {CESM1}, The International Journal of High
  Performance Computing Applications 26~(1) (2012) 31--42.

\bibitem{golaz2019doe}
J.-C. Golaz, P.~M. Caldwell, L.~P. Van~Roekel, M.~R. Petersen, Q.~Tang, J.~D.
  Wolfe, G.~Abeshu, V.~Anantharaj, X.~S. Asay-Davis, D.~C. Bader, et~al., The
  {DOE} {E3SM} coupled model version 1: Overview and evaluation at standard
  resolution, Journal of Advances in Modeling Earth Systems 11~(7) (2019)
  2089--2129.

\bibitem{KANG2021113988}
S.~Kang, E.~M. Constantinescu, H.~Zhang, R.~L. Jacob, Mass-conserving
  implicit-explicit methods for coupled compressible {N}avier--{S}tokes
  equations, Computer Methods in Applied Mechanics and Engineering 384 (2021)
  113988.

\bibitem{constantinescu2007multirate}
E.~M. Constantinescu, A.~Sandu, Multirate timestepping methods for hyperbolic
  conservation laws, Journal of Scientific Computing 33~(3) (2007) 239--278.

\bibitem{sandu2019class}
A.~Sandu, A class of multirate infinitesimal {GARK} methods, SIAM Journal on
  Numerical Analysis 57~(5) (2019) 2300--2327.

\bibitem{roberts2020coupled}
S.~Roberts, A.~Sarshar, A.~Sandu, Coupled multirate infinitesimal {GARK}
  schemes for stiff systems with multiple scales, SIAM J. Sci. Comput 42~(3)
  (2020) A1609.

\bibitem{hachtel2021multirate}
C.~Hachtel, A.~Bartel, M.~G{\"u}nther, A.~Sandu, Multirate implicit {Euler}
  schemes for a class of differential--algebraic equations of index-1, Journal
  of Computational and Applied Mathematics 387 (2021) 112499.

\bibitem{gunther2021multirate}
M.~G{\"u}nther, A.~Sandu, Multirate linearly-implicit {GARK} schemes, BIT
  Numerical Mathematics 62 (2022) 869--901.

\bibitem{abdulle2022explicit}
A.~Abdulle, M.~Grote, G.~Rosilho~de Souza, Explicit stabilized multirate method
  for stiff differential equations, Mathematics of Computation 91~(338) (2022)
  2681--2714.

\bibitem{SkamarockKlemp08}
W.~C. Skamarock, J.~B. Klemp, A time-split nonhydrostatic atmospheric model for
  weather research and forecasting applications, Journal of Computational
  Physics 227~(7) (2008) 3465--3485.

\bibitem{seny2013multirate}
B.~Seny, J.~Lambrechts, R.~Comblen, V.~Legat, J.-F. Remacle, Multirate time
  stepping for accelerating explicit discontinuous {G}alerkin computations with
  application to geophysical flows, International Journal for Numerical Methods
  in Fluids 71~(1) (2013) 41--64.

\bibitem{schlegel2012implementation}
M.~Schlegel, O.~Knoth, M.~Arnold, R.~Wolke, Implementation of multirate time
  integration methods for air pollution modelling, Geoscientific Model
  Development 5~(6) (2012) 1395--1405.

\bibitem{lohner1984use}
R.~L{\"o}hner, K.~Morgan, O.~Zienkiewicz, The use of domain splitting with an
  explicit hyperbolic solver, Computer Methods in Applied Mechanics and
  Engineering 45~(1-3) (1984) 313--329.

\bibitem{kirby2003convergence}
R.~Kirby, On the convergence of high resolution methods with multiple time
  scales for hyperbolic conservation laws, Mathematics of Computation 72~(243)
  (2003) 1239--1250.

\bibitem{wensch2009multirate}
J.~Wensch, O.~Knoth, A.~Galant, Multirate infinitesimal step methods for
  atmospheric flow simulation, BIT Numerical Mathematics 49~(2) (2009)
  449--473.

\bibitem{mikida2019multi}
C.~Mikida, A.~Kl{\"o}ckner, D.~Bodony, Multi-rate time integration on overset
  meshes, Journal of Computational Physics 396 (2019) 325--346.

\bibitem{seny2014efficient}
B.~Seny, J.~Lambrechts, T.~Toulorge, V.~Legat, J.-F. Remacle, An efficient
  parallel implementation of explicit multirate {Runge--Kutta} schemes for
  discontinuous {G}alerkin computations, Journal of Computational Physics 256
  (2014) 135--160.

\bibitem{kopriva2002computation}
D.~A. Kopriva, S.~L. Woodruff, M.~Y. Hussaini, Computation of electromagnetic
  scattering with a non-conforming discontinuous spectral element method,
  International journal for numerical methods in engineering 53~(1) (2002)
  105--122.

\bibitem{bui2012analysis}
T.~Bui-Thanh, O.~Ghattas, Analysis of an hp-nonconforming discontinuous
  {G}alerkin spectral element method for wave propagation, SIAM Journal on
  Numerical Analysis 50~(3) (2012) 1801--1826.

\bibitem{kopriva2019free}
D.~A. Kopriva, F.~J. Hindenlang, T.~Bolemann, G.~J. Gassner, Free-stream
  preservation for curved geometrically non-conforming discontinuous {G}alerkin
  spectral elements, Journal of Scientific Computing 79 (2019) 1389--1408.

\bibitem{trask2020compatible}
N.~Trask, P.~Kuberry, Compatible meshfree discretization of surface pdes,
  Computational Particle Mechanics 7~(2) (2020) 271--277.

\bibitem{mahadevan2022metrics}
V.~S. Mahadevan, J.~E. Guerra, X.~Jiao, P.~Kuberry, Y.~Li, P.~Ullrich,
  D.~Marsico, R.~Jacob, P.~Bochev, P.~Jones, Metrics for intercomparison of
  remapping algorithms (mira) protocol applied to earth system models,
  Geoscientific Model Development 15~(17) (2022) 6601--6635.

\bibitem{gmd-16-1537-2023}
D.~H. Marsico, P.~A. Ullrich, Strategies for conservative and non-conservative
  monotone remapping on the sphere, Geoscientific Model Development 16~(5)
  (2023) 1537--1551.

\bibitem{toro2013riemann}
E.~F. Toro, Riemann solvers and numerical methods for fluid dynamics: A
  practical introduction, Springer Science \& Business Media, 2013.

\bibitem{nishikawa2011two}
H.~Nishikawa, Two ways to extend diffusion schemes to {N}avier--{S}tokes
  schemes: Gradient formula or upwind flux, in: 20th AIAA Computational Fluid
  Dynamics Conference, 2011, p. 3044.

\bibitem{AbhyankarEtAl2018}
S.~Abhyankar, J.~Brown, E.~M. Constantinescu, D.~Ghosh, B.~F. Smith, H.~Zhang,
  {PETSc/TS}: {A} modern scalable {ODE/DAE} solver library, arXiv e-preprints,
  1806.01437.

\bibitem{petsc-user-ref}
S.~Balay, S.~Abhyankar, M.~F. Adams, S.~Benson, J.~Brown, P.~Brune,
  K.~Buschelman, E.~Constantinescu, L.~Dalcin, A.~Dener, V.~Eijkhout, W.~D.
  Gropp, V.~Hapla, T.~Isaac, P.~Jolivet, D.~Karpeev, D.~Kaushik, M.~G. Knepley,
  F.~Kong, S.~Kruger, D.~A. May, L.~C. McInnes, R.~T. Mills, L.~Mitchell,
  T.~Munson, J.~E. Roman, K.~Rupp, P.~Sanan, J.~Sarich, B.~F. Smith,
  S.~Zampini, H.~Zhang, H.~Zhang, J.~Zhang, {PETSc/TAO} users manual, Tech.
  Rep. ANL-21/39 - Revision 3.16, Argonne National Laboratory (2021).

\bibitem{hindmarsh2005sundials}
A.~C. Hindmarsh, P.~N. Brown, K.~E. Grant, S.~L. Lee, R.~Serban, D.~E.
  Shumaker, C.~S. Woodward, {SUNDIALS}: Suite of nonlinear and
  differential/algebraic equation solvers, ACM Transactions on Mathematical
  Software (TOMS) 31~(3) (2005) 363--396.

\bibitem{gardner2022enabling}
D.~J. Gardner, D.~R. Reynolds, C.~S. Woodward, C.~J. Balos, Enabling new
  flexibility in the {SUNDIALS} suite of nonlinear and differential/algebraic
  equation solvers, ACM Transactions on Mathematical Software (TOMS) 48~(3)
  (2022) 1--24.

\bibitem{drazin2004hydrodynamic}
P.~G. Drazin, W.~H. Reid, Hydrodynamic stability, Cambridge University Press,
  2004.

\bibitem{carpenter1994fourth}
M.~H. Carpenter, C.~A. Kennedy, Fourth-order 2{N}-storage {R}unge--{K}utta
  schemes, Tech. Rep. {NASA} {TM} 109112, NASA Langley Research Center (1994).

\bibitem{syrakos2017critical}
A.~Syrakos, S.~Varchanis, Y.~Dimakopoulos, A.~Goulas, J.~Tsamopoulos, A
  critical analysis of some popular methods for the discretisation of the
  gradient operator in finite volume methods, Physics of Fluids 29~(12) (2017)
  127103.

\bibitem{liu1979bulk}
W.~T. Liu, K.~B. Katsaros, J.~A. Businger, Bulk parameterization of air--sea
  exchanges of heat and water vapor including the molecular constraints at the
  interface, Journal of the Atmospheric Sciences 36~(9) (1979) 1722--1735.

\bibitem{smith1988coefficients}
S.~D. Smith, Coefficients for sea surface wind stress, heat flux, and wind
  profiles as a function of wind speed and temperature, Journal of Geophysical
  Research: Oceans 93~(C12) (1988) 15467--15472.

\bibitem{fairall1996bulk}
C.~W. Fairall, E.~F. Bradley, D.~P. Rogers, J.~B. Edson, G.~S. Young, Bulk
  parameterization of air--sea fluxes for tropical ocean-global atmosphere
  coupled-ocean atmosphere response experiment, Journal of Geophysical
  Research: Oceans 101~(C2) (1996) 3747--3764.

\bibitem{bao2000numerical}
J.~Bao, J.~Wilczak, J.~Choi, L.~Kantha, Numerical simulations of air--sea
  interaction under high wind conditions using a coupled model: A study of
  hurricane development, Monthly Weather Review 128~(7) (2000) 2190--2210.

\end{thebibliography}

 \begin{center}
	\scriptsize \framebox{\parbox{4in}{Government License (will be removed at publication):
			The submitted manuscript has been created by UChicago Argonne, LLC,
			Operator of Argonne National Laboratory (``Argonne").  Argonne, a
			U.S. Department of Energy Office of Science laboratory, is operated
			under Contract No. DE-AC02-06CH11357.  The U.S. Government retains for
			itself, and others acting on its behalf, a paid-up nonexclusive,
			irrevocable worldwide license in said article to reproduce, prepare
			derivative works, distribute copies to the public, and perform
			publicly and display publicly, by or on behalf of the Government. The Department of Energy will provide public access to these results of federally sponsored research in accordance with the DOE Public Access Plan. http://energy.gov/downloads/doe-public-access-plan.
}}
	\normalsize
\end{center}

\end{document}